\input amstex

\expandafter\ifx\csname mathdefs.tex\endcsname\relax
  \expandafter\gdef\csname mathdefs.tex\endcsname{}
\else \message{Hey!  Apparently you were trying to
  \string\input{mathdefs.tex} twice.   This does not make sense.} 
\errmessage{Please edit your file (probably \jobname.tex) and remove
any duplicate ``\string\input'' lines}\endinput\fi




\catcode`\X=12\catcode`\@=11

\def\n@wcount{\alloc@0\count\countdef\insc@unt}
\def\n@wwrite{\alloc@7\write\chardef\sixt@@n}
\def\n@wread{\alloc@6\read\chardef\sixt@@n}
\def\r@s@t{\relax}\def\v@idline{\par}\def\@mputate#1/{#1}
\def\l@c@l#1X{\firstpart.#1}\def\gl@b@l#1X{#1}\def\t@d@l#1X{{}}

\def\crossrefs#1{\ifx\all#1\let\tr@ce=\all\else\def\tr@ce{#1,}\fi
   \n@wwrite\cit@tionsout\openout\cit@tionsout=\jobname.cit 
   \write\cit@tionsout{\tr@ce}\expandafter\setfl@gs\tr@ce,}
\def\setfl@gs#1,{\def\@{#1}\ifx\@\empty\let\next=\relax
   \else\let\next=\setfl@gs\expandafter\xdef
   \csname#1tr@cetrue\endcsname{}\fi\next}
\def\m@ketag#1#2{\expandafter\n@wcount\csname#2tagno\endcsname
     \csname#2tagno\endcsname=0\let\tail=\all\xdef\all{\tail#2,}
   \ifx#1\l@c@l\let\tail=\r@s@t\xdef\r@s@t{\csname#2tagno\endcsname=0\tail}\fi
   \expandafter\gdef\csname#2cite\endcsname##1{\expandafter
     \ifx\csname#2tag##1\endcsname\relax?\else\csname#2tag##1\endcsname\fi
     \expandafter\ifx\csname#2tr@cetrue\endcsname\relax\else
     \write\cit@tionsout{#2tag ##1 cited on page \folio.}\fi}
   \expandafter\gdef\csname#2page\endcsname##1{\expandafter
     \ifx\csname#2page##1\endcsname\relax?\else\csname#2page##1\endcsname\fi
     \expandafter\ifx\csname#2tr@cetrue\endcsname\relax\else
     \write\cit@tionsout{#2tag ##1 cited on page \folio.}\fi}
   \expandafter\gdef\csname#2tag\endcsname##1{\expandafter
      \ifx\csname#2check##1\endcsname\relax
      \expandafter\xdef\csname#2check##1\endcsname{}%
      \else\immediate\write16{Warning: #2tag ##1 used more than once.}\fi
      \multit@g{#1}{#2}##1/X%
      \write\t@gsout{#2tag ##1 assigned number \csname#2tag##1\endcsname\space
      on page \number\count0.}%
   \csname#2tag##1\endcsname}}

\def\multit@g#1#2#3/#4X{\def\t@mp{#4}\ifx\t@mp\empty%
      \global\advance\csname#2tagno\endcsname by 1 
      \expandafter\xdef\csname#2tag#3\endcsname
      {#1\number\csname#2tagno\endcsnameX}%
   \else\expandafter\ifx\csname#2last#3\endcsname\relax
      \expandafter\n@wcount\csname#2last#3\endcsname
      \global\advance\csname#2tagno\endcsname by 1 
      \expandafter\xdef\csname#2tag#3\endcsname
      {#1\number\csname#2tagno\endcsnameX}
      \write\t@gsout{#2tag #3 assigned number \csname#2tag#3\endcsname\space
      on page \number\count0.}\fi
   \global\advance\csname#2last#3\endcsname by 1
   \def\t@mp{\expandafter\xdef\csname#2tag#3/}%
   \expandafter\t@mp\@mputate#4\endcsname
   {\csname#2tag#3\endcsname\lastpart{\csname#2last#3\endcsname}}\fi}
\def\t@gs#1{\def\all{}\m@ketag#1e\m@ketag#1s\m@ketag\t@d@l p
\let\realscite\scite
\let\realstag\stag
   \m@ketag\gl@b@l r \n@wread\t@gsin
   \openin\t@gsin=\jobname.tgs \re@der \closein\t@gsin
   \n@wwrite\t@gsout\openout\t@gsout=\jobname.tgs }
\outer\def\localtags{\t@gs\l@c@l}
\outer\def\globaltags{\t@gs\gl@b@l}
\outer\def\newlocaltag#1{\m@ketag\l@c@l{#1}}
\outer\def\newglobaltag#1{\m@ketag\gl@b@l{#1}}

\newif\ifpr@ 
\def\m@kecs #1tag #2 assigned number #3 on page #4.%
   {\expandafter\gdef\csname#1tag#2\endcsname{#3}
   \expandafter\gdef\csname#1page#2\endcsname{#4}
   \ifpr@\expandafter\xdef\csname#1check#2\endcsname{}\fi}
\def\re@der{\ifeof\t@gsin\let\next=\relax\else
   \read\t@gsin to\t@gline\ifx\t@gline\v@idline\else
   \expandafter\m@kecs \t@gline\fi\let \next=\re@der\fi\next}
\def\pretags#1{\pr@true\pret@gs#1,,}
\def\pret@gs#1,{\def\@{#1}\ifx\@\empty\let\n@xtfile=\relax
   \else\let\n@xtfile=\pret@gs \openin\t@gsin=#1.tgs \message{#1} \re@der 
   \closein\t@gsin\fi \n@xtfile}

\newcount\sectno\sectno=0\newcount\subsectno\subsectno=0
\newif\ifultr@local \def\ultralocal{\ultr@localtrue}
\def\firstpart{\number\sectno}
\def\lastpart#1{\ifcase#1 \or a\or b\or c\or d\or e\or f\or g\or h\or 
   i\or k\or l\or m\or n\or o\or p\or q\or r\or s\or t\or u\or v\or w\or 
   x\or y\or z \fi}

\def\resetall{\global\advance\sectno by 1\subsectno=0
   \gdef\firstpart{\number\sectno}\r@s@t}
\def\resetsub{\global\advance\subsectno by 1
   \gdef\firstpart{\number\sectno.\number\subsectno}\r@s@t}
\def\newsection#1\par{\resetall\vskip0pt plus.3\vsize\penalty-250
   \vskip0pt plus-.3\vsize\bigskip\bigskip
   \message{#1}\leftline{\bf#1}\nobreak\bigskip}
\def\subsection#1\par{\ifultr@local\resetsub\fi
   \vskip0pt plus.2\vsize\penalty-250\vskip0pt plus-.2\vsize
   \bigskip\smallskip\message{#1}\leftline{\bf#1}\nobreak\medskip}


\newdimen\marginshift

\newdimen\margindelta
\newdimen\marginmax
\newdimen\marginmin

\def\margininit{       
\marginmax=3 true cm                  
				      
\margindelta=0.1 true cm              
\marginmin=0.1true cm                 
\marginshift=\marginmin
}    

\def\t@gsjj#1,{\def\@{#1}\ifx\@\empty\let\next=\relax\else\let\next=\t@gsjj
   \def\@@{p}\ifx\@\@@\else
   \expandafter\gdef\csname#1cite\endcsname##1{\citejj{##1}}
   \expandafter\gdef\csname#1page\endcsname##1{?}
   \expandafter\gdef\csname#1tag\endcsname##1{\tagjj{##1}}\fi\fi\next}
\newif\ifshowstuffinmargin
\showstuffinmarginfalse
\def\jjtags{\ifx\shlhetal\relax 
  \else
\ifx\shlhetal\undefinedcontrolseq
\else
\showstuffinmargintrue
\ifx\all\relax\else\expandafter\t@gsjj\all,\fi\fi \fi
}

\def\tagjj#1{\realstag{#1}\oldmginpar{\zeigen{#1}}}
\def\citejj#1{\rechnen{#1}\mginpar{\zeigen{#1}}}     

\def\rechnen#1{\expandafter\ifx\csname stag#1\endcsname\relax ??\else
                           \csname stag#1\endcsname\fi}

\newdimen\theight

\def\marginfont{\sevenrm}

\def\trymarginbox#1{\setbox0=\hbox{\marginfont\hskip\marginshift #1}%
		\global\marginshift\wd0 
		\global\advance\marginshift\margindelta}

\def \oldmginpar#1{%
\ifvmode\setbox0\hbox to \hsize{\hfill\rlap{\marginfont\quad#1}}%
\ht0 0cm
\dp0 0cm
\box0\vskip-\baselineskip
\else 
             \vadjust{\trymarginbox{#1}%
		\ifdim\marginshift>\marginmax \global\marginshift\marginmin
			\trymarginbox{#1}%
                \fi
             \theight=\ht0
             \advance\theight by \dp0    \advance\theight by \lineskip
             \kern -\theight \vbox to \theight{\rightline{\rlap{\box0}}%
\vss}}\fi}

\newdimen\upordown
\global\upordown=8pt
\font\tinyfont=cmtt8 
\def\mginpar#1{\smash{\hbox to 0cm{\kern-10pt\raise7pt\hbox{\tinyfont #1}\hss}}}
\def\mginpar#1{{\hbox to 0cm{\kern-10pt\raise\upordown\hbox{\tinyfont #1}\hss}}\global\upordown-\upordown}


\def\t@gsoff#1,{\def\@{#1}\ifx\@\empty\let\next=\relax\else\let\next=\t@gsoff
   \def\@@{p}\ifx\@\@@\else
   \expandafter\gdef\csname#1cite\endcsname##1{\zeigen{##1}}
   \expandafter\gdef\csname#1page\endcsname##1{?}
   \expandafter\gdef\csname#1tag\endcsname##1{\zeigen{##1}}\fi\fi\next}
\def\verbatimtags{\showstuffinmarginfalse
\ifx\all\relax\else\expandafter\t@gsoff\all,\fi}
\def\zeigen#1{\hbox{$\scriptstyle\langle$}#1\hbox{$\scriptstyle\rangle$}}


\def\margintag#1{\ifshowstuffinmargin\oldmginpar{\zeigen{#1}}\fi}

\def\(#1){\edef\dot@g{\ifmmode\ifinner(\hbox{\noexpand\etag{#1}})
   \else\noexpand\eqno(\hbox{\noexpand\etag{#1}})\fi
   \else(\noexpand\ecite{#1})\fi}\dot@g}

\newif\ifbr@ck
\def\eat#1{}
\def\[#1]{\br@cktrue[\br@cket#1'X]}
\def\br@cket#1'#2X{\def\temp{#2}\ifx\temp\empty\let\next\eat
   \else\let\next\br@cket\fi
   \ifbr@ck\br@ckfalse\br@ck@t#1,X\else\br@cktrue#1\fi\next#2X}
\def\br@ck@t#1,#2X{\def\temp{#2}\ifx\temp\empty\let\neext\eat
   \else\let\neext\br@ck@t\def\temp{,}\fi
   \def\teemp{#1}\ifx\teemp\empty\else\rcite{#1}\fi\temp\neext#2X}
\def\resetbr@cket{\gdef\[##1]{[\rtag{##1}]}}
\def\references{\resetbr@cket\newsection References\par}

\newtoks\symb@ls\newtoks\s@mb@ls\newtoks\p@gelist\n@wcount\ftn@mber
    \ftn@mber=1\newif\ifftn@mbers\ftn@mbersfalse\newif\ifbyp@ge\byp@gefalse
\def\defm@rk{\ifftn@mbers\n@mberm@rk\else\symb@lm@rk\fi}
\def\n@mberm@rk{\xdef\m@rk{{\the\ftn@mber}}%
    \global\advance\ftn@mber by 1 }
\def\rot@te#1{\let\temp=#1\global#1=\expandafter\r@t@te\the\temp,X}
\def\r@t@te#1,#2X{{#2#1}\xdef\m@rk{{#1}}}
\def\b@@st#1{{$^{#1}$}}\def\str@p#1{#1}
\def\symb@lm@rk{\ifbyp@ge\rot@te\p@gelist\ifnum\expandafter\str@p\m@rk=1 
    \s@mb@ls=\symb@ls\fi\write\f@nsout{\number\count0}\fi \rot@te\s@mb@ls}
\def\byp@ge{\byp@getrue\n@wwrite\f@nsin\openin\f@nsin=\jobname.fns 
    \n@wcount\currentp@ge\currentp@ge=0\p@gelist={0}
    \re@dfns\closein\f@nsin\rot@te\p@gelist
    \n@wread\f@nsout\openout\f@nsout=\jobname.fns }
\def\m@kelist#1X#2{{#1,#2}}
\def\re@dfns{\ifeof\f@nsin\let\next=\relax\else\read\f@nsin to \f@nline
    \ifx\f@nline\v@idline\else\let\t@mplist=\p@gelist
    \ifnum\currentp@ge=\f@nline
    \global\p@gelist=\expandafter\m@kelist\the\t@mplistX0
    \else\currentp@ge=\f@nline
    \global\p@gelist=\expandafter\m@kelist\the\t@mplistX1\fi\fi
    \let\next=\re@dfns\fi\next}
\def\symbols#1{\symb@ls={#1}\s@mb@ls=\symb@ls} 
\def\bigsymbol{\textstyle}
\symbols{\bigsymbol\ast,\dagger,\ddagger,\sharp,\flat,\natural,\star}
\def\ftnumbers{\ftn@mberstrue} \def\ftsymbols{\ftn@mbersfalse}
\def\paginal{\byp@ge} \def\resetftnumbers{\ftn@mber=1}
\def\ftnote#1{\defm@rk\expandafter\expandafter\expandafter\footnote
    \expandafter\b@@st\m@rk{#1}}

\long\def\jump#1\endjump{}
\def\ssum{\mathop{\lower .1em\hbox{$\textstyle\Sigma$}}\nolimits}

\def\qed{\nobreak\kern 1em \vrule height .5em width .5em depth 0em}
\def\newneq{\hbox{\rlap{\hbox to 1\wd9{\hss$=$\hss}}\raise .1em 
   \hbox to 1\wd9{\hss$\scriptscriptstyle/$\hss}}}
\def\subsetne{\setbox9 = \hbox{$\subset$}\mathrel{\hbox{\rlap
   {\lower .4em \newneq}\raise .13em \hbox{$\subset$}}}}
\def\supsetne{\setbox9 = \hbox{$\subset$}\mathrel{\hbox{\rlap
   {\lower .4em \newneq}\raise .13em \hbox{$\supset$}}}}

\def\vbar{\mathchoice{\vrule height6.3ptdepth-.5ptwidth.8pt\kern-.8pt}
   {\vrule height6.3ptdepth-.5ptwidth.8pt\kern-.8pt}
   {\vrule height4.1ptdepth-.35ptwidth.6pt\kern-.6pt}
   {\vrule height3.1ptdepth-.25ptwidth.5pt\kern-.5pt}}
\def\f@dge{\mathchoice{}{}{\mkern.5mu}{\mkern.8mu}}
\def\b@c#1#2{{\rm \mkern#2mu\vbar\mkern-#2mu#1}}
\def\b@b#1{{\rm I\mkern-3.5mu #1}}
\def\b@a#1#2{{\rm #1\mkern-#2mu\f@dge #1}}
\def\bb#1{{\count4=`#1 \advance\count4by-64 \ifcase\count4\or\b@a A{11.5}\or
   \b@b B\or\b@c C{5}\or\b@b D\or\b@b E\or\b@b F \or\b@c G{5}\or\b@b H\or
   \b@b I\or\b@c J{3}\or\b@b K\or\b@b L \or\b@b M\or\b@b N\or\b@c O{5} \or
   \b@b P\or\b@c Q{5}\or\b@b R\or\b@a S{8}\or\b@a T{10.5}\or\b@c U{5}\or
   \b@a V{12}\or\b@a W{16.5}\or\b@a X{11}\or\b@a Y{11.7}\or\b@a Z{7.5}\fi}}

\catcode`\X=11 \catcode`\@=12




\let\thischap\jobname

\def\partof#1{\csname returnthe#1part\endcsname}
\def\CHAPOF#1{\csname returnthe#1chap\endcsname}

\def\chapof#1{\CHAPOF{#1}}

\def\setchapter#1,#2,#3;{%
  \expandafter\def\csname returnthe#1part\endcsname{#2}%
  \expandafter\def\csname returnthe#1chap\endcsname{#3}%
}

\def\setprevious#1 #2 {%
  \expandafter\def\csname set#1page\endcsname{\input page-#2}
}


 \setchapter  E53,B,N;       \setprevious E53 null
 \setchapter  300z,B,A;       \setprevious 300z E53
 \setchapter  88r,B,I;       \setprevious 88r 300z
 \setchapter  600,B,II;       \setprevious  600 88r
 \setchapter  705,B,III;       \setprevious   705 600
 \setchapter  734,B,IV;        \setprevious   734 705
 \setchapter  300x,B,;      \setprevious   300x 734

 \setchapter 300a,A,V.A;      \setprevious 300a 88r
 \setchapter 300b,A,V.B;       \setprevious 300b 300a
 \setchapter 300c,A,V.C;       \setprevious 300c 300b
 \setchapter 300d,A,V.D;       \setprevious 300d 300c
 \setchapter 300e,A,V.E;       \setprevious 300e 300d
 \setchapter 300f,A,V.F;       \setprevious 300f 300e
 \setchapter 300g,A,V.G;       \setprevious 300g 300f

  \setchapter  E46,B,VI;      \setprevious    E46 734
  \setchapter  838,B,VII;      \setprevious   838 E46

\newwrite\pageout
\def\rememberpagenumber{\let\setpage\relax
\openout\pageout page-\jobname  \relax \write\pageout{\setpage\the\pageno.}}

\def\recallpagenumber{\csname set\jobname page\endcsname
\def\headmark##1{\rightheadtext{\chapof{\jobname}.##1}}\WRITETOC}
\def\setupchapter#1{\leftheadtext{\chapof{\jobname}. #1}}

\def\setpage#1.{\pageno=#1\relax\advance\pageno1\relax}

\def\cprefix#1{
\edef\theotherpart{\partof{#1}}\edef\theotherchap{\chapof{#1}}%
\ifx\theotherpart\thispart
   \ifx\theotherchap\thischap 
    \else 
     \theotherchap%
    \fi
   \else 
     \theotherchap\fi}

\def\sectioncite[#1]#2{%
     \cprefix{#2}#1}

\edef\thispart{\partof{\thischap}}
\edef\thischap{\chapof{\thischap}}

\def\lastpage of '#1' is #2.{\expandafter\def\csname lastpage#1\endcsname{#2}}

\def\yCITE[#1]#2{\cprefix{#2}.\scite{#2-#1}}

\newwrite\writetoc
\immediate\openout\writetoc \jobname.toc
\def\addcontents#1{\def\WRITETOC{\immediate\write\writetoc{\noexpand\tocentry{\chapof{\jobname}}{#1}{\number\pageno}}}}



\def\spuriousreset{}


\expandafter\ifx\csname citeadd.tex\endcsname\relax
\expandafter\gdef\csname citeadd.tex\endcsname{}
\else \message{Hey!  Apparently you were trying to
\string\input{citeadd.tex} twice.   This does not make sense.} 
\errmessage{Please edit your file (probably \jobname.tex) and remove
any duplicate ``\string\input'' lines}\endinput\fi

\sectno=-1   
\localtags
\jjtags
\NoBlackBoxes
 \newbox\noforkbox \newdimen\forklinewidth
\forklinewidth=0.3pt   
\setbox0\hbox{$\textstyle\bigcup$}
\setbox1\hbox to \wd0{\hfil\vrule width \forklinewidth depth \dp0
                        height \ht0 \hfil}
\wd1=0 cm
\setbox\noforkbox\hbox{\box1\box0\relax}
\def\unionstick{\mathop{\copy\noforkbox}\limits}
\def\nonfork#1#2_#3{#1\unionstick_{\textstyle #3}#2}
\def\nonforkin#1#2_#3^#4{#1\unionstick_{\textstyle #3}^{\textstyle #4}#2}     
%
\setbox0\hbox{$\textstyle\bigcup$}
\setbox1\hbox to \wd0{\hfil{\sl /\/}\hfil}
\setbox2\hbox to \wd0{\hfil\vrule height \ht0 depth \dp0 width
                                \forklinewidth\hfil}
\wd1=0cm
\wd2=0cm
\newbox\doesforkbox
\setbox\doesforkbox\hbox{\box1\box0\relax}
\def\nunionstick{\mathop{\copy\doesforkbox}\limits}

\def\fork#1#2_#3{#1\nunionstick_{\textstyle #3}#2}
\def\forkin#1#2_#3^#4{#1\nunionstick_{\textstyle #3}^{\textstyle #4}#2}     

\define\mr{\medskip\roster}
\define\sn{\smallskip\noindent}
\define\mn{\medskip\noindent}
\define\bn{\bigskip\noindent}
\define\ub{\underbar}
\define\wilog{\text{without loss of generality}}
\define\ermn{\endroster\medskip\noindent}

\define\rest{\restriction}
\define\dbcu{\dsize\bigcup}
\define \nl{\newline}
\magnification=\magstep 1
\documentstyle{amsppt}

{    
\catcode`@11

\ifx\alicetwothousandloaded@\relax
  \endinput\else\global\let\alicetwothousandloaded@\relax\fi

\gdef\subjclass{\let\savedef@\subjclass
 \def\subjclass##1\endsubjclass{\let\subjclass\savedef@
   \toks@{\def\usualspace{{\rm\enspace}}\eightpoint}%
   \toks@@{##1\unskip.}%
   \edef\thesubjclass@{\the\toks@
     \frills@{{\noexpand\rm2000 {\noexpand\it Mathematics Subject
       Classification}.\noexpand\enspace}}%
     \the\toks@@}}%
  \nofrillscheck\subjclass}
} 


\expandafter\ifx\csname alice2jlem.tex\endcsname\relax
  \expandafter\xdef\csname alice2jlem.tex\endcsname{\the\catcode`@}
\else \message{Hey!  Apparently you were trying to
\string\input{alice2jlem.tex}  twice.   This does not make sense.}
\errmessage{Please edit your file (probably \jobname.tex) and remove
any duplicate ``\string\input'' lines}\endinput\fi

\expandafter\ifx\csname bib4plain.tex\endcsname\relax
  \expandafter\gdef\csname bib4plain.tex\endcsname{}
\else \message{Hey!  Apparently you were trying to \string\input
  bib4plain.tex twice.   This does not make sense.}
\errmessage{Please edit your file (probably \jobname.tex) and remove
any duplicate ``\string\input'' lines}\endinput\fi

\def\renewcommand{\newcommand}	       
\edef\cite{\the\catcode`@}%
\catcode`@ = 11
\let\@oldatcatcode = \cite
\chardef\@letter = 11
\chardef\@other = 12
%
%
%
%
\def\@innerdef#1#2{\edef#1{\expandafter\noexpand\csname #2\endcsname}}%
%
%
\@innerdef\@innernewcount{newcount}%
\@innerdef\@innernewdimen{newdimen}%
\@innerdef\@innernewif{newif}%
\@innerdef\@innernewwrite{newwrite}%
%
%
%
\def\@gobble#1{}%
%
%
%
\ifx\inputlineno\@undefined
   \let\@linenumber = \empty 
\else
   \def\@linenumber{\the\inputlineno:\space}%
\fi
%
%
%
\def\@futurenonspacelet#1{\def\cs{#1}%
   \afterassignment\@stepone\let\@nexttoken=
}%
\begingroup 
\def\\{\global\let\@stoken= }%
\\ 
\endgroup
\def\@stepone{\expandafter\futurelet\cs\@steptwo}%
\def\@steptwo{\expandafter\ifx\cs\@stoken\let\@@next=\@stepthree
   \else\let\@@next=\@nexttoken\fi \@@next}%
\def\@stepthree{\afterassignment\@stepone\let\@@next= }%
%
%
%
\def\@getoptionalarg#1{%
   \let\@optionaltemp = #1%
   \let\@optionalnext = \relax
   \@futurenonspacelet\@optionalnext\@bracketcheck
}%
%
%
\def\@bracketcheck{%
   \ifx [\@optionalnext
      \expandafter\@@getoptionalarg
   \else
      \let\@optionalarg = \empty
      \expandafter\@optionaltemp
   \fi
}%
\def\@@getoptionalarg[#1]{%
   \def\@optionalarg{#1}%
   \@optionaltemp
}%
%
%
%
\def\@nnil{\@nil}%
\def\@fornoop#1\@@#2#3{}%
\def\@for#1:=#2\do#3{%
   \edef\@fortmp{#2}%
   \ifx\@fortmp\empty \else
      \expandafter\@forloop#2,\@nil,\@nil\@@#1{#3}%
   \fi
}%
\def\@forloop#1,#2,#3\@@#4#5{\def#4{#1}\ifx #4\@nnil \else
       #5\def#4{#2}\ifx #4\@nnil \else#5\@iforloop #3\@@#4{#5}\fi\fi
}%
\def\@iforloop#1,#2\@@#3#4{\def#3{#1}\ifx #3\@nnil
       \let\@nextwhile=\@fornoop \else
      #4\relax\let\@nextwhile=\@iforloop\fi\@nextwhile#2\@@#3{#4}%
}%
%
%
%
\@innernewif\if@fileexists
\def\@testfileexistence{\@getoptionalarg\@finishtestfileexistence}%
\def\@finishtestfileexistence#1{%
   \begingroup
      \def\extension{#1}%
      \immediate\openin0 =
         \ifx\@optionalarg\empty\jobname\else\@optionalarg\fi
         \ifx\extension\empty \else .#1\fi
         \space
      \ifeof 0
         \global\@fileexistsfalse
      \else
         \global\@fileexiststrue
      \fi
      \immediate\closein0
   \endgroup
}%
%
%
%
%
\def\bibliographystyle#1{%
   \@readauxfile
   \@writeaux{\string\bibstyle{#1}}%
}%
\let\bibstyle = \@gobble
%
%
\let\bblfilebasename = \jobname
\def\bibliography#1{%
   \@readauxfile
   \@writeaux{\string\bibdata{#1}}%
   \@testfileexistence[\bblfilebasename]{bbl}%
   \if@fileexists
      \nobreak
      \@readbblfile
   \fi
}%
\let\bibdata = \@gobble
%
%
\def\nocite#1{%
   \@readauxfile
   \@writeaux{\string\citation{#1}}%
}%
\@innernewif\if@notfirstcitation
%
%
\def\cite{\@getoptionalarg\@cite}%
%
%
\def\@cite#1{%
   \let\@citenotetext = \@optionalarg
   \printcitestart
   \nocite{#1}%
   \@notfirstcitationfalse
   \@for \@citation :=#1\do
   {%
      \expandafter\@onecitation\@citation\@@
   }%
   \ifx\empty\@citenotetext\else
      \printcitenote{\@citenotetext}%
   \fi
   \printcitefinish
}%
\newif\ifweareinprivate
\weareinprivatetrue
\ifx\shlhetal\undefinedcontrolseq\weareinprivatefalse\fi
\ifx\shlhetal\relax\weareinprivatefalse\fi
\def\@onecitation#1\@@{%
   \if@notfirstcitation
      \printbetweencitations
   \fi
   \expandafter \ifx \csname\@citelabel{#1}\endcsname \relax
      \if@citewarning
         \message{\@linenumber Undefined citation `#1'.}%
      \fi
     \ifweareinprivate
      \expandafter\gdef\csname\@citelabel{#1}\endcsname{%
\strut 
\vadjust{\vskip-\dp\strutbox
\vbox to 0pt{\vss\parindent0cm \leftskip=\hsize 
\advance\leftskip3mm
\advance\hsize 4cm\strut\openup-4pt 
\rightskip 0cm plus 1cm minus 0.5cm ?  #1 ?\strut}}
         {\tt
            \escapechar = -1
            \nobreak\hskip0pt\pfeilsw
            \expandafter\string\csname#1\endcsname
             \pfeilso
            \nobreak\hskip0pt
         }%
      }%
     \else  
      \expandafter\gdef\csname\@citelabel{#1}\endcsname{%
            {\tt\expandafter\string\csname#1\endcsname}
      }%
     \fi  
   \fi
   \csname\@citelabel{#1}\endcsname
   \@notfirstcitationtrue
}%
%
%
\def\@citelabel#1{b@#1}%
%
%
\def\@citedef#1#2{\expandafter\gdef\csname\@citelabel{#1}\endcsname{#2}}%
%
%
%
\def\@readbblfile{%
   \ifx\@itemnum\@undefined
      \@innernewcount\@itemnum
   \fi
   \begingroup
      \def\begin##1##2{%
         \setbox0 = \hbox{\biblabelcontents{##2}}%
         \biblabelwidth = \wd0
      }%
      \def\end##1{}
      %
      %
      \@itemnum = 0
      \def\bibitem{\@getoptionalarg\@bibitem}%
      \def\@bibitem{%
         \ifx\@optionalarg\empty
            \expandafter\@numberedbibitem
         \else
            \expandafter\@alphabibitem
         \fi
      }%
      \def\@alphabibitem##1{%
         \expandafter \xdef\csname\@citelabel{##1}\endcsname {\@optionalarg}%
         \ifx\biblabelprecontents\@undefined
            \let\biblabelprecontents = \relax
         \fi
         \ifx\biblabelpostcontents\@undefined
            \let\biblabelpostcontents = \hss
         \fi
         \@finishbibitem{##1}%
      }%
      \def\@numberedbibitem##1{%
         \advance\@itemnum by 1
         \expandafter \xdef\csname\@citelabel{##1}\endcsname{\number\@itemnum}%
         \ifx\biblabelprecontents\@undefined
            \let\biblabelprecontents = \hss
         \fi
         \ifx\biblabelpostcontents\@undefined
            \let\biblabelpostcontents = \relax
         \fi
         \@finishbibitem{##1}%
      }%
      \def\@finishbibitem##1{%
         \biblabelprint{\csname\@citelabel{##1}\endcsname}%
         \@writeaux{\string\@citedef{##1}{\csname\@citelabel{##1}\endcsname}}%
         \ignorespaces
      }%
      %
      %
      \let\em = \bblem
      \let\newblock = \bblnewblock
      \let\sc = \bblsc
      \frenchspacing
      \clubpenalty = 4000 \widowpenalty = 4000
      \tolerance = 10000 \hfuzz = .5pt
      \everypar = {\hangindent = \biblabelwidth
                      \advance\hangindent by \biblabelextraspace}%
      \bblrm
      \parskip = 1.5ex plus .5ex minus .5ex
      \biblabelextraspace = .5em
      \bblhook
      \input \bblfilebasename.bbl
   \endgroup
}%
%
%
\@innernewdimen\biblabelwidth
\@innernewdimen\biblabelextraspace
%
%
%
\def\biblabelprint#1{%
   \noindent
   \hbox to \biblabelwidth{%
      \biblabelprecontents
      \biblabelcontents{#1}%
      \biblabelpostcontents
   }%
   \kern\biblabelextraspace
}%
%
%
%
\def\biblabelcontents#1{{\bblrm [#1]}}%
%
%
\def\bblrm{\rm}%
%
%
\def\bblem{\it}%
%
%
\def\bblsc{\ifx\@scfont\@undefined
              \font\@scfont = cmcsc10
           \fi
           \@scfont
}%
%
%
\def\bblnewblock{\hskip .11em plus .33em minus .07em }%
%
%
\let\bblhook = \empty
%
%
%
\def\printcitestart{[}
\def\printcitefinish{]}
\def\printbetweencitations{, }
\def\printcitenote#1{, #1}
%
%
%
\let\citation = \@gobble
%
%
%
\@innernewcount\@numparams
%
%
\def\newcommand#1{%
   \def\@commandname{#1}%
   \@getoptionalarg\@continuenewcommand
}%
%
%
\def\@continuenewcommand{%
   \@numparams = \ifx\@optionalarg\empty 0\else\@optionalarg \fi \relax
   \@newcommand
}%
%
%
\def\@newcommand#1{%
   \def\@startdef{\expandafter\edef\@commandname}%
   \ifnum\@numparams=0
      \let\@paramdef = \empty
   \else
      \ifnum\@numparams>9
         \errmessage{\the\@numparams\space is too many parameters}%
      \else
         \ifnum\@numparams<0
            \errmessage{\the\@numparams\space is too few parameters}%
         \else
            \edef\@paramdef{%
               \ifcase\@numparams
                  \empty  No arguments.
               \or ####1%
               \or ####1####2%
               \or ####1####2####3%
               \or ####1####2####3####4%
               \or ####1####2####3####4####5%
               \or ####1####2####3####4####5####6%
               \or ####1####2####3####4####5####6####7%
               \or ####1####2####3####4####5####6####7####8%
               \or ####1####2####3####4####5####6####7####8####9%
               \fi
            }%
         \fi
      \fi
   \fi
   \expandafter\@startdef\@paramdef{#1}%
}%
%
%
%
%
\def\@readauxfile{%
   \if@auxfiledone \else 
      \global\@auxfiledonetrue
      \@testfileexistence{aux}%
      \if@fileexists
         \begingroup
            \endlinechar = -1
            \catcode`@ = 11
            \input \jobname.aux
         \endgroup
      \else
         \message{\@undefinedmessage}%
         \global\@citewarningfalse
      \fi
      \immediate\openout\@auxfile = \jobname.aux
   \fi
}%
%
%
\newif\if@auxfiledone
\ifx\noauxfile\@undefined \else \@auxfiledonetrue\fi
%
%
%
%
\@innernewwrite\@auxfile
\def\@writeaux#1{\ifx\noauxfile\@undefined \write\@auxfile{#1}\fi}%
%
%
%
\ifx\@undefinedmessage\@undefined
   \def\@undefinedmessage{No .aux file; I won't give you warnings about
                          undefined citations.}%
\fi
%
%
\@innernewif\if@citewarning
\ifx\noauxfile\@undefined \@citewarningtrue\fi
%
%
%
\catcode`@ = \@oldatcatcode

\def\pfeilso{\leavevmode
            \vrule width 1pt height9pt depth 0pt\relax
           \vrule width 1pt height8.7pt depth 0pt\relax
           \vrule width 1pt height8.3pt depth 0pt\relax
           \vrule width 1pt height8.0pt depth 0pt\relax
           \vrule width 1pt height7.7pt depth 0pt\relax
            \vrule width 1pt height7.3pt depth 0pt\relax
            \vrule width 1pt height7.0pt depth 0pt\relax
            \vrule width 1pt height6.7pt depth 0pt\relax
            \vrule width 1pt height6.3pt depth 0pt\relax
            \vrule width 1pt height6.0pt depth 0pt\relax
            \vrule width 1pt height5.7pt depth 0pt\relax
            \vrule width 1pt height5.3pt depth 0pt\relax
            \vrule width 1pt height5.0pt depth 0pt\relax
            \vrule width 1pt height4.7pt depth 0pt\relax
            \vrule width 1pt height4.3pt depth 0pt\relax
            \vrule width 1pt height4.0pt depth 0pt\relax
            \vrule width 1pt height3.7pt depth 0pt\relax
            \vrule width 1pt height3.3pt depth 0pt\relax
            \vrule width 1pt height3.0pt depth 0pt\relax
            \vrule width 1pt height2.7pt depth 0pt\relax
            \vrule width 1pt height2.3pt depth 0pt\relax
            \vrule width 1pt height2.0pt depth 0pt\relax
            \vrule width 1pt height1.7pt depth 0pt\relax
            \vrule width 1pt height1.3pt depth 0pt\relax
            \vrule width 1pt height1.0pt depth 0pt\relax
            \vrule width 1pt height0.7pt depth 0pt\relax
            \vrule width 1pt height0.3pt depth 0pt\relax}

\def\pfeilsw{ \leavevmode 
            \vrule width 1pt height0.3pt depth 0pt\relax
            \vrule width 1pt height0.7pt depth 0pt\relax
            \vrule width 1pt height1.0pt depth 0pt\relax
            \vrule width 1pt height1.3pt depth 0pt\relax
            \vrule width 1pt height1.7pt depth 0pt\relax
            \vrule width 1pt height2.0pt depth 0pt\relax
            \vrule width 1pt height2.3pt depth 0pt\relax
            \vrule width 1pt height2.7pt depth 0pt\relax
            \vrule width 1pt height3.0pt depth 0pt\relax
            \vrule width 1pt height3.3pt depth 0pt\relax
            \vrule width 1pt height3.7pt depth 0pt\relax
            \vrule width 1pt height4.0pt depth 0pt\relax
            \vrule width 1pt height4.3pt depth 0pt\relax
            \vrule width 1pt height4.7pt depth 0pt\relax
            \vrule width 1pt height5.0pt depth 0pt\relax
            \vrule width 1pt height5.3pt depth 0pt\relax
            \vrule width 1pt height5.7pt depth 0pt\relax
            \vrule width 1pt height6.0pt depth 0pt\relax
            \vrule width 1pt height6.3pt depth 0pt\relax
            \vrule width 1pt height6.7pt depth 0pt\relax
            \vrule width 1pt height7.0pt depth 0pt\relax
            \vrule width 1pt height7.3pt depth 0pt\relax
            \vrule width 1pt height7.7pt depth 0pt\relax
            \vrule width 1pt height8.0pt depth 0pt\relax
            \vrule width 1pt height8.3pt depth 0pt\relax
            \vrule width 1pt height8.7pt depth 0pt\relax
            \vrule width 1pt height9pt depth 0pt\relax
      }


\def\widestnumber#1#2{}

\def\citewarning#1{\ifx\shlhetal\relax 
    \else
    \par{#1}\par
    \fi
}

\def\rm{\fam0 \tenrm}

\def\fakesubhead#1\endsubhead{\bigskip\noindent{\bf#1}\par}



%
%
%

%

\font\textrsfs=rsfs10
\font\scriptrsfs=rsfs7
\font\scriptscriptrsfs=rsfs5

\newfam\rsfsfam
\textfont\rsfsfam=\textrsfs
\scriptfont\rsfsfam=\scriptrsfs
\scriptscriptfont\rsfsfam=\scriptscriptrsfs

\edef\oldcatcodeofat{\the\catcode`\@}
\catcode`\@11

\def\Cal@@#1{\noaccents@ \fam \rsfsfam #1}

\catcode`\@\oldcatcodeofat


\expandafter\ifx \csname margininit\endcsname \relax\else\margininit\fi

\long\def\red#1\endred{}
\long\def\green#1\endgreen{}
\long\def\blue#1\endblue{}
\long\def\private#1\endprivate{}

\def\endred{ \unmatched endred! }
\def\endgreen{ \unmatched endgreen! }
\def\endblue{ \unmatched endblue! }
\def\endprivate{ \unmatched endprivate! }

\ifx\latexcolors\undefinedcs\def\latexcolors{}\fi

\def\emptycs{}
\def\evaluatelatexcolors{%
        \ifx\latexcolors\emptycs\else
        \expandafter\xxevaluate\latexcolors\xxfertig\evaluatelatexcolors\fi}
\def\xxevaluate#1,#2\xxfertig{\setupthiscolor{#1}%
        \def\latexcolors{#2}}


\font\smallfont=cmsl7
\def\rutgerscolor{\ifmmode\else\endgraf\fi\smallfont
\advance\leftskip0.5cm\relax}
\def\setupthiscolor#1{\edef\tmptmpcs{\noexpand\bgroup\noexpand\rutgerscolor
\noexpand\def\noexpand\currentcolor{#1}%
\noexpand}%
\expandafter\let\csname#1\endcsname\tmptmpcs
\def\tmptmpcs{\checkColorUnmatched{#1}\popthecolor}
\expandafter\let\csname end#1\endcsname\tmptmpcs}

\def\checkColorUnmatched#1{\def\expectcolor{#1}%
    \ifx\expectcolor\currentcolor   
    \else \edef\failhere{\noexpand\tryingToClose '\currentcolor' with end\expectcolor}\failhere\fi}

\def\currentcolor{???}

\def\popthecolor{\ifmmode\else\endgraf\fi\egroup}

\expandafter\def\csname#1\endcsname{}

\evaluatelatexcolors

 \let\outerhead\head
 \def\head{\innerhead}
 \let\innerhead\outerhead

 \let\outersubhead\subhead
 \def\subhead{\innersubhead}
 \let\innersubhead\outersubhead

 \let\outersubsubhead\subsubhead
 \def\subsubhead{\innersubsubhead}
 \let\innersubsubhead\outersubsubhead

 \let\outerproclaim\proclaim
 \def\proclaim{\innerproclaim}
 \let\innerproclaim\outerproclaim

 %
 %
 %
 %

\def\demo#1{\medskip\noindent{\it #1.\/}}
\def\enddemo{\smallskip}

\def\remark#1{\medskip\noindent{\it #1.\/}}
\def\endremark{\smallskip}

\pageheight{8.5truein}
\topmatter
\title{Dependent first order theories, continued} \endtitle
\author {Saharon Shelah \thanks {The author would like to thank the
Israel Science Foundation for partial support of this research (Grant
No. 242/03). Publication 783. \null\newline
 I would like to thank 
Alice Leonhardt for the beautiful typing.} \endthanks} \endauthor 

\affil{ The Hebrew University of Jerusalem \\
 Einstein Institute of Mathematics\\
 Edmond J. Safra Campus, Givat Ram \\
 Jerusalem 91904, Israel
 \medskip
 Department of Mathematics \\
 Hill Center-Busch Campus \\
 Rutgers, The State University of New Jersey \\
 110 Frelinghuysen Road \\
 Piscataway, NJ 08854-8019  USA} \endaffil

\abstract  A dependent theory is a (first order complete theory) $T$
 which does not have the independence property.  A major result here is:
 if we expand a model of $T$ by the traces on it of sets definable in a
 bigger model then we preserve its being dependent.  Another one justifies the
 cofinality restriction in the theorem (from a previous work) saying
 that pairwise perpendicular indiscernible sequences, can have
 arbitrary dual-cofinalities in some models containing them.   We
 introduce ``strongly dependent $T$" and look at definable groups
 in such models; also look at forking and relatives.
\endabstract

\endtopmatter
\document

\newpage

\head {Annotated Content} \endhead
 \spuriousreset
\bn
\ub{Recall}:   Dependent $T=T$ without the independence property.
\bn
\S0 Introduction, pg.3
\bn
\S1 Expanding by making a type definable, p.4 
\mr
\item "{${{}}$}"  [Suppose we expand $M \prec {\frak C}$ by a relation
for each set of the form 
$\{\bar b:\bar b \in {}^m M$ and $\models \varphi[\bar b,\bar a]\}$,
where $\bar a \in {}^{\omega >} {\frak C},\varphi(\bar x,\bar y) \in
\Bbb L(\tau_T)$ and $m = \ell g(\bar x)$.  We prove that the theory
of this model is dependent and has elimination of quantifiers.]
\endroster
\bn
\S2  More on indiscernible sequences, p.15 
\mr
\item "{${{}}$}"  [This is complimentary to \cite[\S5]{Sh:715}.
Dedekind cuts with cofinality from both sides $\le \kappa + |T| 
= \kappa$ inside $\kappa$-saturated models (of a dependent theory $T$)
tend to be filled together.]
\endroster
\bn
\S3  Strongly dependent theories, p.26
\mr
\item "{${{}}$}"  [Being strongly dependent is related to being
superstable; however, strongly dependent theories which are stable (called
strongly stable) are not necessarily superstable.  We start the
investigation of this class of first order theories.  In particular,
for such a theory there is no non-algebraic types $p,q$ with definable
functions essentially from $q({\frak C})$ onto ${}^\omega(p({\frak
C}))$.  Also there is no equivalence relation on $p(\bar x)$ with
infinitely many equivalence classes, each class has essentially one to
one definable correspondence with the whole.]
\endroster
\bn
\S4  Definable groups, p.31
\mr
\item "{${{}}$}"  [We start to investigate definable groups for
dependent and strongly dependent theories, in particular with the size
of the commutator of most members.]
\endroster
\bn
\S5 Non-forking, p.39-56
\mr
\item "{${{}}$}"  [We try to see what does non-forking satisfy for
dependent theories.]
\endroster
\newpage

\head {\S0 Introduction} \endhead  \resetall \sectno=0
 \spuriousreset
\bn
The work in \cite{Sh:715} tries to deal with the investigation of a (first
order complete) theories $T$ which has the dependence property, i.e., 
does not have the independence property.

If $T$ is stable, expanding a model $M$ of $T$ by $p \restriction
\varphi(\bar x,\bar y)$ for $p \in \bold S^m(M)$, that is expanding $M$ by the
relation $R^{\varphi(\bar x,\bar y)}_{p,M} = \{\bar a \in 
{}^{\ell g(\bar y)}M:\varphi(\bar x,\bar a) \in p\}$ is an 
inessential one, i.e., by a relation on $M$ definable in $M$ with 
parameters.  This fails for unstable
theories but in \S1 we prove a weak relative: if $T$ is a dependent
theory then so is the expansion above, i.e., 
Th$(M,R^{\varphi(\bar x,\bar y)}_{p,M})_{p,\varphi(\bar x,\bar y)}$.

In \cite[\S5]{Sh:715} it is shown that for any model $N$ of
a dependent unstable $T$, we can find a $\kappa$-saturated model $M$
extending $N$ such 
that the following set is quite arbitrary: pairs of cofinalities of a
cut in $M$ for some definable partial order in $N$
(so not fulfilled in $M$) or even the set of pairs
$(\kappa_1,\kappa_2)$ of regular cardinals for which there is an
indiscernible sequence $\langle a_\alpha:\alpha < \kappa_1 \rangle
{}^\frown \langle b_\beta:\beta < \kappa^*_2 \rangle$ such that the
$(\kappa_1,\kappa_2)$-cut is respected in $M$, that is, we cannot find
an element in $M$ which we can add after the $a_\alpha$'s but before
the $b_\beta$'s linearly ordered by some $\varphi(x,y;\bar c)$ which is
a partial order.  However, there were
restrictions on the cofinalities being not too small.  
In section 2 we show that, to a large extent, yes these restrictions
are necessary.

The family of dependent theories is parallel to the family of stable
theories.  But actually better in some sense is the family of
superstable, i.e. the balance of the ``size" of the family of
such theories and what we can tell about them seem better for the family
of superstable ones.  In \S3 a related family, in a sense parallel to
superstable, called strongly dependent theories,
is defined.  Now every superstable $T$ is strongly stable (defined as
stable, strongly dependent), but the inverse fails (see
also \cite{Sh:839}, \cite{Sh:F660}).  We then observe some basic
properties.  This is continued in \cite{Sh:863}. 

In \S4 we look at groups definable in models of dependent theories,
and also in strongly dependent theories.  
In \S5 we try to look systematically at a parallel to non-forking.
\bn
This work is continued in \cite{Sh:876}, \cite{Sh:863}, \cite{Sh:886},
\cite{FiSh:E50}, \cite{CoSh:919}, \cite{Sh:F705}, \cite{Sh:877},
\cite{Sh:906}, \cite{Sh:900}.  
More specifically on a parallel to
uni-dimensional for the theory of the real field
see work of E. Firstenberg and S. Shelah
\cite{FiSh:E50}.  Concerning strongly dependent theories (see section
3) we try to investigate them in \cite{Sh:863}.
\nl
We should add to the history in \cite{Sh:715} that Keisler 
\cite{Ke87} connects
dependent theories and measures on the set of definable subsets of a
model.  Also Poizat \cite{Po81} (and then \cite[pg.202,3]{Sh:93}
answering positively a question of Poizat).  Poizat dealing with the number
of complete types in $\bold S(N)$ finitely satisfied in $M \prec N$; 
prove that the number is $\le 2^{\{M\|}$ (when $|T| \le \|M\|$)
and his question was whether it is $\le (\text{Ded}(\|M\|^{|T|})$ was
answered positively in \cite{Sh:93}.  Here in
\scite{ns.5} we follow \cite{Sh:93} proving that we can replace
finitely satisfiable but does not split.  Also \cite[3.2]{Sh:715} is
5.2 of Baldwin-Benedikt \cite{BalBl00}.
\nl
Note that Baisalov and Poizat \cite{BaPo98} proved that if $T$ is
o-minimal a theorem which is a consequence of \S1.
\nl
We thank Eyal Firstenberg, Aviv Tatarsky and the referee for many
helpful corrections and lately Itay Kaplan and friends for pointing
our deficiencies in \S5.
\nl
\ub{Notation}:
\sn
As in \cite{Sh:715} and, in addition
\definition{\stag{0.istP} Definition}  1) For $\bar{\bold b} = \langle
\bar b_t:t \in I \rangle$ an infinite indiscernible sequence, let
tp$'(\bar{\bold b}) = \langle \text{ tp}(\bar b_{t^n_0} \char 94
\ldots \char 94 \bar b_{t^n_{n-1}},\emptyset,{\frak C}):n < \omega \rangle$
where $t^n_\ell <_I t^n_{\ell,k}$ for $\ell < k < n < \omega$; the
choice of the $t^n_\ell$'s is immaterial. \nl
2) Let ``$M$ is $n$-saturated" mean ``$M$ is $\aleph_0$-saturated" for
$n < \omega$. \nl
3) Let $A/B$ mean tp$(A,B)$, inside ${\frak C}$ or ${\frak C}^{\text{eq}}$,
of course.
\enddefinition
\newpage

\head {\S1 Expanding by making a type definable} \endhead  \resetall \sectno=1
 \spuriousreset
\bigskip

What, in short, do we show here?  We say $A$ is full over $M$ if every
$p \in \bold S^{< \omega}(M)$ is realized in $A$, (Definition
\scite{t.1.3}).  We let ${\frak B}_{A,M}$ be the expansion of $M$, for each
$\varphi(\bar x,\bar a),\bar a \in {}^{\omega >} A$, by the following $\ell
g(\bar x)$-place relation:  all realizations of $\varphi(-,\bar a)$,
i.e., by $\varphi(M,\bar a)$ (see Definition \scite{t.3}(2)).  
We prove here that if $A$ is full over $M$
then Th$({\frak B}_{M,A})$ has elimination of quantifiers (see Claim
\scite{t.4A}(1), its proof depends only on \scite{t.1}, \scite{t.1B}(2)).  By
this we prove that Th$({\frak B}_{M,A})$ is dependent (in
\scite{t.5A} depending on \scite{t.3}(4), \scite{t.4}(1),(5) only), 
so for this conclusion ``$A$ is full over $M$" is not needed.
\bigskip

\demo{\stag{t.0} Context}  1) $T$ is a (first order complete) dependent
theory in the language $\Bbb L(\tau_T)$.
\nl
2) ${\frak C} = {\frak C}_T$ is a monster model for $T$.
\enddemo
\bigskip

\proclaim{\stag{t.1} Claim}   Assume
\mr
\item "{$(a)$}"  $M$ a model
\sn
\item "{$(b)$}"  $D$ an ultrafilter on $M$, i.e. on the Boolean
Algebra ${\Cal P}(M)$.
\ermn
\ub{Then} for any $\bar c \in {}^{\omega >}{\frak C}$ and formula
$\varphi(x,y,\bar c)$ we have: if the set $\{a \in M:(\exists y \in
M)({\frak C} \models \varphi[a,y,\bar c])\}$ belongs to $D$ then it
belongs to ${\text{\rm def\/}}_2(D)$, see definition below.
\endproclaim
\bigskip

\definition{\stag{t.1.4} Definition}  1) When $D$ is an ultrafilter on a
set $B \subseteq {\frak C}$ let ${\text{\rm def\/}}_2(D) = \{A \in
D$: some member of def$_1(D)$ is included in $A\}$ where def$_1(D)
= \{A \in D$: for some $\bar c \in {}^{\omega
>}{\frak C}$ and formula $\psi(x,\bar c)$ the set $\psi(M,\bar c) = 
\{a \in M:{\frak C} \models \psi(a,\bar c)\}$ 
belongs to $D$ and is equal to $A\}$.
\nl
2) Similarly when $D$ is an ultrafilter on ${}^m B,m < \omega$.
\enddefinition
\bigskip

\remark{\stag{t.1A} Remark}  Note the following easy comments.
\nl
1) Of course, Claim \scite{t.1} holds also for 
$\varphi = \varphi(\bar x,\bar y,\bar c)$ when $D$ an ultrafilter 
on ${}^m M$ and $m = \ell g(\bar y)$ because, e.g. we can just work 
in ${\frak C}^{\text{eq}}$. 
\nl
2) $T$ is dependent iff $T^{\text{eq}} = 
\text{\rm Th}({\frak C}^{\text{eq}})$ is so (1) is dependent this
justify the statement above in part (1) and 
Th$({\frak C})$ is dependent iff Th$({\frak C},c)_{c \in C}$ is (for
any $C \subseteq {\frak C}$) and $T$ dependent $\Rightarrow \text{\rm
Th}({\frak C} \restriction \tau')$ is dependent when $\tau' \subseteq
\tau_T$.
\nl
3) Note that def$_1(D)$ is a filter on $A$.
\nl
4) In the proof of \scite{t.1} the hypothesis 
``$T$ dependent" is used only for deducing ``$\varphi(x,y,\bar c)$
is dependent" which is naturally defined.
\nl
5) Recall the following (which is used in the proof):
\mr
\item "{$(a)$}"  $\Delta \subseteq \Bbb L(\tau_T)$ means $\Delta$ is a
set of objects of the form $\varphi(\bar x),\varphi$ a (first order)
formula from $\Bbb L(\tau_T),\bar x$ a sequence of
variables with no repetitions including the free variables, but
changing the variables of $\varphi$ 
is allowed here, i.e., there is no difference
between $\varphi(x)$ and $\varphi(y)$; we may write 
$\varphi(\bar x,\bar y)$ instead of $\varphi(\bar x \char 94 \bar y)$
\sn 
\item "{$(b)$}"  tp$_\Delta(\bar a,A) = \{\varphi(\bar x,\bar b):\bar
x = \langle x_\ell:\ell < \ell g(\bar a)\rangle,\varphi(\bar x,\bar y)
\in \Delta$ and ${\frak C} \models \varphi[\bar a,\bar b]$ and $b \in
{}^{\omega >}A\}$
\sn
\item "{$(c)$}"  $\langle \bar b_t:t \in I\rangle$ is
$\Delta$-indiscernible over $B$ means that $I$ is a linear order and
if $\varphi(\bar x_1,\dotsc,\bar x_n,\bar y) 
\in \Delta,\ell g(\bar x_\ell) = \ell
g(\bar b_t)$ for $\ell=1,\dotsc,n$ and $t \in I$ and $\bar c \in
{}^{\ell g(\bar y)}B$ then for any $s_1 <_I \ldots <_I s_n$ and $t_1 <_I
\ldots <_I t_n$ we have ${\frak C} \models ``\varphi[\bar
b_{t_1},\dotsc,\bar b_{t_n},\bar c] \equiv \varphi[\bar
b_{s_1},\dotsc,\bar b_{s_n},\bar c]"$. 
\ermn
6) In the proof of Claim \scite{t.1} we 
do not need to close $\Delta_1$ to $\Delta_2$, i.e. we can let
   $\Delta_2 = \Delta_1$ \ub{provided that} we redefine
   tp$_{\Delta_1}(a,A)$ as

$$
\align
\text{tp}(a,A) \cap
   \{\varphi(a_0,\dotsc,a_{m-1},x,a_{m+1},\dotsc,a_n):&\varphi(x_0,
\dotsc,x_{m-1},x_m,   \\
  &x_{m+1},\dotsc,x_{n-1}) \in \Delta\}
\endalign
$$
\mn
or more specifically, in $(*)_1$ from $\boxtimes_1$, inside the proof
of $\boxtimes_1$, we replace ``$a_\ell$ realizes
tp$_{\Delta_2}(a_w,\ldots)$" by ``$a_\ell$ realizes
$\{\varphi(a_{\ell_0},\dotsc,a_{\ell_{m-1}},x,a_{\omega
+1},\dotsc,a_{\omega +n-1+m},\bar b):\varphi(x_0,\dotsc,x_{n-1},\bar
y) \in \Delta_1$ and $\bar b \in {}^{\ell g(\bar y)}B,m<n,\ell_0 <
\ldots < \ell_{m-1} < \ell$ and ${\frak C} \models
\varphi[a_{\ell_0},\dotsc,a_{\ell_{m-1}},a_\omega,a_{\omega
+n+1},\dotsc,a_{\omega +n-1-m},\bar b)\}$. 
\endremark
\bigskip

\demo{Proof}  We shall use ``$T$ is dependent" only in the last sentence
of the proof toward contradiction.  Assume that $\bar c,
\varphi(x,y,\bar c)$ form a counterexample.
\nl
So 
\mr
\item "{$\circledast_0$}"  $(i) \quad$ the set 
$A^* = \{a \in M$: for some $b \in M$ we have $\models
\varphi[a,b,\bar c]\}$ belongs \nl

$\quad$ to $D$, 
\sn
\item "{${{}}$}"  $(ii) \quad A^* \notin \text{ def}_2(D)$, that is, no
$A' \in \text{ def}_1(D)$ is included in $A^*$.
\ermn
By the choice of $A^*$ we can, for each $a \in A^*$, choose
$b_a \in M$ such that $\models \varphi[a,b_a,\bar c]$.  Let $D_1 = D$
and let $D_2$ be
the following ultrafilter on ${}^2 M:X \in D_2$ iff $X \subseteq
{}^2 M$ and for some $A \in D$ we have $\{(a,b_a):a \in A \cap A^*\}
\subseteq X$.

We can choose $\langle (a_{\omega +n},b_{\omega +n}):n < \omega
\rangle$ from ${\frak C}$ such that
\mr
\item "{$\circledast_1$}"   for $n_1 < n_2 < \omega$ the pair
$(a_{\omega + n_1},b_{\omega +n_1})$ realizes the type Av$(M \cup
\{a_{\omega + \ell},b_{\omega + \ell}:\ell \in (n_1,n_2]\},D_2)$.
\ermn
It follows that
$a_{\omega+n_1}$ realizes the type Av$(M \cup\{a_{\omega +
\ell},b_{\omega + \ell}:\ell \in (n_1,n_2]\},D_1)$ and
\mr
\item "{$\circledast_2$}"   for $n_1 < n_2$, the element $a_{\omega +
n_1}$ realizes the type Av$(M \cup\{a_{\omega + \ell}:\ell \in
(n_1,n_2]\},D)$
\sn
\item "{$\circledast_3$}"   for $n_1 < n_2$ the triple
$(a_{2n_1},a_{2n_1+1},b_{2n_1+1})$ realizes the type Av$(M \cup
\{a_{\omega + 2\ell},a_{\omega + 2 \ell +1},b_{\omega +2
\ell+1}:\ell \in (n_1,n_2]\},D_3)$ for some ultrafilter $D_3$ on
$M^3$.
\nl
[Why?  We define $D_3 := \{X \subseteq M^3:\{a \in M:\{(b,c) \in M
\times M:(a,b,c) \in X\} \in D_2\} \in D_1\}$.]
\nl
(we use mainly $\circledast_1$).
\ermn
Now clearly
\mr
\item "{$\boxtimes_0$}"   $\langle(a_{\omega +n},b_{\omega +n}):n <
\omega \rangle$ is an indiscernible sequence over $M$
\sn
\item "{$\boxtimes_1$}"  if $\Delta_1 \subseteq \Bbb L(\tau_T)$ is finite 
\ub{then} we can find $n(*) < \omega$ and finite
$\Delta_2 \subseteq \Bbb L(T)$ such that
{\roster
\itemitem{ $(*)_1$ }  if $n_1 < \omega$ and $B \subseteq
M$ is finite and for each $\ell < n_1$ the element $a_\ell \in M$
realizes the type
\nl
tp$_{\Delta_2}(a_\omega,\{a_0,\dotsc,a_{\ell-1}\} \cup \{a_{\omega
+1},\dotsc,a_{\omega +n(*)}\} \cup B)$ \ub{then} 
$\langle a_\ell:\ell < n_1 \rangle
\char 94 \langle a_{\omega + \ell}:\ell < \omega\rangle$ is a
$\Delta_1$-indiscernible sequence over $B$ (and even
$\Delta_2$-indiscernible).  
\endroster}
Note that this is close to \cite[1.16=np1.5tex]{Sh:715}; note that it 
follows from the result that even for $n_1 = \omega$ this holds.
\ermn
[Why does $\boxtimes_1$ hold?  Let $n(*)$ be 
arity$(\Delta_1)$, i.e., the maximal
number of free variables of a formula from $\Delta_1$, it is finite as
$\Delta_1$ is finite so without loss of generality 
each $\varphi \in \Delta_1$ is
$\varphi(\bar x)$, Rang$(\bar x) \subseteq \{x_\ell:\ell < n(*)\}$.
Let $\Delta_2$ be the closure of $\Delta_1$
under identifying and
 permuting the variables and let $\Delta_{2,k}$ be defined like
$\Delta_2$ but we allow to add from $\{x_0,\dotsc,x_k\}$ 
dummy variables to each formula (we can use below $\Delta_2 =
\cup\{\Delta_{2,k}:k < \omega\}$).
We have to prove that for this choice of $n(*)$ and $\Delta_2$ the
assertion $(*)_1$ holds.

So assume $n_1 < \omega$ and $B,a_\ell$ (for $\ell < n_1$) are as
required in the assumption of $(*)_1$.
Now we prove by induction on $k \le n_1$ that
\mr
\item "{$(*)^1_k$}"   the sequences $\langle
a_{\omega + \ell}:\ell < n_1 + n(*)\rangle$ and $\langle a_\ell:\ell < k
\rangle \char 94 \langle a_{\omega + \ell}:\ell < n_1 + n(*)-k\rangle$
realize the same $\Delta_{2,n_1+n(*)}$-type over $B$
which means that: if $m \le n(*),\bar d \in {}^m B$ and $\varphi(\bar
y_1,\bar y_2) \in \Delta_{2,n_1 +n(*)+m},\ell g(\bar y_1) = n_1 +
n(*),\ell g(\bar y_2) = m$ then ${\frak C} \models \varphi[\langle
a_{\omega + \ell}:\ell < n_1+n(*)\rangle,\bar d]$ iff ${\frak C}
\models \varphi[\langle a_\ell:\ell < k\rangle \char 94 \langle
a_{\omega + \ell}:\ell < n_1 + n(*)-k\rangle,\bar d]$; note that we
can allow $m \le n(*)$.
\ermn
For $k=0$, the two expressions gives the same sequence.  
Assume this holds for $k$ and we shall prove it for $k+1$.  First
$\langle a_\ell:\ell < k+1 \rangle \char 94 \langle a_{\omega +
\ell}:\ell <n_1 + n(*)-(k+1)\rangle$ realize the same type as
$\langle a_0,\dotsc,a_k,a_{\omega +1},\dotsc,a_{\omega +n_1+n(*)-(k+1)}\rangle$
simply as $\langle a_{\omega +\ell}:\ell < \omega\rangle$ is an
indiscernible sequence over $M$, by $\boxtimes_0$.  Now by the
assumption of $(*)_1$ we know that $a_k,a_\omega$
realizes the same $\Delta_2$-type over $B \cup \{a_0,\dotsc,a_{k-1}\}
\cup \{a_{\omega +1},\dotsc,a_{\omega + n_1+n(*)-k}\}$.

As $n(*)$ is
the arity of $\Delta_1$ hence also of $\Delta_2$ and the definition of
$\Delta_{2,n_1+n(*)}$ it follows that the sequence $\langle
a_0,\dotsc,a_{k-1},a_k,a_{\omega +1},\dotsc,a_{\omega
+n_1+n(*)-k-1}\rangle$ realizes over $B$ the same
$\Delta_{2,n_1+n(*)}$-type as the sequence \nl
$\langle
a_0,\dotsc,a_{k-1},a_\omega,a_{\omega+1},\dotsc,a_{\omega +n_1+n(*)-k-1}
\rangle$ but by the induction hypothesis on $k$ the latter
realizes over $B$ the same $\Delta_{2,n_1+n(*)}$-type as the sequence $\langle
a_\omega,a_{\omega+1},\dotsc,a_{\omega+n_1+n(*)-1}\rangle$,
hence $(*)^1_{k+1}$ holds  so we have
carried the induction on $k \le n_1$.
 Now the desired conclusion follows from $(*)^1_m$
by $\boxtimes_0$ as each formula in $\Delta_1$ and even $\Delta_2$ has
$\le n(*)$ free variables.]
\mr
\item "{$\boxtimes_2$}"  if $\Delta_1 \subseteq \Bbb
L(\tau_T)$ is finite \ub{then} we can find $n(*) < \omega$ and finite
$\Delta_2 \subseteq \Bbb L(\tau_T)$ such that
{\roster
\itemitem{ $(*)_2$ }  \ub{if} $n_1 < \omega,B \subseteq M$ is finite and
for each $\ell < n_1,a_{2 \ell} \in M$ realizes
\nl
tp$_{\Delta_2}(a_\omega,\{a_{2m},a_{2m+1},b_{2m+1}:m < \ell\} \cup 
\{a_{\omega + \ell},b_{\omega + \ell}:\ell =1,\dotsc,n(*)\} \cup B)$ and 
$\langle a_{2 \ell+1},b_{2 \ell +1} \rangle$ realizes 
tp$_{\Delta_2}((a_\omega,b_\omega),\{a_{2m},a_{2m+1},b_{2m+1}:m <
\ell\} \cup \{a_{2 \ell}\} \cup \{a_{\omega + \ell},b_{\omega + \ell}:
\ell =1,\dotsc,n(*)\} \cup B)\})$ \ub{then} $\langle(a_{2 \ell},a_{2 \ell +1},
b_{2 \ell +1}):\ell < n_1 \rangle \char 94 
\langle a_{\omega + 2 \ell},a_{\omega +2\ell +1},b_{\omega + 2 \ell
+1}:\ell < \omega \rangle)$ is $\Delta_1$-indiscernible over $B$ (and
even $\Delta_2$-indiscernible)
\endroster}
[Why?  The proof is similar to the proof of $\boxtimes_1$ mainly
replacing the use of $\circledast_1$ by $\circledast_3$.]
\sn
\item "{$\boxtimes_3$}"  if $B \subseteq M$ is finite, $n^* < \omega$
and $\Delta \subseteq \Bbb L(\tau_T)$ is finite, 
\ub{then} we can find $a \in M$ realizing the finite type $q = \text{
tp}_\Delta(a_\omega,B \cup \{a_{\omega + \ell},b_{\omega + \ell}:
\ell=1,\dotsc,n^*\})$ such that $\models \neg(\exists y \in
M)\varphi(a,y,\bar c)$. \nl
[Why?  The set $A := \{a \in M:a$ realizes $q$, equivalently satisfies
the formula $\wedge q \in \text{ Av}({\frak C},D)\}$ belongs to $D$
because $q$ is finite and the choice of 
$\langle a_{\omega + \ell},b_{\omega + \ell}:\ell < \omega \rangle$; 
moreover, it belongs to def$_1(D)$ by the definition of def$_1(D)$ as
$\wedge q$ is a formula.  But def$_1(D) \subseteq \text{ def}_2(D)$
hence $A \in \text{ def}_2(D)$.

So by the assumption toward contradiction and
choice of $A^*$, i.e., by $(*)_0$, we have $\neg(A \subseteq A^*)$ 
so there is $a \in A$ such that $a \notin A^*$ which means that
$\neg(\exists y \in M)\varphi(a,y,\bar c)$, so we are done.]
\ermn
By the above and compactness (or use an ultrapower)
\mr
\item "{$\boxtimes_4$}"  there are $N,a_{2n},a_{2n+1},b_{2n+1}$ 
(for $n < \omega$) such that
{\roster
\itemitem{ $(a)$ }  $N$ is $|T|^+$-saturated
\sn
\itemitem{ $(b)$ }  $a_{2n},a_{2n+1},b_{2n+1} \in N$
\sn
\itemitem{ $(c)$ }  $\langle a_n:n < \omega \rangle$ is an
indiscernible sequence
\sn
\itemitem{ $(d)$ }  $\langle (a_{2n},a_{2n+1},b_{2n+1}):n < \omega
\rangle \char 94 \langle (a_{\omega+2n},a_{\omega +2n+1},b_{\omega
+2 n+1}):n < \omega \rangle$ is an indiscernible sequence
\sn
\itemitem{ $(e)$ }  ${\frak C} \models \varphi[a_{2n+1},b_{2n+1},\bar c]$
\sn
\itemitem{ $(f)$ }  for no $n < \omega$ and $b \in N$ do we have
${\frak C} \models \varphi[a_{2n},b,\bar c]$.
\endroster}
\ermn
[Why?  By compactness it is enough to prove the following: 
for every $n_1 < \omega$ and finite
$\Delta_1 \subseteq \Bbb L(\tau_T)$ to which $\varphi$ belongs there are
$a_{2n},a_{2n+1},b_{2n+1} \in M$ for $n<n_1$ such that 
clauses (a)-(f) holds when
we restrict ourselves to $n<n_1$ and $\Delta_1$-types replacing $N$ by
$M$.  We first choose a finite $\Delta_2 \subseteq \Bbb L(\tau_T)$ as in
$\boxtimes_2$, and then choose $(a_{2n},a_{2n+1},b_{2n+1})$ by
induction on $n$ such that the demand in $(*)_2$ of $\boxtimes_2$
holds.  Arriving to $n$ choose $a_{2n} \in M$ such that in addition, clause
(f) holds, this is possible by $\boxtimes_3$, and then choose
$(a_{2n+1},b_{2n+1}) \in {}^2 M$ recalling $(a_\omega,b_\omega)$ realizes Av$(M
\cup \{a_{\omega +n},b_{\omega +n}:1 \le n < \omega\},D_2)$.  So we
are done proving $\boxtimes_4$.]
\bn
Next by clause (d) of $\boxtimes_4$
\mr
\item "{$\boxtimes_5$}"  there is an automorphism $F$ of ${\frak C}$
such that \nl
$n < \omega$ implies $F((a_{\omega +2n},a_{\omega+2n+1},b_{\omega+2n+1})) =
(a_{2n},a_{2n+1},b_{2n+1})$.
\ermn
Hence we can find $b_{2n} \in {\frak C}$ for $n < \omega$
such that $\langle (a_n,b_n):n < \omega \rangle$ is an indiscernible
sequence (over $\emptyset$, not necessarily over $\bar c$!) 
and as $N$ is $|T|^+$-saturated \wilog \, $b_{2n} \in N$ for
$n < \omega$.  But ${\frak C} \models \varphi
[a_{2n+1},b_{2n+1},\bar c]$ for $n < \omega$ by clause (e) of $\boxtimes_4$
so as $T$ is dependent for every large enough $n < \omega$ we have 
${\frak C} \models \varphi[a_{2n},b_{2n},\bar c]$.  But
as $b_{2n} \in N$ clearly $\{a_n,b_n:n < \omega\} \subseteq N$ hence
$n < \omega \Rightarrow {\frak C} \models \varphi[a_{2n},b_{2n},\bar c]$
contradicting clause (f) of $\boxtimes_4$.
\hfill$\square_{\scite{t.1}}$
\enddemo
\bn
Recall
\definition{\stag{t.1.3} Definition}  For $A \subseteq C (\subseteq
{\frak C})$ we say that $C$ is full over $A$ when: for every $m
<\omega$ and $p \in \bold S^m(A)$, there is $\bar c \in {}^m C$ which 
realizes $p$.
\enddefinition
\bigskip

\demo{\stag{t.1.7} Observation}  
We have def$_\ell(D_1) = \text{ def}_\ell(D_2)$ for $\ell =1,2$
\underbar{when}:
\mr
\item "{$(a)$}"   $D_1,D_2$ are ultrafilters on ${}^m A$,
\sn
\item "{$(b)$}"  $A \subseteq C$,
\sn
\item "{$(c)$}"   $C$ is full over $A$, 
\sn
\item "{$(d)$}"   Av$(C,D_1) = \text{ Av}(C,D_2)$.
\endroster
\enddemo
\bigskip

\demo{Proof}  Easy.  \hfill$\square_{\scite{t.1.7}}$
\enddemo
\bigskip

\proclaim{\stag{t.1B} Claim}  1) Assume
\mr
\item "{$(a)$}"  $M \subseteq C$
\sn
\item "{$(b)$}"  $D_0$ is an ultrafilter on ${}^{m_0} M$
\sn
\item "{$(c)$}"  $\bar b_0$ realizes {\rm Av}$(C,D_0)$
\sn
\item "{$(d)$}"  {\rm tp}$(\bar b_0 \char 94 \bar b_1,C)$ is f.s. in $M$ and
$m_1 = \ell g(\bar b_1)$
\sn
\item "{$(e)$}"  $C$ is full over $M$.
\ermn
\ub{Then} for some ultrafilter $D$ on ${}^{m_0+m_1}M$ we have
\mr
\item "{$(\alpha)$}"  {\rm Av}$(C,D) = \text{\rm tp}
(\bar b_0 \char 94 \bar b_1,C)$
\sn
\item "{$(\beta)$}"  the projection of $D$ on ${}^{m_0}M$ is $D_0$.
\ermn
2)  Assume that clauses (a) and (e) of part (1) hold.
\ub{Then} for any $\bar c \in {}^{\omega >}{\frak C}$ and formula
$\varphi(\bar x,y,\bar z) \in \Bbb L(\tau_T),\ell g(\bar z) 
= \ell g(\bar c)$ there are $\psi(\bar x,\bar z') \in \Bbb L(\tau_T)$ and 
$\bar d$ of length $\ell g(\bar z')$ from
${\frak C}$, (and even from $C$) such that
$\{\bar a \in M:(\exists y \in M)({\frak C} 
\models \varphi[\bar a,y,\bar c])\} = 
\{\bar a \in M: {\frak C} \models \psi(\bar a,\bar d)\}$. 
\endproclaim
\bigskip

\demo{Proof}  1) Let

$$
{\Cal E}_0 = \bigl\{ \{\bar a \in {}^{m_0+m_1} M:\bar a \restriction
m_0 \in X\}:X \in D_0 \bigr\}
$$

$$
\align
{\Cal E}_1 = \bigl\{ \{\bar a \in {}^{m_0+m_1}M:&{\frak C} \models \varphi[\bar
a;\bar c]\}:\varphi(\bar x;\bar y) \in \Bbb L(\tau_T) \\
  &\ell g(\bar x) = m_0 + m_1,\ell g(\bar y) = \ell
  g(\bar c), \\
  &\bar c \in {}^{\omega >} C \text{ and } {\frak C}
\models \varphi[\bar b_0 {}^\frown \bar b_1; \bar c] \bigr\}.
\endalign
$$
\mn
Clearly it suffices to prove that there is an ultrafilter on ${}^{m_0
+ m_1}M$ extending ${\Cal E}_0 \cup {\Cal E}_1$.  For this it suffices
to show that any finite subfamily of ${\Cal E}_0 \cup {\Cal E}_1$
has a non-empty intersection.  But ${\Cal E}_0$ is closed under finite
intersections as $D_0$ is an ultrafilter on ${}^{m_0}M$ and 
${\Cal E}_1$ is closed under finite intersections as $\Bbb L(\tau_T)$ is
closed under conjunctions, so it suffices to prove that $X_0 \cap X_1
\ne \emptyset$ when
\mr
\item "{$(i)$}"  $X_0 = \{\bar a \in {}^{m_0+m_1} M:\bar a \restriction
m_0 \in X\} \in {\Cal E}_0 \text{ for some } X \in D_0$
\sn
\item "{$(ii)$}"  $X_1 = \{\bar a \in {}^{m_0+m_1} M:{\frak C} \models
\varphi[\bar b_0 \char 94 \bar b_1;\bar c]\} 
\in {\Cal E}_1$ where $\varphi(\bar x,\bar y)$ and $\bar c$ 
are as in the definition of ${\Cal E}_1$.
\ermn
As tp$(\bar b_0 \char 94 \bar b_1,C)$ is finitely satisfiable in 
$M$ (= assumption $(d)$) clearly there is 
an ultrafilter $D'_1$ on ${}^{m_0+m_1}M$ such that Av$(C,D'_1) = 
\text{ tp}(\bar b_0 \char 94 \bar b_1,C)$. 

Let $D'_0$ be the projection of $D'_1$ to ${}^{m_0}M$, i.e., $\{Y
\subseteq {}^{m_0}M:\{\bar a \in {}^{m_0+m_1}M:\bar a \restriction m_0
\in Y\} \in D'_1\}$.  Clearly $D'_0$ is an ultrafilter over ${}^{m_0}M$.
We have ${\frak C} \models \varphi[\bar b_0,\bar b_1;\bar c]$,
so $X_1 \in D'_1$ hence $X'_0 = \{\bar a \restriction m_0:\bar a
\in X_1\} \in D'_0$; which implies that the set
$X''_0 := \{\bar a_0 \in {}^{m_0}M$: for some $\bar a_1 \in {}^{m_1}M$ 
we have $\bar a_0 \char 94 \bar a_1 \in
X_1$, i.e., $\models \varphi[\bar a_0,\bar a_1;\bar c]\}$ belongs to
$D'_0$. 
\nl
By \scite{t.1} (and \scite{t.1A}(1),(2)) it 
follows that $X''_0$ includes some $Y''_0 
\in \text{ def}_1(D'_0)$.  Now Av$(C,D_0) = \text{ tp}(\bar b_0,C) = 
\text{ Av}(C,D'_0)$, because the first equality holds as by
assumption (b) the sequence $\bar b_0$ realizes Av$(C,D_0)$ and
second equality holds as $\bar b_0 \char 94 \bar b_1$ realizes
Av$(C,D'_1)$ and the choice of $D'_0$.  But by assumption $(e)$ every
$p \in \bold S^{< \omega}(M)$ is realized by some sequence from $C$, hence
by Observation \scite{t.1.7} we have def$_2(D_0) = \text{ def}_2(D'_0)$. 
But $Y''_0 \in \text{\rm def}_1(D'_0)$ so $Y''_0 \in \text{\rm
def}_2(D_0)$ hence $Y''_0 \in D_0$.  By the choice of $Y''_0$ we have
$Y''_0 \subseteq X''_0 \subseteq {}^{m_0}M$ so by the previous sentence
$X''_0 \in D_0$, but by clause (i) above also $X \in D_0$ hence 
$X \cap X''_0 \in D_0$, so we
can find $\bar a_0 \in X \cap X''_0 \subseteq {}^{m_0} M$.  By the
definition of $X''_0$ there is $\bar a_1 \in {}^{m_1}M$ such that
${\frak C} \models \varphi[\bar a_0,\bar a_1;\bar c]$.  Now $\bar a_0 \char 94
\bar a_1 \in X_1$ by the definition of $X_1$ from 
clause (ii) and $\bar a_0 \char 94 \bar a_1 \in X_0$ because $\bar a_0 \in
X$ and $X_0$'s definition from clause $(i)$.  So 
$\bar a_0 \char 94 \bar a_1 \in X_0 \cap X_1$ hence $X_0 \cap
X_1 \ne \emptyset$ and we are done. \nl
2) Let $\varphi^*(\bar x,y,\bar z) \in \Bbb L(\tau_T)$ and $\bar c^*
\in {}^{\ell g(\bar z)}{\frak C}$ and we should find 
$\psi(\bar x,\bar z'),\bar d$ 
as required.  Let $\bar c \in {}^{\ell g(\bar z)}C$ realizes tp$(\bar
c^*,M)$ for our purpose we may assume \wilog \, $\bar c^* = \bar c$.
For any formula $\psi(\bar x,\bar z') \in \Bbb L(\tau_T)$ and 
$\bar d \in {}^{\ell g(\bar z')}{\frak C}$ let 
$Y_{\psi(\bar x,\bar d),M} = \{\bar a \in {}^{\ell g(\bar x)}M:
{\frak C} \models \psi[\bar a,\bar d]\}$ and let
$X_{\varphi(\bar x,y,\bar c),M} = \{\bar a \in {}^{\ell g(\bar x)}M:
{\frak C} \models \varphi[\bar a,b,\bar c]$ for some $b \in M\}$.

Lastly, let ${\Cal P} = \{Y_{\psi(\bar x,\bar d),M}:\psi(\bar x,\bar y) \in
\Bbb L(\tau_T),\bar d \in {}^{\ell g(\bar z)} C$ and $Y_{\psi(\bar
x,\bar d),M} \subseteq X_{\varphi(\bar x,y,\bar c),M}\}$.
\nl
Clearly ${\Cal P}$ is closed under finite unions and is a family of subsets
of $M$.  Also if $X_{\varphi(\bar x,y,\bar c),M}$
is equal to some member of ${\Cal P}$ then we are 
done, so assume toward contradiction
that this fails.  So as $X_{\varphi(\bar x,y,\bar c)} \subseteq M$ there
is an ultrafilter $D$ on $M$ such that $X_{\varphi(\bar x,y,\bar c),M} \in
D$ but $D$ is disjoint to ${\Cal P}$ which contradicts \scite{t.1}.   
\hfill$\square_{\scite{t.1B}}$   
\enddemo
\bigskip

\demo{\stag{t.2} Conclusion}  Assume
\mr
\item "{$(a)$}"  $M \prec M_1$
\sn
\item "{$(b)$}"   $M_1$ is $\|M\|^+$-saturated.
\ermn
\ub{Then} $\{A:A/M_1$ is f.s. in $M\}$ has amalgamation and JEP (the joint
embedding property) by elementary maps from ${\frak C}$ to ${\frak
C}$ which are the identity on $M_1$.
\enddemo
\bigskip

\demo{Proof}  The joint embedding property is trivial.  For the
amalgamation, by compactness we should consider finite sequence $\bar
a_0,\bar a_1,\bar a_2$ such that tp$(\bar a_0 \char 94 \bar
a_\ell,M_1)$ is f.s. in $M$ for $\ell = 1,2$
and we should find sequences $\bar b_0,
\bar b_1,\bar b_2$ such that $\ell g(\bar b_\ell) = \ell g(\bar
a_\ell)$ for $\ell=0,1,2$ and 
tp$(\bar a_0 \char 94 \bar a_\ell,M_1) = \text{ tp}(\bar
b_0 \char 94 \bar b_\ell,M_1)$ for $\ell =1,2$ and tp$(\bar b_0 \char
94 \bar b_1 \char 94 \bar b_2,M_1)$ is f.s. in $M$. \nl
Let $m_\ell = \ell g(a_\ell)$, let $D_0$ be an ultrafilter on
${}^{m_0}M$ such that tp$(\bar a_0,M_1) = \text{ Av}(M_1,D_0)$.  By
\scite{t.1B}(1) for $\ell \in \{1,2\}$ there is an ultrafilter $D_\ell$ on
${}^{m_0+m_\ell}M$ such that
\mr
\item "{$(*)_1$}"  tp$(\bar a_0 \char 94 \bar a_\ell,M_1)$ is
Av$(M_1,D_\ell)$
\sn
\item "{$(*)_2$}"   the projection of $D_\ell$ on ${}^{m_0}M$ is
$D_0$.
\ermn
Let $m = m_0 + m_1 + m_2$ and let $D'_1$ be the filter on ${}^m M$
consisting of $\{Y \subseteq {}^m M$: for some $X \in D_1$ for every
$\bar a \in {}^m M$ we have $\bar a \restriction (m_0 + m_1) \in X
\Rightarrow \bar a \in Y\}$ and let $D'_2$ be the filter on ${}^m M$
consisting of $\{Y \subseteq {}^m M$: for some $X \in D_2$ for every
$\bar a \in {}^m M$ we
have $(\bar a \restriction m_0) \char 94 (\bar a \restriction [m_0 +
m_1,m)) \in X \Rightarrow \bar a \in Y\}$.  Easily $Y_1 \in D'_1 \and
Y_2 \in D'_2 \Rightarrow Y_1 \cap Y_2 \ne \emptyset$ because $D_1,D_2$
has the same projection on ${}^{m_0}M$.

Hence we can find an ultrafilter $D^*$ on ${}^{m_0+m_1+m_2}M$ which
extends $D'_1 \cup D'_2$ hence if 
$\bar b_0 \char 94 \bar b_1 \char 94 \bar b_2$ realizes
Av$(M_1,D^*)$ then $\bar b_0 \char 94 \bar b_\ell$ realizes tp$(\bar a_0
\char 94 \bar a_\ell,M_1)$ for $\ell=1,2$.  So we are done.
\hfill$\square_{\scite{t.2}}$
\enddemo
\bn
\margintag{t.2.7}\ub{\stag{t.2.7} Discussion}:  Next we 
shall deduce the promised results.  If $M^+$ is an
expansion of a model $M \prec {\frak C}$ by the restriction of
relations definable in ${\frak C}$ (with parameters) \ub{then}
Th$(M^+)$ is still dependent.  Moreover, if we do this for close
enough family of such relations then Th$(M^+)$ has elimination of
quantifiers.  Toward formulating this result we define several
extensions of $T$.
\bigskip

\definition{\stag{t.3} Definition}  Let $M \prec {\frak C},
A \subseteq {\frak C}$ and for simplicity $\tau_T$ has predicate
symbols only. \nl
1) We define a universal first order theory $T_{M,A}$ as follows
\mr
\item "{$(a)$}"  the vocabulary is \nl
$\tau_{M,A} = \{P_{\varphi(\bar x,\bar a)}:\varphi \in \Bbb L(\tau_T)$
and $\bar a \in {}^{\ell g(\bar y)}A\} \cup \{c_a:a \in M\}$ \nl
with
{\roster
\itemitem{ $(i)$ }  $c_a$ an individual constant
\sn
\itemitem{ $(ii)$ }  $P_{\varphi(\bar x,\bar a)}$ being a 
predicate with arity $\ell g(\bar x)$; \ub{but} we identify 
$P_{R(\bar x)}$ with $R$ (where $\bar x = \langle x_\ell:\ell <
\text{\rm arity}(R)\rangle$)  so $\tau_T \subseteq \tau_{M,A}$)
\endroster}
\sn
\item "{$(b)$}"  $T_{M,A}$ is the set of universal (first order)
sentences satisfied in ${\frak B}_{M,M,A}$, see part (2).
\ermn
2) Assume $M \subseteq C \prec {\frak C}$ and tp$(C,M \cup A)$ is f.s. in $M$
(e.g., $C=M$).  We define ${\frak B} = {\frak B}_{C,M,A}$ as the
$\tau_{M,A}$-model with universe $C$ such that 
$P^{\frak B}_{\varphi(\bar x,\bar a)} = 
\{\bar b \in {}^{\ell g(\bar x)}C:{\frak C} \models 
\varphi[\bar b,\bar a]\}$ for $\varphi(\bar x,\bar y) \in \Bbb
L(\tau_T),\bar a \in {}^{\ell g(\bar y)}(A)$ and such that
$c^{\frak B}_a = a$  for $a \in M$.  If $C=M$ we may omit $C$.
\nl
3) A model ${\frak B}$ of $T_{M,A}$ is 
called quasi-standard \ub{if} $c^{\frak B}_a=a$ for
$a \in M$.
\nl
3A) A model ${\frak B}$ of $T_{M,A}$ is called standard if it 
is ${\frak B}_{C,M,A}$ for
some $C,M \subseteq C \subseteq {\frak C}$ satisfying tp$(C,M \cup A)$ 
is finitely satisfiable in $M$.
\nl
4) Let $T^*_{M,A}$ be the model completion of $T_{M,A}$ (well defined
only if it exists!)
\enddefinition
\bigskip

\demo{\stag{t.4} Observation}  1) If $M \subseteq C$ and
tp$(C,M \cup A)$ is finitely satisfiable 
in $M$, \ub{then} ${\frak B}_{C,M,A}$ is a model of
$T_{M,A}$. \nl
2) If ${\frak B}$ is a model of $T_{M,A}$, \ub{then} ${\frak B}$ is
isomorphic to the standard model ${\frak B} = {\frak B}_{C,M,A}$ of 
$T_{M,A}$ for some $C$. \nl
3) Moreover, if ${\frak B}_1 \subseteq {\frak B}_2$ are models of
$T_{M,A}$ and ${\frak B}_1$ is standard, \ub{then} ${\frak B}_2$ is
(quasi standard and is)
isomorphic over ${\frak B}_1$ to some standard ${\frak B}'_2$
satisfying ${\frak B}_1 \subseteq {\frak B}'_2$. \nl
4) If $A_1 \subseteq A_2,M \subseteq C$ and tp$(C,M \cup A_2)$ is f.s. in
$M$ \ub{then} ${\frak B}_{C,M,A_1}$ is a 
reduct of ${\frak B}_{C,M,A_2}$. \nl
5) If $M \subseteq C_1 \subseteq C_2$ and tp$(C_2,M \cup A)$ is
f.s. in $M$ \ub{then} ${\frak B}_{C_1,M,A}$ is a submodel of 
${\frak B}_{C_2,M,A}$ (and tp$(C_1,M \cup A)$ is finitely satisfiable
in $M$ hence ${\frak B}_{C_1,M,A}$ is well defined).
\enddemo
\bigskip

\demo{Proof}  Easy.
\enddemo
\bigskip

\proclaim{\stag{t.4A} Claim}  Assume $A$ is full over $M$. \nl
1) ${\frak B}_{M,M,A}$ is a model of $T_{M,A}$ with elimination of 
quantifiers; in fact every subset of ${}^m({\frak B}_{M,M,A})$, i.e.,
of ${}^m|M|$ definable in ${\frak B}_{M,M,A}$ by some
first order formula with parameters is definable by an atomic
formula $R(x_0,\dotsc,x_{m-1})$ in this model. \nl
2) If ${\text{\rm tp\/}}(C,A)$ is f.s. in $M$ \ub{then} we can find $M^+$
such that
\mr
\item "{$(a)$}"  $M \cup C \subseteq M^+ \prec {\frak C}$
\sn
\item "{$(b)$}"  ${\text{\rm tp\/}}(M^+,A)$ is f.s. in $M$
\sn
\item "{$(c)$}"  ${\frak B}_{M^+,M,A}$ is an elementary extension of
${\frak B}_{M,M,A}$.
\ermn
3) $T_{M,A}$ has amalgamation and JEP. \nl
4) ${\text{\rm Th\/}}({\frak B}_{M,M,A})$ is the model completion of
$T_{M,A}$ so is equal to $T^*_{M,A}$ (which is well defined). \nl
5) $T^*_{M,A}$ is a dependent (complete first order) theory.
\endproclaim
\bigskip

\demo{Proof}  1)  By Claim \scite{t.1B}(2), Definition \scite{t.3}(1) 
and $A$ being full over $M$. \nl
2) E.g. use an ultrapower ${\frak C}^\kappa/D$ of ${\frak C}$ with
$\kappa \ge |T| + |C| + |A|,D$ a regular filter on $\kappa$ and let
$\bold j$ be the canonical embedding of ${\frak C}$ into ${\frak
C}^\kappa/D$.   So we
can find $f:C \rightarrow M^\kappa/D$ such that $f \cup (\bold j
\restriction A)$ is an elementary mapping, i.e., a $({\frak C},{\frak
C}^\kappa/D)$-elementary embedding, now it should be clear. \nl
3) The JEP is trivial because of the individual constants $c_a(a \in
M)$.  The amalgamation property holds by \scite{t.2} as we can replace
$M_1$ there by any set which is full over $M$.
\nl
4) By parts (1),(2),(3) where we have already proved.
\nl
5) As ${\frak B}_{M,M,A}$ is a model of it and reflects.

That is, assume $\psi(x,\bar y)$ is a formula with the independence
property in
$T^*_{M,A}$, then by part (1) \wilog \, $\psi$ is an atomic relation
hence for some formula $\varphi(x,\bar y,\bar z) \in \Bbb L(\tau_T)$
and $\bar c \in {}^{\ell g(\bar z)}A$, for every $a,\bar b$ from
$M,{\frak C} \models \varphi[a,\bar b,\bar c]$ iff ${\frak B}_{M,M,A}
\models \psi(a,\bar b)$.  

By the choice of $\psi(x,\bar y)$ for every $n < \omega$ there are
$\bar a^n_\ell \in {}^{\ell g(\bar y)}({\frak B}_{M,M,A}) = {}^{\ell
g(\bar y)}(M)$ for $\ell < \omega$ and $b^n_w \in {\frak B}_{M,M,A}$, i.e.,
$b^n_w \in M$ for $w \subseteq \{0,\dotsc,n-1\}$ such that for every
$w \subseteq \{0,\dotsc,n-1\}$ and $\ell < n$ we have ${\frak
B}_{M,M,A} \models \psi[b^n_w,a^n_\ell]^{\text{if}(\ell \in w)}$,
hence ${\frak C} \models \varphi[b^n_w,\bar a^n_\ell,\bar
c]^{\text{if}(\ell \in w)}$.  So $\varphi(x;\bar y,\bar z)$ has the
independence property in $T$.
${{}}$  \hfill$\square_{\scite{t.4A}}$
\enddemo
\bigskip

\demo{\stag{t.5A} Conclusion}  Assume $M \prec {\frak C}$ and $A \subseteq
{\frak C}$.  Then Th$({\frak B}_{M,M,A})$ is a dependent
(complete first order) theory. 
\enddemo
\bigskip

\demo{Proof}  By \scite{t.4}(4) and \scite{t.4A}(5) it is the reduct
of a dependent (complete first order) theory. 
More fully let $A_1$ be
full over $M$ such that $A \subseteq A_1$ and let $\kappa = |A_1| +
|T|$.  We can find a $\kappa^+$-saturated elementary extension ${\frak
B}'$ of ${\frak B}_{M,M,A_1}$ and by \scite{t.3}(1) \wilog \, it is
${\frak B}_{C,M,A_1}$ for some $C$, so $M \subseteq C$ and tp$(C,M
\cup A_1)$ is finitely satisfiable in $M$.  Clearly if Th$({\frak
B}_{M,MA})$ is dependent then so is Th$({\frak B}_{M,M,A_1}) =
\text{\rm Th}({\frak B}')$.  By \scite{t.4A}(1),(5) we are done.
  \hfill$\square_{\scite{t.5A}}$
\enddemo
\bigskip

\definition{\stag{t.6} Definition}  1) For any model ${\frak B}$ (not
necessarily of $T$) and $A \subseteq {\frak B}$ let $\Bbb B^m[A,{\frak B}]$
be the family of subsets of ${}^m A$ of the form $\{\bar a \in {}^m
A:\varphi(\bar x,\bar a) \in p\}$ for some $p \in \bold S^m(A,{\frak B})$. \nl
2) If ${\frak B} \prec {\frak C}$ we may omit ${\frak B}$.
\enddefinition
\bigskip

\remark{Remark}  If ${\frak B} = {\frak C}$  (or just if
${\frak B}$ is $|A|^+$-saturated) \ub{then} $\Bbb B^m[A,{\frak B}] =
\{\{\bar a:{\frak B} \models \varphi[\bar b,\bar a]\}:\varphi(\bar
x,\bar y) \in \Bbb L(\tau_{\frak B})$ and $\bar b \in {}^{\ell g(\bar y)}
{\frak B}\}$. 
\endremark
\bn
\margintag{t.7}\ub{\stag{t.7} Question}:  Assume $M \subseteq A \subseteq {\frak C}$ 
and ${\frak B}$ a
standard model of $T_{M,A}$ and $N = {\frak B} \restriction \tau_T$.
\ub{Then} do we have
\medskip

\hskip10pt $(*)_{T,T_{M,A}} \quad$  for any ultrafilter $D_0$ on 
$\Bbb B[N,N]$, the number of ultrafilters $D_1$ 
\nl

\hskip60pt on $\Bbb B[N,{\frak B}]$
extending it is at most $2^{|T|+|A|}$?
\bigskip

\remark{\stag{t.8} Remark}  1) For complete (first order theories) $T \subseteq
T_1$, the condition $(*)_{T,T_1}$ of \scite{t.7} has affinity to
conditions like ``any model of $T$ has $< 1$ or $\le \aleph_0$ or $< \|M\|$
expansions to a model of $T_1$".  What is the syntactical characterization? \nl
2) When is ${\frak B}_{N,M,A}$ a model of $T^*_{M,A}$?  Assume $T^*$
has elimination of quantifiers does the condition implies it, i.e.,
implies ${\frak B}_{N,M,A} \models T^*_{M,A}$?
\sn

\hskip10pt $\boxdot_{N,M,A} \quad$  every formula over $N \cup A$ which does
not fork over $N$ 
\nl

\hskip60pt is realized in $N$.
\endremark
\bn
\margintag{t.1.99}\ub{\stag{t.1.99} Discussion}:  1) Note that in the proof \scite{t.1} we
use ``$T$ is dependent" just to deduce that the formula
$\varphi(x,y,\bar z)$ is dependent, i.e., for some $n =
n_{\varphi(x,y,\bar z)}$
\mr
\item "{$\circledast$}"  ${\frak C} \models \neg(\exists x_0 y_0,\dotsc,x_{n-1}
y_{n-1}) \dsize \bigwedge_{w \subseteq n} (\exists \bar z) \dsize
\bigwedge_{\ell < n} \varphi(x_\ell,y_\ell,\bar z)^{\text{if}(\ell \in w)}$.
\ermn
We can in the proof use finite $\Delta_1,\Delta_2$ large enough for
$\varphi(x,y,\bar c)$, i.e., such that for a suitable $n$:
\mr
\item "{$\circledast_2$}"  $\Delta_1 = \{(\exists \bar z)(\dsize
\bigwedge_{\ell < n}
\varphi(x_0,y_0,\dotsc,x_{n-1},y_{n-1},\bar z)^{\text{if}(\ell \in w)}:
w \subseteq n\}$.
\nl
In particular we need
\item "{$\circledast_3$}"  there is $\Delta_1$-indiscernible sequence
$\langle (a_\ell,b_\ell):\ell < 2n\rangle$ and $\bar c'$ such that
${\frak C} \models \varphi[a_\ell,b_\ell,\bar c']$ iff $\ell$ is odd
\sn
\item "{$\circledast_4$}"  $\Delta_2 = \{(\exists y_0,
y_2,\dotsc,y_{2n-2})(\wedge q(x_0,y_0,\dotsc,x_{2n-1},y_{2n-1}):q$ is
a complete $\Delta_1$-type of a $\Delta_1$-indiscernible sequence of
pairs of length $2n\}$
\nl
hence
\sn
\item "{$\circledast_5$}"  there is no $\Delta_2$ indiscernible
sequence $\langle(a_{2 \ell},a_{2 \ell +1},b_{2 \ell +1}):\ell <
n\rangle \char 94 \langle (a_{\omega + 2 \ell},a_{\omega +2 \ell+1},
b_{\omega + 2 \ell +1}):\ell < n\rangle$ such that ${\frak C} \models
\varphi[a_{2 \ell+1},b_{2 \ell +1},\bar c]$ for $\ell < n$ and 
$\{\bar a_{2 \ell},a_{2 \ell},b_{2 \ell +1}:\ell < n\} \subseteq M$ and for
each $\ell < n$ for no $b' \in M$ do we have 
$\models \varphi [a_{2 \ell},b',\bar c]$.
\ermn
2) So looking at the proof and \scite{t.1B}(2)
\mr
\item "{$\circledast_6$}"   there is a finite set $\Delta =
\Delta^*_\varphi$ of formulas of the form $\psi(x,\bar z)$ computable
from $\varphi(x,y,\bar z)$ (and $n_\varphi$) such that:
{\roster
\itemitem{ $(a)$ }  if $M,\bar c,D$ are as in \scite{t.1} then for
some $\bar c'$ the set $\psi(M,\bar c')$ belongs to $D$ and is 
included in $\{a \in M$: for
no $b \in M$ do we have $\models \varphi[a,b,\bar c]\}$
\sn
\itemitem{ $(b)$ }   $\{a \in M:(\exists b \in M)(\varphi(a,b,\bar
c)\}$ is a finite union of sets from $\{\psi(M,\bar C):\bar c' \in
{}^{\bar z'}{\frak C}$ and $\psi(x,\bar z') \in \Delta\}$.
\endroster}
\ermn
If in $\circledast_6(b)$ there is a bound $n$ on the size of
the set not depending on $(M,\bar c)$ let $\Delta^*_\varphi =
\{\psi_\ell(\bar x,\bar z_\ell):\ell < n_*\}$ and let $\psi^*(\bar
x,\bar z) = \dsize \bigwedge_{\ell > n} z^n = z^\ell \rightarrow
\psi_\ell(x,\bar z_\ell)$ can serve instead so in $\circledast_6$
\wilog \, $\Delta^*_\varphi = \{\psi^*_\varphi(\bar x,\bar z^*)\}$.
\nl
3) We elaborate; we know that if $\bold I = \{a \in 
{}^m M$: there is $b \in M$ such that ${\frak C} \models
\varphi[a,b,\bar c]\}$ where 
$\varphi = \varphi(x,y,\bar z) \in \Bbb L(\tau_T),\bar c \in 
{\frak C},M \prec {\frak C}$ \ub{then} for some $\psi(x,\bar z') 
\in \Bbb L(\tau_T)$ and $\bar c' \in {}^{\ell g(\bar z')}{\frak C}$ we
have $\bold I = \psi(M,\bar c')$.  Can 
we characterize $\psi$?  Yes, but not so well.
Toward this first let $n(*)$ be minimal such that there are no
$a_\ell,b_\ell,(\ell < n(*),\bar c_\eta,(\eta \in {}^{\eta(*)}2)$
from ${\frak C}$ such that $M \models \varphi(a_\ell,b_\ell,
\bar c_\eta)$ iff $\eta(\ell)=1$.

Let $\psi_n(x_0,y_0,\dotsc,x_{n(*)-1} \char 94 y_{n(*)-1}) = (\exists
\bar z) \dsize \bigwedge_{\ell < n(*)} \varphi(x_\ell,y_\ell,
\bar z)^{\eta(\ell)}$ for $\eta \in {}^{n(*)} 2$ and $\Delta_1 
= \{\psi_\eta(\bar x_0,\bar y_0,\dotsc,\bar x_{n(*)},\bar
y_{n(*)-1})\}$.  Let $\Delta_2$
be the closure of $\Delta_1$ under permuting the variables.

Let $\Delta_{3,k}$ be the set of formulas of the form
$\vartheta(y_{2k(*)};x_0,y_0,x_1,y_1,\dotsc,x_{2k-1},y_{2k-1} -
y_{2k};x_{2k+1},y_{2k+1},\dotsc,x_{2n(*)-2},x_{2n(*)-1},y_{2n(*)-1}) =
(\exists y_{2k+2}) \ldots (\exists y_{2n(*)_2}) \psi^*$ where $\psi^*$
is a conjunction or formula from $\Delta_2$ and their negation.

Now $\psi$ belongs to $\Delta_{3,k}$ for some $k < n(*)$.  (In fact we
could be somewhat more specific).
\nl
Why?  We work with $\cup \Delta_{3,\ell}$ choose $a_{2 \ell},a_{2
\ell+1},k_{2 \ell +1}$ as in the proof for it.  Then we choose $b_{2
\ell+1} \in M$ by induction on $\ell < n(*)$ such that $\langle
(a_\ell,b_\ell):\ell < 2n(*)\rangle$ is $\Delta_1$-indiscernible.  So
 for every $\eta \in {}^{(*)} 2$ we have $(\exists \bar z) \dsize
\bigwedge_{\ell < n(*)} \varphi(a_\ell,b_\ell,\bar z)^{\eta(\ell)}$.

\newpage

\head {\S2 More on indiscernible sequences} \endhead  \resetall \sectno=2
 \spuriousreset
\bigskip

\demo{\stag{np8.1.7} Context}  1) $T$ is a (first order complete)
dependent theory.
\nl
2) ${\frak C}$ is the monster model of $T$.
\enddemo
\bigskip

This section is complimentary to \cite[\S5]{Sh:715} so recall the
definition.
\definition{\stag{np8.1.8} Definition}  Let $\bar{\bold a}^\ell = \langle
\bar a^\ell_t:t \in I_\ell \rangle$ be an indiscernible sequence which is
endless (i.e., $I_\ell$ having no last element) for $\ell = 1,2$. \nl
1) We say that $\bar{\bold a}^1,\bar{\bold a}^2$ are \ub{perpendicular} when:
\mr
\item "{$(*)$}"  \ub{if} $\bar b^\ell_n$ realizes Av$(\{\bar b^k_m$:
we have $m < n \and k \in \{1,2\}$ or we have 
$m = n \and k < \ell\} \cup \bar{\bold a}^1 \cup
\bar{\bold a}^2,\bar{\bold a}^\ell)$ for $\ell =
1,2$ \ub{then} $\bar{\bold b}^1,
\bar{\bold b}^2$ are mutually indiscernible (i.e., each is
indiscernible over the set of elements appearing in the other)
 where $\bar{\bold b}^\ell = \langle \bar b^\ell_n:n < \omega \rangle$ for
$\ell = 1,2$.
\ermn
We define ``$\Delta$-perpendicular" in the obvious way. \nl
2) We say $\bar{\bold a}^1,\bar{\bold a}^2$ are equivalent and write
$\approx$ \ub{if} for every
$A \subseteq {\frak C}$ we have Av$(A,\bar{\bold a}^1) = \text{Av}
(A,\bar{\bold a}^2)$. \nl
3) If $\bar{\bold a}^1 \subseteq A$ we let dual-cf$(\bar{\bold a}^1,A) =
\text{ Min}\{|B|:B \subseteq A$ and no $\bar c \in {}^{\omega >}A$ realizes
Av$(B,\bar{\bold a}^1)\}$; we usually apply this when $A = M$.
\enddefinition
\bigskip

\proclaim{\stag{np8.1} Claim}  Assume
\mr
\item "{$(\alpha)$}"  $\bar{\bold b} = \langle \bar b_t:t \in I_0
\rangle$ is an infinite indiscernible sequence over $A$
\sn
\item "{$(\beta)$}"  $B \subseteq {\frak C}$.
\ermn
\ub{Then} we can find $I_1,J$ and $\bar b_t$ for $t \in I_1 \backslash
I_0$ such that:
\mr
\item "{$(a)$}"  $I_0 \subseteq I_1,I_1 \backslash I_0 \subseteq J
\subseteq I_1$ and $|I_1 \backslash I_0| \le |J| \le |B| + |T|$
\sn
\item "{$(b)$}"  $\bar{\bold b}' = \langle \bar b_t:t \in I_1 \rangle$
is an indiscernible sequence over $A$
\sn
\item "{$(c)$}"  if $I_2$ is a $J$-free extension of $I_1$ (see below) 
and $\bar b_t$ for $t \in I_2 \backslash I_1$ are such that 
$\bar{\bold b}'' = \langle \bar b_t:t \in I_2 \rangle$
is an indiscernible sequence over $A$ \ub{then} 
{\roster
\itemitem{ $\circledast$ }  if $n < \omega,\bar s,\bar t \in
{}^n(I_2)$ and $\bar s \sim_J \bar t$ (see below), 
\ub{then} $\bar b_{\bar s},
\bar b_{\bar t}$ realize the same type over $A \cup B$
where $\bar b_{\langle t_\ell:\ell < n \rangle} = \bar b_{t_0} \char 94 \bar
b_{t_1} \char 94 \ldots \char 94 \bar b_{t_{n-1}}$.
\endroster}
\endroster
\endproclaim
\bigskip

\definition{\stag{np8.1A} Definition}  1) For linear orders $J,I_1,I_2$
we say that $I_2$ is a $J$-free extension of $I_1$ \ub{when}: 
$J \subseteq I_1 \subseteq I_2$ and
\mr
\item "{$\circledast$}"  if $t \in I_2 \backslash I_1$ and $s \in J$
then for some $t' \in I_1$ we have $I_2 \models s < t' < t$ or 
$I_2 \models t < t' < s$.
\ermn
2) For linear orders $J,I_1,I_2$ we say that $I_2$ is a strong
$J$-free extension of $I_1$ when $J \subseteq I_1 \subseteq I_2$ and:
\mr
\item "{$\circledast$}"  if $t \in I_2 \backslash I_1$ then for some
$s_1,s_2 \in I_1$ we have $s_1 <_{I_2} t <_{I_2} s_2$ and
$[s_1,s_2]_{I_1} \cap J = \emptyset$.
\ermn
3) For linear orders $J \subseteq I$ and $\bar s,t \in {}^n I$ let
$\bar s \sim_J \bar t$ mean that $(s_\ell <_I s_k) 
\equiv (t_\ell <_I t_k)$ and $(s_\ell <_I r) \equiv
(t_\ell <_I r)$ and $(r <_I s_\ell \equiv r <_I t_\ell)$
whenever $\ell,k < n$ and $r \in J$).  
Similarly for $\bar s,\bar t \in {}^\alpha I$.
\enddefinition
\bigskip

\remark{\stag{np8.1B} Remark}  In \scite{np8.1} why do we 
need ``$J$-free"?  Let $M =
(\Bbb R,<,Q^M),Q^M = \Bbb Q,B = \{0\},A = \emptyset,I_0$ the
irrationals, $b_t=t$ for $t \in I_0$.
\endremark
\bigskip

\demo{Proof}  

Proof of \scite{np8.1}.

We try to choose by induction on $\zeta < \lambda^+$
where $\lambda = |T| + |B|$ a sequence $\bar{\bold b}^\zeta = 
\langle \bar b_t:t
\in J_\zeta \rangle$ and together with $\bar{\bold b}^{\zeta +1}$
we choose $n_\zeta,\bar s_\zeta,\bar t_\zeta,
J'_\zeta,\varphi_\zeta,\bar c_\zeta,\bar d_\zeta$ such that
\mr
\item "{$(a)$}"   $J_\zeta$ is a linear order, increasing continuous
with $\zeta$
\sn
\item "{$(b)$}"  $J_0 = I_0$ (so $\bar{\bold b}^0 = \bar{\bold b}),
J_{\varepsilon +1} \backslash
J_\varepsilon$ is finite so $|J_\varepsilon \backslash I_0| <
|\varepsilon|^+ + \aleph_0$
\sn
\item "{$(c)$}"  $\bar{\bold b}^\zeta$ is an indiscernible sequence over $A$
\sn
\item "{$(d)$}"  $J'_\zeta \subseteq J_\zeta,J_\zeta = I_0 \cup
J'_\zeta,J'_\zeta$ is increasing continuous with $\zeta$ and $|J'_\zeta|
< |\zeta|^+ + \aleph_0$
\sn
\item "{$(e)$}"  if $\zeta = \varepsilon +1$ \ub{then} $n_\varepsilon <
\omega,\bar s_\varepsilon \in {}^{n_\varepsilon}(J'_\zeta),\bar
t_\varepsilon \in {}^{n_\varepsilon}(J'_\zeta),\varphi_\varepsilon =
\varphi_\varepsilon(\bar x_0,\dotsc,\bar x_{n_\varepsilon-1},\bar
c_\varepsilon,\bar d_\varepsilon),\bar c_\varepsilon \subseteq B,\bar
d_\varepsilon \subseteq A$ and $J'_\zeta = J'_\varepsilon \cup (\bar
s_\varepsilon \char 94 \bar t_\varepsilon)$
\sn
\item "{$(f)$}"  $\bar s_\varepsilon \sim_{J'_\varepsilon} \bar
t_\varepsilon$ and $\models \varphi_\varepsilon
[\bar b_{\bar s_\varepsilon},\bar
c_\varepsilon,\bar d_\varepsilon] \wedge \neg \varphi_\varepsilon[\bar b_{\bar
t_\varepsilon},\bar c_\varepsilon,\bar d_\varepsilon]$
\sn
\item "{$(g)$}"  $J_{\zeta +1}$ is a $J'_\zeta$-free extension of $J_\zeta$. 
\ermn
If we succeed, for some unbounded $w \subseteq \lambda^+$ and
$n_*,\varphi_\xi,\bar c^*$ and $u$ for every 
$\varepsilon \in w$ we have $n_\varepsilon = n_*,\varphi_\varepsilon =
\varphi_*,\bar c_\varepsilon = \bar c^*$ 
and $u = \{\ell < n_*:s_{\varepsilon,\ell} \in J'_\zeta\}$.  
Now let $J^* = \cup\{J'_\zeta:\zeta < \lambda^+\}$, so every 
$J' \subseteq J^*$ of
cardinality $\le \lambda$ is included in $J'_\zeta$ for some $\zeta <
\lambda^+$ and we get contradiction to clause (b) of 
\cite[3.2]{Sh:715}(=3.4t), hence we fail, i.e., we cannot choose for
some $\zeta$.  But we can choose $\bar{\bold b}_\zeta = \langle b_t:t
\in J_\zeta \rangle$, if $\zeta = 0$ by clause (b) and if $\zeta$ is a
limit ordinal by clause (a).  So $\zeta = \varepsilon +1$, we have
chosen $\bar{\bold b} = \langle \bar b_t:t \in J_\zeta \rangle$ but we
cannot choose $J_{\zeta +1},\bar{\bold b}^{\zeta +1},n_\zeta,\bar
s_\zeta,\bar t_\zeta,J'_\zeta,\varphi_\zeta,\bar c_\zeta,\bar d_\zeta$
as required.
Then $\bar{\bold b}^\varepsilon$ is as required.  
\hfill$\square_{\scite{np8.1}}$
\enddemo
\bn
The aim of \scite{np8.2} + \scite{np8.3} below is to show a complement of
\cite[\S5]{Sh:715}; that is, 
in the case of small cofinality what occurs in one cut
is the ``same" as what occurs in others. 
\proclaim{\stag{np8.2} Claim}  Assume
\mr
\item "{$(a)$}"  $\mu \ge |T|$ 
\sn
\item "{$(b)$}"  $I_\ell$ for $\ell < 4$ are pairwise disjoint linear orders
\sn
\item "{$(c)$}"  $I_\ell = \dbcu_{\beta < \mu^+}
I^\beta_\ell,I^\beta_\ell$ (strictly) increasing with $\beta$ and
$|I^\beta_\ell| \le \mu$ for $\ell < 4$
\sn
\item "{$(d)$}"  $(\alpha) \quad 
\ell \in \{0,2\} \Rightarrow I^\beta_\ell$ an end
segment of $I_\ell$ 
\sn
\item "{${{}}$}"  $(\beta) \quad \ell \in 
\{1,3\} \Rightarrow I^\beta_\ell$ is an initial segment of $I_\ell$
\sn
\item "{$(e)$}"  $I = I_0 + I_1 + I_2 + I_3$ and $I^\beta = I^\beta_0
+ I^\beta_1 + I^\beta_2 + I^\beta_3$
\sn
\item "{$(f)$}"  $\langle \bar b_t:t \in I \rangle$ is an
indiscernible sequence.
\ermn
\ub{Then} we can find a limit ordinal 
$\beta(*) < \mu^+$ and $\langle \bar b^*_t:t \in
I \rangle$ such that:
\mr
\item "{$(A)$}"  $\bar b^*_t = \bar b_t$ if $t \in I \backslash
I^{\beta(*)}$
\sn
\item "{$(B)_1$}"  $\langle \bar b^*_t:t \in I \backslash
I^{\beta(*)}_0 \backslash I^{\beta(*)}_1 \rangle$ is an indiscernible sequence
\sn
\item "{$(B)_2$}"  $\langle \bar b^*_t:t \in I \backslash
I^{\beta(*)}_2 \backslash I^{\beta(*)}_3 \rangle$ is an indiscernible sequence
\sn  
\item "{$(C)_1$}"  ${\text{\rm tp\/}}_*(\langle 
\bar b^*_t:t \in I^{\beta(*)}_0 \cup
I^{\beta(*)}_1 \rangle,\cup \{ \bar b^*_t:t \in (I \backslash
I^\beta) \cup I^{\beta(*)}_2 \cup I^{\beta(*)}_3 \cup
(I^{\beta(*)+ \omega}_0 \backslash I^{\beta(*)}_0) \cup (I^{\beta(*)+
\omega}_1 \backslash I^{\beta(*)}_1)\}) \vdash$ 
\nl
${\text{\rm tp\/}}_*(\langle \bar b^*_t:t \in I^{\beta(*)}_0 \cup
I^{\beta(*)}_1 \rangle,\cup \{\bar b^*_t:t \in (I \backslash
I^{\beta(*)}) \cup I^{\beta(*)}_2 \cup I^{\beta(*)}_3\})$ 
\nl
for any $\beta\in [\beta(*) + \omega,\mu^+)$
\sn
\item "{$(C)_2$}"  ${\text{\rm tp\/}}_*(\langle 
\bar b^*_t:t \in I^{\beta(*)}_2 \cup
I^{\beta(*)}_3 \rangle,\cup \{ \bar b^*_t:t \in (I \backslash
I^\beta) \cup I^{\beta(*)}_0 \cup I^{\beta(*)}_1 \cup
(I^{\beta(*)+\omega}_2 \backslash I^{\beta(*)}_2) \cup
(I^{\beta(*)+\omega}_3 \backslash I^{\beta(*)}_3)\}) \vdash$ 
\nl
${\text{\rm tp\/}}(\langle b^*_t:t \in I^{\beta(*)}_2 \cup
I^{\beta(*)}_3 \rangle,\cup \{\bar b^*_t:t \in (I \backslash
I^{\beta(*)}) \cup I^{\beta(*)}_0 \cup I^{\beta(*)}_1\})$ 
\nl
for any $\beta \in [\beta(*) + \omega,\mu^+)$
\sn
\item "{$(D)_1$}"  $\langle \bar b^*_t:t \in I_0 \backslash
I^{\beta(*)}_0 \rangle$ is an indiscernible sequence over
$\cup \{\bar b^*_t:t \in I^{\beta(*)}_0 \cup I_1 \cup I_2 \cup I_3\}$
\sn
\item "{$(D)_2$}"  $\langle \bar b^*_t:t \in (I_1 \backslash
I^{\beta(*)}_1) + (I_2 \backslash I^{\beta(*)}_2) \rangle$ is an
indiscernible sequence over
$\cup \{\bar b^*_t:t \in I_0 \cup I^{\beta(*)}_1 \cup I^{\beta(*)}_2 
\cup I_3\}$
\sn
\item "{$(D)_3$}"  $\langle \bar b^*_t:t \in I_3 \backslash
I^{\beta(*)}_3 \rangle$ is an indiscernible sequence over
$\cup \{\bar b^*_t:t \in I_0 \cup I_1 \cup I_2 \cup I^{\beta(*)}_3\}$.
\endroster
\endproclaim
\bigskip

\remark{\stag{np8.1.Z} Remark}  What occurs if $T$ is stable (or just
$\bar{\bold b}$ is)?  We get something like $\{\bar b^*_t:t \in
I^{\beta(*)}_0 \cup I^{\beta(*)}_1\} = \{\bar b^*_t:t \in
I^{\beta(*)}_2 \cup I^{\beta(*)}_3\}$.
\endremark
\bigskip

\demo{Proof}  For simplicity assume $I^0_\ell = \emptyset$.

We choose by induction on $n < \omega$ an ordinal $\beta(n)$ and
$\langle \bar b^n_t:t \in I \rangle$ such that:
\mr
\item "{$(\alpha)$}"   $\beta(n) < \mu^+,\beta(0) = 0,\beta(n) +
\omega \le \beta(n+1)$
\sn
\item "{$(\beta)$}"   $\bar b^n_t = \bar b_t$ if $t \in I \backslash
I^{\beta(n)}$ or if $n=0$
\sn
\item "{$(\gamma)_1$}"   $\langle \bar b^n_t:t \in I \backslash I^{\beta(n)}_0
 \backslash I^{\beta(n)}_1 \rangle$ realizes the same type as
$\langle \bar b_t:
t \in I \backslash I^{\beta(n)}_0 \backslash I^{\beta(n)}_1\rangle$
\sn
\item "{$(\gamma)_2$}"  $\langle \bar b^n_t:t \in I \backslash
I^{\beta(n)}_2 \backslash I^{\beta(n)}_3\rangle$ realizes the same
type as $\langle \bar b_t:
t \in I \backslash I^{\beta(n)}_2 \backslash I^{\beta(n)}_3 \rangle$
\sn
\item "{$(\delta)_1$}"   if $n$ is even then: \nl

$(1) \quad \bar b^{n+1}_t = \bar b^n_t$ for $t \in I \backslash
I^{\beta(n)}_2 \backslash I^{\beta(n)}_3$ \nl

$(2) \quad$ if $\beta(n+1) < \beta < \mu^+$ then the type which
$\langle \bar b^{n+1}_t:t \in I^{\beta(n)}_2 \cup I^{\beta(n)}_3 \rangle$
\nl

\hskip35pt realizes over 
$\cup\{\bar b^n_t:t \in (I_0 \backslash I^\beta_0) \cup
I^{\beta(n+1)}_0 \cup (I_1 \backslash I^\beta_1)$ \nl

\hskip35pt   $\cup I^{\beta(n+1)}_1
\cup (I_2 \backslash I^\beta_2) \cup (I^{\beta(n+1)}_2 \backslash
I^{\beta(n)}_2) \cup (I_3 \backslash I^\beta_3) \cup (I^{\beta(n+1)}_3
\backslash I^{\beta(n)}_3)\}$  \nl

\hskip35pt  has a unique extension over \nl

\hskip35pt $\cup\{\bar b^n_t:t \in I
\backslash I^{\beta(n)}_2 \backslash I^{\beta(n)}_3\}$ \nl

$(3) \quad \bar b^{n+1}_t = b^n_t$ if $t \in I^{\beta(k)}_2 \cup
I^{\beta(k)}_3,k < n$
\sn
\item "{$(\delta)_2$}"   if $n$ is odd
like $(\delta_1)$ inverting the roles of
$(I_0,I_1),(I_2,I_3)$
\sn
\item "{$(\varepsilon)$}"   $\langle \bar b^n_t:t \in I \rangle$
satisfies clauses $(D)_1, (D)_2, (D)_3$ of the claim with
$\beta(n)$ instead of $\beta(*)$.
\ermn
The induction step is as in the proof of \scite{np8.1} (though we use
the finite character for the middle clause (2) of clauses
$(\delta)_1,(\delta)_2$). 

Alternatively, letting $n$ be even we try to choose
$\beta_n(\varepsilon),\bar{\bold b}^{n,\varepsilon} = 
\langle \bar b^{n,\varepsilon}_t:t \in I^{\beta(n)}_2 + I^{\beta(n)}_3 \rangle$ by
induction on $\varepsilon \le \mu^+$ such that:
\mr
\item "{$\bigodot(a)$}"  $\beta_n(\varepsilon) < \mu^+$
\sn
\item "{$(b)$}"  $\beta_n(0) = \beta(n)$
\sn
\item "{$(c)$}"  $\beta_n(\varepsilon)$ is increasing and continuous
\sn
\item "{$(d)$}"  $\zeta < \varepsilon \Rightarrow \text{
tp}(\bar{\bold b}^{n,\varepsilon}),\cup\{b^n_t:t \in (I \backslash
I^\varepsilon_{\beta_n(\varepsilon)}) \cup I_{\beta_n(\zeta)}\}) \vdash$ \nl
tp$(\bar{\bold b}^{n,\zeta},\cup\{\bar b^n_t:t \in (I \backslash
I_{\beta_n(\varepsilon)}) \cup I_{\beta_n(\zeta)}\}$
\sn
\item "{$(e)$}"  if $\varepsilon = \zeta +1$, then
$(\delta)_1(2)$ fails if we let \nl

$\bar b^{n+1}_t = \cases
b^n_t &\text{ \ub{if} } t \in I \backslash I^{\beta(n)}_2 \backslash
I^{\beta(n)}_3 \\
\bar b^{n,\zeta}_t &\text{ \ub{if} } t \in I^{\beta(n)}_2 \cup
I^{\beta(n)}_3 \endcases$.
\ermn
If we succeed to carry the induction, by \cite{Sh:715}, for
some $\varepsilon$ the sequences $\langle \bar b^n_t:t \in
I^{\beta_n(\varepsilon)}_0 \rangle,\langle \bar b^n_t:t \in
I^{\beta_n(\varepsilon)}_1 + I^{\beta_n(\varepsilon)}_2
\rangle,\langle \bar b^n_t:t \in I^{\beta_n(\varepsilon)}_3 \rangle$
are mutually indiscernible over $\cup\{\bar b^{n,\mu^+}_t:t \in
I^{\beta(n)}_2 + I^{\beta(n)}_3\} \cup \{b^n_t:t \in (I \backslash
I_{\beta_n(\varepsilon)})\}$ (because $\langle \bar b_t:t \in I_0
\backslash I^{\beta_n(\varepsilon)}_0 \rangle,\langle \bar b_t:t \in
(I_1 \backslash I_{\beta_n(\varepsilon)}) + I_2 \backslash
I^{\beta_n(\varepsilon)}_2 \rangle,\langle \bar b_t:t \in I_3
\backslash I^{\beta_n(\varepsilon)}_3 \rangle$ are mutually
indiscernible, recalling $(\beta)$.

This contradicts $(e)$.  So we cannot complete the induction.  We
certainly succeed for $\varepsilon=0$, and there is no problem for
limit $\varepsilon \le \mu^+$.  So for some $\varepsilon = \zeta +1$ we
have succeed for $\zeta$ and cannot choose for $\varepsilon$.  We
define $\bar b^{n+1}_i$ as in $(e)$ of $\bigodot$ above, and choose
$\beta(n+1) \in [\beta_n(\varepsilon),\mu^+)$ such that clauses
$(\varepsilon)$ holds.

Let $\beta(*) = \cup\{\beta(n):n < \omega\} < \mu^+,\bar b^*_t$ is
$\bar b^n_t$ for every $n$ large enough (exists by clause $(\beta)$ if
$t \in I \backslash I^{\beta(*)}$ and by $(\delta)_\ell(1) + (3)$ if
$t \in I^{\beta(*)}$). \nl
Clearly we are done. \hfill$\square_{\scite{np8.2}}$
\enddemo
\bigskip

\proclaim{\stag{np8.2A} Claim}  Assume
\mr
\item "{$(a)$}"  $I,I^\beta,I_\ell,I^\beta_\ell$ for $\ell < 4,\beta <
\mu^+$ are as in the assumption of claim \scite{np8.2}
\sn
\item "{$(b)$}"  $\beta(*)$ and $\langle \bar b^*_t:t \in I \rangle$ are as
in the conclusion of claim \scite{np8.2}
\sn
\item "{$(c)$}"  $J^+ = J^+_0 + J^+_1 + J^+_2 + J^+_3 + J^+_4$ linear
orders
\sn
\item "{$(d)$}"  $J = J_0 + J_1 + J_2 + J_3 + J_4$ linear orders
\sn
\item "{$(e)$}"  $J_1 = J^+_1 + I^{\beta(*)}_0 + I^{\beta(*)}_1$ and
$J_3 = I^{\beta(*)}_2 + I^{\beta(*)}_3$
\sn
\item "{$(f)$}"  $J_0 \subseteq J^+_0$ and $I_0 \backslash
I^{\beta(*)}_0 \subseteq J^+_0$
\sn
\item "{$(g)$}"  $J_2 \subseteq J^+_2$ and $(I_1 \backslash
I^{\beta(*)}_1) + (I_2 \backslash I^{\beta(*)}_2) \subseteq J_2$
\sn
\item "{$(h)$}"  $J_4 \subseteq J^+_4$ and $(I_3 \backslash
I^{\beta(*)}_3) \subseteq J^+_4$
\sn
\item "{$(i)$}"  $\langle \bar b^*_t:t \in J^+ \rangle$ is an
indiscernible sequence.
\ermn
1) If $J'_0,J'_2,J'_4$ are infinite initial segments of $J_0,J_2,J_4$
respectively \ub{then}
\mr
\item "{$(\alpha)$}"  ${\text{\rm tp\/}}
(\langle \bar b^*_t:t \in J_3 \rangle,\cup
\{\bar b_s:s \in J'_0 \cup J_1 \cup J'_2 \cup J'_4) \vdash 
{ \text{\rm tp\/}}
(\langle \bar b^*_t:t \in J_3 \rangle, \cup \{\bar b^*_s:s \in J_0
\cup J_1 \cup J_2 \cup J_4)$
\sn
\item "{$(\beta)$}"  like $(\alpha)$ interchanging $J_3,J_1$.
\ermn
2) If $J_0$ has no first element, $J'_0 \subseteq J_0$ is unbounded
from below, $J'_2 \subseteq J_2$ is infinite and $J_4$ has no last
element and $J'_4 \subseteq J_4$ is unbounded from above, \ub{then}
the conclusions of (1) holds
\mr
\item "{$(\alpha)$}"  ${\text{\rm tp\/}}
(\langle \bar b^*_t:t \in J_3 \rangle, \bigcup
\{\bar b_s:s \in J'_0 \cup J_1 \cup J'_2 \cup J'_4) \vdash { \text{\rm tp\/}}
(\langle \bar b^*_t:t \in J_3 \rangle,\cup \{\bar b^*_s:s \in J_0
\cup J_1 \cup J_2 \cup J_4\})$
\sn
\item "{$(\beta)$}"  ${\text{\rm tp\/}}
(\langle \bar b^*_t:t \in J_1 \rangle,\bigcup
\{\bar b_s:s \in J'_0 \cup J'_2 \cup J_3 \cup J'_4) \vdash { \text{\rm tp\/}}
(\langle \bar b^*_t:t \in J_1 \rangle,\cup \{\bar b^*_s:s \in J_0
\cup J_2 \cup J_3 \cup J_4\})$.
\ermn
3) If $J^*_0,J^*_2,J^*_4$ has neither first element nor last element
and $J'_0,J'_2,J'_4$ are subsets of $J_0,J_2,J_4$ respectively
unbounded from below and $J''_0,J''_2,J''_4$ are subsets of
$J_0,J_2,J_4$ respectively unbounded from above, \ub{then} the
conclusion of part (1) holds.
\endproclaim
\bigskip

\demo{Proof}  The result follows by the local character of $\vdash$ and by the
indiscernibility demands in \scite{np8.2}, i.e., clauses $(D)_1,(D)_2,(D)_3$.
\hfill$\square_{\scite{np8.2A}}$
\enddemo
\bigskip

\demo{\stag{np8.3} Conclusion}  1) If $\mu \ge \kappa \ge |T|$, 
\ub{then} for some linear order $J^*$ of cardinality $\kappa$ we have
\mr
\item "{$\boxtimes_{\bar{\bold b}^*,J^*}$}"   we can find 
$\bar b_t \in {}^m M$ for $t \in J_3$ such that
$\langle \bar b_t:t \in J \backslash J_1 \rangle$ 
is an indiscernible sequence \underbar{when}:
\sn
\item "{${{}}$}"  $(a) \quad J = J_0 + J_1 + J_2 + J_3 + J_4$
\sn
\item "{${{}}$}"  $(b) \quad$ the cofinalities of
$J_0,J_2,J_4$ and their inverse are $\le \mu$ \nl

\hskip40 pt but are infinite
\sn
\item "{${{}}$}"  $(c) \quad J_1 \cong J^*$ and $J_3 \cong J^*$
(hence $J_1,J_3$ have cardinality $\le \kappa$)
\sn
\item "{${{}}$}"  $(d) \quad \langle \bar b_t:t \in J \backslash
J_3 \rangle$ is an indiscernible sequence (of $m$-tuples)
\sn
\item "{${{}}$}"  $(e) \quad M$ is a $\mu^+$-saturated model
\sn
\item "{${{}}$}"  $(f) \quad \cup\{\bar b_t:t \in J \backslash
J_3\} \subseteq M$.
\ermn
2) If we allow $J^*$ to depend on tp$'(\bar{\bold b}^*)$, see
 Definition \scite{0.istP}(1),  \ub{then} we
can use $J^*$ of the form $\delta^* + \delta,\delta < \gamma$

($\delta^*$ - the inverse of $\delta$).
\enddemo
\bigskip

\demo{Proof}  Let $\bar{\bold b}^*$ be an infinite
indiscernible sequence.

Let $J_0,J_2,J_4$ be disjoint linear orders as in (b).  
Apply \scite{np8.2} with $I_1,I_3$ isomorphic to $(\mu^+,<)$ and
$I_0,I_2$ isomorphic to $(\mu^+,>)$, say $I_\ell =
\{t^\ell_\alpha:\alpha < \mu^+\}$ with $t^\ell_\alpha$ increasing with
$\alpha$ if $\ell \in \{1,3\}$ and decreasing with $\alpha$ if $\ell
\in \{0,2\}$, we get $\bar{\bold b}^* = \langle b^*_t:t \in 
\underset {\ell < 4} {}\to \Sigma I_\ell \rangle,\beta(*)$ 
as there with tp$'(\bar{\bold b}^*
\restriction I_0) = \text{ tp}'(\bar{\bold b}^\circledast)$, see
Definition \scite{0.istP}.  Let $J^+_0 = J_0 + (I_0 \backslash
I^{\beta(*)}_0),J^+_1 = J_1 = I^{\beta(*)}_0 + I^{\beta(*)}_1,
J^+_2 = J_2 + (I_1 \backslash I^{\beta(*)}_1) + 
(I_2 \backslash I^{\beta(*)}_2),
J^+_3 = I^{\beta(*)}_2 + I^{\beta(*)}_3 + J_3$ and
$J^+_4 = J_4 + (I_3 \backslash I^{\beta(*)}_3)$ and
$J^+ = J^+_0 + J^+_1 + J^+_2 + J^+_3 + J^+_4$.  
All $J_\ell$ are infinite linear
orders, choose $J^* = J_1$, clearly $J_3 \cong J^*$.  Now
\mr
\item "{$(*)$}"  $\langle \bar b^*_t:t \in J \backslash J_3 \rangle$
is an indiscernible sequence and
\sn
\item "{$(**)$}"  if $M \supseteq \cup \{\bar b^*_t:
t \in J \backslash J_3\}$ is $\mu^+$-saturated then we can find
$\bar b'_t \in {}^m M$ for $t \in J_3$ such that
$$
\langle \bar b^*_t:t \in J_0 \rangle \char 94
\langle \bar b'_t:t \in J_3 \rangle \char 94
\langle \bar b^*_t:t \in J_4 \rangle
$$
is an indiscernible sequence.
\nl
[Why?  Choose $J'_0 \subseteq J_0$ unbounded from below of cardinality
cf$(J_0,>_{J_0})$ which is $\le \mu$ but $\ge
\aleph_0$ and similarly $J'_2 \subseteq J_2,J'_4 \subseteq J_4$ and
choose $J''_0 \subseteq J_0$ unbounded from above of cardinality
cf$(J_0)$ which is $\le \mu$ and similarly $J''_2 \subseteq J_2,J''_4
\subseteq J_4$ (all O.K. by clause (b) of the assumption). \nl
Now $p = \text{ tp}(\langle b^*_t:t \in J_3 \rangle,\cup\{\bar b_s:s
\in J'_0 \cup J''_0 \cup J'_2 \cup J''_2 \cup J'_4 \cup J'_4\})$ is a
type of cardinality $\le |T| + |J'_0| + |J''_0| + |J'_2| + |J''_2| +
|J'_4| + |J''_4| \le \mu$ hence is realized by some sequence $\langle
\bar b'_t:t \in J_3 \rangle$ from $M$. \nl
By Claim \scite{np8.2A} the desired conclusion in $(**)$ holds.]
\ermn
So we have gotten the desired conclusion for any $\langle J_\ell:\ell \le 4
\rangle$ and indiscernible sequence, $\bar{\bold b} = \langle \bar
b_t:t \in J \backslash J_5 \rangle$ as long as tp$'(\bar{\bold b}) =
\text{ tp}'(\bar{\bold b}^*)$ and the order type of $J_1,J_3$ is as
required for $\bar{\bold b}^*)$.  This is enough for part
(2), we are left with (1).

Note that by the proof of \scite{np8.1}, the set of $\beta(*)$ as
required contains $E \cap \{\delta < \mu^+:\text{cf}(\delta) =
\aleph_0\}$ for some club $E$ (in fact even contains $E$).  So if $\mu
\ge 2^{|T|}$, as $\{\text{tp}'(\bar{\bold b}):\bar{\bold b}$ an
infinite indiscernible sequence$\}$ has cardinality $\le 2^{|T|}$ we
are done.

Otherwise choose $J^*$ a linear order of cardinality $\mu$ isomorphic
to its inverse, to $J^* \times \omega$ and to $J^* \times (\gamma+1)$
ordered lexicographically for every $\gamma \le \mu$ hence for every
$\gamma < \mu^+$, (e.g. note if $J^{**}$ is dense with no first and last
element and saturated (or special) of cardinality $> \mu$, then
$J^{**} \times \omega$ satisfies this and use the L.S. argument).  So
we can in \scite{np8.2} hence \scite{np8.2A} use 
$I_\ell(\ell < 4)$ such that $I^{\beta
+1}_\ell \cong J^*$ for $\beta < \mu^+,\ell < 4$.  So $I^{\beta(*)}_0
+ I^{\beta(*)}_1 \cong J^* \cong I^{\beta(*)}_2 + I^{\beta(*)}_3$.
\hfill$\square_{\scite{np8.3}}$     
\enddemo
\bigskip

\demo{\stag{np8.4} Conclusion}  In \scite{np8.3}:
\mr
\item "{$(A)$}"  we can choose $J^* = \mu^* + \mu$ i.e. $\{0\}
\times (\mu,>) + \{1\} \times (\mu,<)$
\sn
\item "{$(B)$}"  if $J$ is a linear order $(\ne \emptyset)$ of
cardinality $\le \mu$, we can use $J^* = (\mu^* + \mu) \times J$
ordered lexicographically
\sn
\item "{$(C)$}"   we can change the conclusion of \scite{np8.3} 
to make it symmetrical between $J_3$ and $J_1$
\sn
\item "{$(D)$}"  we use only clause $(E)_2$ of \scite{np8.2} or we could
use only clause $(E)_1$. 
\endroster
\enddemo
\bigskip

\demo{Proof}  (A),(B) combine the proofs of \scite{np8.1} and
\scite{np8.2} trying to contradict each formula, by bookkeeping trying
for it enough times.
 \hfill$\square_{\scite{np8.4}}$
\enddemo
\bn
We may look at it differently, part (2) is close in formulation to be
a complement to \cite[\S5]{Sh:715}.
\demo{\stag{np8.5} Conclusion}  1) Assume
\mr
\item "{$(a)$}"  $J = I \times J^*$ is ordered 
lexicographically, $J^*,\mu$ are as in \scite{np8.3}, $I$ infinite
\sn
\item "{$(b)$}"  $\langle \bar b_t:t \in J \rangle$ an indiscernible
sequence, $\ell g(\bar b_t) = m$ or just $\ell g(\bar b_t) < \mu^+$
\sn
\item "{$(c)$}"  for $s \in I$ let $\bar c_s$ be $\langle \bar b_t:t
\in \{s\} \times J^* \rangle$, more exactly the concatanation of the
sequences in $\bar b_t$ for $t \in \{s\} \times J^*$.
\ermn
\ub{Then}
\mr
\item "{$(\alpha)$}"  $\langle \bar c_s:s \in I \rangle$ is an
infinite indiscernible sequence
\sn
\item "{$(\beta)$}"  if $s_0 <_I \ldots <_I s_7$ then there is $\bar
c$ realizing tp$(\bar c_{s_2},\bigcup\{\bar c_{s_\ell}:\ell \le 7,\ell
\ne 2\})$ such that tp$(\bar c,\cup\{\bar c_{s_\ell}:\ell \le 7,\ell
\ne 2\}) \vdash \text{ tp}(\bar c_{s_2},\bigcup\{\bar c_s:s_0 \le_I s \le_I
s_1$ or $s_3 \le_I s \le_I s_4$ or $s_6 \le_I s \le_I s_7\})$
\sn
\item "{$(\gamma)$}"  similarly inverting the order
(i.e. interchanging the roles of $s_2,s_5$ in clause $(\beta))$.
\ermn
2) Assume the sequence $\langle \bar c_s:s \in I \rangle$ from part (1)
satisfies $M \supseteq \bigcup\{\bar c_s:s \in I\}$ and
$(I_1,I_2),(I_3,I_4)$ are Dedekind cuts of $I$, each of
$I_1,(I_2)^*,I_3,(I_4)^*$ is non-empty of cofinality $\le \mu$ and let
$I^+ \supseteq I,t_2,t_5 \in I^+_1$ realize the cuts
$(I_1,I_2),(I_3,I_4)$ respectively and $\bar c_t$ for $t \in I^+
\backslash I$ are such that $\langle \bar c_t:t \in I^+ \rangle$ is
indiscernible (then for notational simplicity). \ub{Then}
\mr
\item "{$\boxdot$}"  there is a sequence in $M$ realizing tp$(\bar
c_{t_2},\bigcup \{\bar c_s:s \in I\})$ \ub{iff} there is a sequence in $M$
realizing tp$(\bar c_{t_5},\bigcup \{\bar c_s:s \in I\})$.
\endroster
\enddemo
\bigskip

\remark{Concluding Remark}  There is a gap between
\cite[5.11=np5.5]{Sh:715} and the results in \S2,
some light is thrown by
\endremark
\bigskip

\proclaim{\stag{npa.8.6} Claim}  In 
\cite[5.11=np5.5]{Sh:715}; we can omit the demand
${\text{\rm cf\/}}({\text{\rm Dom\/}}(\bar{\bold a}^\zeta)) \ge
\kappa_1$ (= clause (f) there) if we add $\zeta < \zeta^* \Rightarrow
(\theta^1_\zeta)^+ = \lambda$.
\endproclaim
\bigskip

\demo{Proof}  By the omitting type argument.
\enddemo
\bn
\margintag{npa.8.7}\ub{\stag{npa.8.7} Question}:  Assume:
\mr
\item "{$(a)$}"  $\langle (N_i,M_i):i \le \kappa \rangle$ is
$\prec$-increasing (as pairs), $M_{i+1},N_{i+1}$ are
$\lambda^+_i$-saturated, $\|N_i\| \le \lambda_i,\langle \lambda_i:i <
\kappa \rangle$ increasing, $\kappa < \lambda_0$
\sn
\item "{$(b)$}"  $p(\bar x)$ is a partial type over $N_0 \cup
M_\kappa$ of cardinality $\le \lambda_0$.
\ermn
1) Does $p(\bar x)$ have a $\lambda^+_0$-isolated extension?
\nl
2) Does this help to clarify DOP?
\nl
3) Does this help to clarify ``if any $M$ is a benign set" (see
\cite{BBSh:815}).
\bigskip

\proclaim{\stag{npa.8.8} Claim}  Assume
\mr
\item "{$(a)$}"  $M$ is $\lambda^+$-saturated
\sn
\item "{$(b)$}"  $p(\bar x)$ is a type of cardinality $\le \kappa,\ell
g(\bar x) \le \kappa$
\sn
\item "{$(c)$}"  ${\text{\rm Dom\/}}(p) \subseteq A \cup M,|A| \le
\kappa \le \lambda$
\sn
\item "{$(d)$}"  $B \subseteq M,|B| \le \lambda$.
\ermn
\ub{Then} there is a type $q(\bar x)$ over $A \cup M$ of cardinality
$< \kappa$ and $r(\bar x) \in \bold S^{\ell g(\bar x)}(A \cup B)$ such
that

$$
p(\bar x) \subseteq q(\bar x)
$$

$$
q(\bar x) \vdash r(\bar x)
$$
\endproclaim
\bigskip

\remark{Remark}  This defines a natural quasi order (type definable)
is it directed?
\endremark
\newpage

\head{\S3 strongly dependent theories} \endhead  \resetall \sectno=3
 \spuriousreset
\bigskip

\demo{\stag{ss.0} Context}  $T$ complete first order, ${\frak C}$ a
monster model of $T$.
\enddemo
\bigskip

\definition{\stag{ss.1} Definition}  1) $T$ is strongly$^1$ dependent
(we may omit the 1) \ub{if}: \nl
there are no $\bar \varphi = \langle \varphi_n(\bar x,\bar y_n):n <
\omega \rangle$ and $\langle \bar a^n_\alpha:n < \omega,\alpha <
\lambda \rangle$ such that
\mr
\item "{$(*)$}"  for every $\eta \in {}^\omega \lambda$ the set
$p_\eta = \{\varphi_n(\bar x,\bar a^n_\alpha)^{\text{if}(\eta(n)=\alpha)}:
\alpha < \lambda\}$ is
consistent; so $\ell g(\bar a^n_\alpha) = \ell g(\bar y_n)$.
\ermn
2) $T$ is strongly stable if it is stable and strongly dependent.
\nl
3) $\kappa_{\text{ict}}(T)$ is the first $\kappa$ such that there is no $\bar
\varphi = \langle \varphi_\alpha(\bar x,\bar y_\alpha):\alpha < \kappa
\rangle$ satisfying the parallel of part (1), in this case we say that $\bar
\varphi$ witnesses $\kappa < \kappa_{\text{ict}}(T)$ and let $m(\bar
\varphi) = \ell g(\bar x)$.
\enddefinition
\bigskip

\proclaim{\stag{ss.2} Claim}  1) If $T$ is superstable, \ub{then} $T$
is strongly dependent.
\nl
2) If $T$ is strongly dependent, \ub{then} $T$ is dependent.
\nl
3) There are stable $T$ which are not strongly dependent. \nl
4) There are stable not superstable $T$ which are strongly
dependent. \nl
5) There are unstable strongly dependent theories. \nl
6) The theory of real closed fields is strongly dependent; moreover
every o-minimal (complete first-order) $T$ is strongly
dependent.
\nl
7)  If $T$ is stable \ub{then} $\kappa_{\text{ict}}(T) \le \kappa(T)$.
\nl
8) If $T$ is dependent \ub{then} we may add, in \scite{ss.1}(1)
\mr
\item "{$(**)$}"  for each $n < \omega$ for some $k_n$ any $k_n$ of the
formulas $\{\varphi_n(\bar x,\bar a^n_\alpha):\alpha < \lambda\}$ 
are contradictory.
\endroster
\endproclaim
\bigskip

\demo{Proof}  
1),2),7),8)  Easy.  
\nl
3) E.g. $T = \text{ Th}({}^\omega \omega,E^1_n)_{n < \omega}$ where $\eta
E^1_n \nu \Leftrightarrow \eta(n) = \nu(n)$ and use $\varphi_n(x,y_n)
= x E^1_n y_n$ for $n < \omega$. \nl
4) E.g.,  $T = \text {Th}({}^\omega \omega,
E^2_n)_{n < \omega}$ where $(\eta E_n \nu)
\Leftrightarrow (\eta \restriction n = \nu \restriction n)$. \nl
5) E.g., $T = \text{ Th}(\Bbb Q,<)$, the theory of dense linear orders with no
first and no last element.
\nl
6) For simplicity we use $\bar x = \langle x\rangle$, (justified in
\cite[Observation,1.7]{Sh:863}(1)). 
Assume $\langle \varphi_n(x,\bar y_n):n < \omega \rangle$ and
$\langle \bar a^n_\alpha:\alpha < \lambda \rangle$ are as in
Definition \scite{ss.1}.  Clearly we can replace $\varphi_n(x,\bar
y_n),\bar a^n_\alpha$ by $\varphi'(x,\bar y'_n),\bar b^n_\alpha$ when
$\bar y_n \trianglelefteq \bar y'_n,\bar a^n_\alpha \trianglelefteq
\bar b^n_\alpha$ and $\varphi_n(x,\bar a^n_\alpha) \equiv
\varphi'_n(x,\bar b^n_\alpha)$.  Also we can restrict ourselves to
$\langle \bar a^n_\alpha:n < \omega,\alpha \in u_n\rangle$ where $u_n
\subseteq \lambda$ is infinite for $n <  \omega$.  Hence by 
the elimination of quantifiers and density of the linear order \wilog \,
$\varphi_n(x,\bar y_n) = (\varphi_{\eta,1}(x,\bar y_n) \vee
\varphi_{n,2}(x,\bar y_n)) \wedge \varphi_{n,3}(\bar y_n)$ where
(\wilog \, $\bar y_n = \langle y_\ell:\ell=0,\dotsc,k(n)\rangle,u(n,1)
\subseteq \{0,\dotsc,k(n)-1\},u(n,2) \subseteq \{0,\dotsc,k(n)-1\}$

$$
\varphi_{n,1}(x,\bar y_n) = \dsize \bigvee_{\ell \in u(n,1)}
(y_{n,\ell} < x < y_{n,\ell +1})
$$

$$
\varphi_{n,2}(x,\bar y_n) = \dsize \bigvee_{\ell \in u(n,2)} x = y_{n,\ell}
$$
\mn
and

$$
\varphi_{n,3}(\bar y) = \dsize \bigwedge_{\ell < k(n)} y_{n,\ell} <
y_{n,\ell +1}.
$$
\mn
As $\varphi_{n,3}(\bar a^n_\alpha)$ is satisfied for every $n,\alpha$
as said above we can omit $\varphi_{n,3}(\bar y)$.

For each $\eta \in {}^\omega \lambda,p_\eta$ is consistent (and
$\eta \ne \nu \in {}^\omega \lambda \Rightarrow p_\eta,p_\nu$ are
contradictory), hence clearly each $p_\eta$ is not algebraic.  From
this it follows that $(*)$ of \scite{ss.1}(1)
is true also if we replace $\langle \varphi_n(x,\bar y_n):n < \omega
\rangle$ by $\langle \varphi_{n,1}(x,\bar y_n):n < \omega \rangle$.
Also \wilog \, $\langle \bar a^n_\alpha:\alpha < \lambda \rangle$ is
indiscernible over $\bigcup\{a^m_\beta:m \ne n,m < \omega$ and $\beta <
\lambda\}$.  Now for some $\langle \ell_n:n < \omega \rangle \in \dsize
\prod_{n < \omega} k(n)$, we can replace $\langle \varphi_n(x,\bar
y):n < \omega \rangle$ by $\varphi'_n(x,\bar y) = y_{n,\ell_n} < x <
y_{n,\ell_n +1}$.  So \wilog \, $n < \omega \Rightarrow k(n) =1,\ell_n =
0,\bar y_n = (y_{n,0},y_{n,1})$.

Now $\langle \bar a^n_\alpha:\alpha < \lambda \rangle$ is an
indiscernible sequence, and $\varphi_n({\frak C},\bar a^n_\alpha)$
is the open convex set, actually an interval 
which $\bar a^n_\alpha$ define.  So checking by
cases (they are $a^n_{\alpha,0} < a^n_{\alpha,1} < a^n_{\alpha +1,0} <
a^n_{\alpha +1,1},a^n_{\alpha,0} < a^n_{\alpha +1,0} < a^n_{\alpha,1}
< a^n_{\alpha +1,1},a^n_{\alpha +1,0} < a^n_{\alpha +1,1} <
a^n_{\alpha,0} < a^n_{\alpha,1},a^n_{\alpha +1,0} < a^n_{\alpha,0} <
a^n_{\alpha +1,1} < a^n_{\alpha,1}$ letting 
$p^n_\beta(x) := \{\varphi(x,\bar a^n_\alpha)^{\text{if}(\alpha=\beta)}:
\alpha < \lambda\}$ we note that it is a type
such that $p^n_\beta({\frak C})$ is a convex set; obviously it is 
disjoint to $p^n_\gamma({\frak C})$ for 
$\gamma \in \lambda \backslash \{\beta\}$.

Clearly there are 
$\alpha \ne \beta < \lambda$ such that $p^0_\alpha({\frak C}) <
p^0_\beta({\frak C})$ and choose $a^*$ such that $p^0_\alpha({\frak C})
< a^* < p^0_\beta({\frak C})$.  Now for every $\gamma < \lambda$ we
have $p^1_\gamma({\frak C}) \bigcap p^0_\alpha({\frak C}) \ne \emptyset$
and $p^1_\gamma({\frak C}) \bigcap p^0_\beta({\frak C}) \ne
\emptyset$, i.e., $p^1_0({\frak C})$ is disjoint neither to
$p^0_\alpha({\frak C})$ nor to $p^0_\beta({\frak C})$ (by the choice
of $\bar \varphi,\langle \bar a^n_\alpha:n < \omega,\alpha <
\lambda\rangle$).  As $p^1_\gamma({\frak C})$ is convex, by the choice
of $a^*$ necessarily $a^* \in p^1_\gamma({\frak C})$. As $\gamma$ was
any ordinal $< \lambda$ it follows that
$a^* \in \bigcap\{p^1_\gamma({\frak C}):
\gamma < \lambda\}$, clear contradiction. 
(In fact we get
contradiction even if we use only $n=0,1$, see on this
\cite{Sh:863}).  The o-minimal case holds by the same proof.  
     \hfill$\square_{\scite{ss.2}}$
\enddemo
\bigskip

\definition{\stag{ss.3} Definition}  1) We say a pair of types
$(p(x),q(\bar y))$ is a $(1 = \aleph_0)$-pair of types 
(or $(p(\bar x),q(\bar y))$ satisfies $1 = \aleph_0$) \ub{if} there is a set
$A$ such that: for every countable set $B \subseteq p({\frak C})$,
there is an element $\bar a \in q({\frak C})$ satisfying 
$B \subseteq \text{ acl}(\{\bar a\} \cup A)$.  We say $p(x)$ is a $(1=
\aleph_0)$-type if this holds for some $q(\bar y)$.
\nl
1A) If $A = \text{ Dom}(p)$ we add purely.  We call $A$ a witness to
$p(x)$ being a $(1 = \aleph_0)$-type.
\nl
2) We say $T$ is a local 
$(1 = \aleph_0)$-theory \ub{if} for some $A$ (the witness) 
some non-algebraic type $p$ over $A$ is a $(1 = \aleph_0)$-type.  If $A =
\emptyset$ we say purely. 
\nl
2A) We say $T$ is a global $(1 = \aleph_0)$-theory when the type $x=x$
is a $(1 = \aleph_0)$-type.
\nl
3) We say that a pair $(p(x),q(\bar y))$ of types is a semi 
$(1=\aleph_0)$-pair of types \ub{if}: for some set $A$ 
for every indiscernible sequence
$\langle a_n:n < \omega \rangle$ over $A$ satisfying $n < \omega
\Rightarrow \bar a_n \in p({\frak C})$ 
there is $\bar a \in q({\frak C})$ such that $\{\bar a_n:n < \omega\} 
\subseteq \text{ acl}(\bar a \cup A)$.  We say $p(\bar x)$ is semi
$(\aleph_0=1)$-type if this holds for some $q(\bar y)$.
\nl
4) We say that the pair $(p(x),q(\bar y))$ of types is a 
weakly $(1=\aleph_0)$-pair of types
if there are $A \supseteq \text{\rm Dom}(p)$ and an infinite
indiscernible sequence $\langle a_n:n < \omega \rangle$ over $A$ with
each $a_n$ realizing $p$ such that for some $\bar c \in q({\frak C})$
we have $\{a_n:n < \omega\} \subseteq \text{\rm acl}(A \cup \bar c)$.  
\nl
5) We say $p(x)$ is semi/weakly $(\aleph_0=1)$-type if some pair
$(q,p)$ is semi/weakly $(\aleph_0 =1)$-pair of types.
\nl
6) In (3),(4) we let ``purely", ``witness"
``local"; ``global" be defined similarly.
\nl
7) Above we can allow $p=p(\bar x),\ell g(\bar x) = m$.
\enddefinition
\bigskip

\demo{\stag{ss.3A} Observation}  1) Every algebraic type $p(x)$ is 
a $(1 = \aleph_0)$-type.  If $p \subseteq q$ and $p$ is a $(1 =
\aleph_0)$-type \ub{then} $q$ is an $(1 = \aleph_0)$-type.
\nl
2) If $p(x)$ is a $(1 = \aleph_0)$-type \ub{then} $p(x)$ is a semi $(1
   = \aleph_0)$-type.
\nl
3) If $p(x)$ is a semi $(1,\aleph_0)$-type then $p(x)$ is a weakly $(1 =
\aleph_0)$-type.
\nl
4) If $(p(x),q(\bar y))$ is [semi][weakly]-$(1=\aleph_0)$ type in
${\frak C}$ \ub{then} the same holds in ${\frak C}^{\text{eq}}$.
If $p(x)$ is [semi][weakly]-$(1 = \aleph_0)$-type in ${\frak
C}^{\text{eq}}$ such that $p({\frak C}^{\text{eq}}) \subseteq {\frak
C}$ \ub{then} so is the case in ${\frak C}$.  We can also keep track
of the witness.
\nl
5) For some $T,T$ is not locally $(1 = \aleph_0)$-theory but
$T^{\text{eq}}$ is.
\enddemo
\bigskip

\demo{Proof}  Easy.  \hfill$\square_{\scite{ss.3A}}$
\enddemo
\bigskip

\proclaim{\stag{ss.4} Claim}  1) If $T$ is strongly dependent, \ub{then}
no non-algebraic type is a $(1=\aleph_0)$-type.
\nl
2) Moreover no non-algebraic type is a weakly $(1 = \aleph_0)$-type.
\endproclaim
\bigskip

\remark{Remark}  We can weaken the assumption of \scite{ss.4} to:
for some $\omega$-sequence of
non-algebraic types $\langle p_n(x):n < \omega \rangle$ over $A$, for
every $\langle b_n:n < \omega\rangle \in \dsize \prod_{n< \omega}
p_n({\frak C})$, for some
$\bar c$ we have $\{b_n:n < \omega\} \subseteq c \ell(A \cup \bar c)$.
\endremark
\bigskip

\demo{Proof}  Let $\lambda > |T|^+$.  Assume toward contradiction
that $p(x)$ is a non-algebraic $(1 = \aleph_0)$-type and 
$A$ a witness for it.  As
$p(x)$ is not algebraic, we can find $\bar b^n = \langle
b^n_\alpha:\alpha < \lambda \rangle$ for $n < \omega$ such that
\mr
\item "{$(*)_1$}"  $b^n_\alpha$ realizes $p$
\sn
\item "{$(*)_2$}"  $b^n_\alpha \ne b^n_\beta$ for $\alpha < \beta <
\lambda,n < \omega$
\sn
\item "{$(*)_3$}"  $\langle b^n_\alpha:(n,\alpha) \in \omega \times \lambda
\rangle$ is an indiscernible sequence over 
$A$ where $\omega \times \lambda$ is ordered lexicographically.
\ermn
Let $\bar a \in {}^{\omega >}({\frak C})$ be such that 
$\{b^n_0:n < \omega\} \subseteq
\text{ acl}(A \cup \bar a)$ so for each $n$ we can find $k_n < \omega,\bar c_n
\in {}^{\omega >} A$ and a formula $\varphi_n(x,\bar y,\bar z)$ such that
${\frak C} \models \varphi(b^n_0,\bar a,\bar c_n) \and (\exists^{\le k_n}
x)\varphi(x,\bar a,\bar c_n)$.  By omitting some $b^n_\alpha$'s we have
$(n,\alpha) \in \omega \times \lambda \backslash \{(m,\omega):m <
\omega\} \Rightarrow {\frak C} \models \neg \varphi_n[b^n_\alpha,
\bar a,\bar c_n]$.
\nl
Let $\bar a^n_\alpha = \langle b^n_\alpha \rangle \char 94 \bar c_n$ 
and $\varphi_n$ have already been chosen.  

Now check Definition \scite{ss.1}.  \hfill$\square_{\scite{ss.4}}$
\enddemo
\bigskip

\definition{\stag{ss.4.3} Definition}  1) We say $T$ is strongly$^2$ 
(or strongly$^+$) dependent \ub{when}: there is \ub{no} sequence $\langle
\varphi_n(\bar x,\bar y_0,\dotsc,\bar y_n):n< \omega\rangle$ and $\bar
a^n_\alpha \in {}^{\ell g(y_n)}{\frak C}$ for $n < \omega,\alpha <
\lambda$ (any infinite $\lambda$) such that for every $\eta \in
{}^\omega \lambda$ the set $\{\varphi_n(\bar x,\bar
a^0_{\eta(0)},\dotsc,a^{n-1}_{\eta(n-1)},
a^n_\alpha)^{\text{if}(\alpha=\eta(n))}:n < \omega,\alpha < \lambda\}$
is consistent.
\nl
2) Let $\ell \in \{1,2\}$.  
We say that $T$ is strongly$^{\ell,*}$ dependent \ub{when}: if $\langle
\bar{\bold a}_t:t \in I \rangle$ is an indiscernible sequence over
$A,t \in I \Rightarrow \ell g(\bar{\bold a}_t) = \alpha$ (so
constant but not necessarily finite) and $m < \omega$ and $\bar b_n \in
{}^m{\frak C}$ for $n < \omega,\langle \bar b_n:n < \omega \rangle$ is
an indiscernible sequence over $A \cup\{\bar{\bold a}_t:t \in I\}$
\ub{then} we can divide $I$ to finitely many convex sets $\langle
I_m:m < k \rangle$ such that for each 
$m < k,\langle \bar{\bold a}_t:t \in I_m\rangle$ is an 
indiscernible sequence over $\cup\{\bar b_\alpha:\alpha < \omega\}
\cup A \cup \{\bar a_s:s \in I \backslash I_m$ and $\ell=2\}$.
\nl
3) $T$ is strongly$^\ell$ stable (or strongly$^{\ell,*}$ stable)
when it is strongly$^\ell$ dependent (or strongly$^{\ell,*}$
dependent) and stable.
\enddefinition
\bigskip

\proclaim{\stag{ss.4.5} Claim}   If $T$ is strongly$^+$ dependent \ub{then}:
\mr
\item "{$\circledast_1$}"  for any $A \subseteq {\frak C}$, infinite complete
linear order $I$ and indiscernible sequence $\langle \bar a_t:t \in I
\rangle$ over $A,\ell g(\bar a_t)$ possibly infinite, for any finite
$B \subseteq {\frak C}$, there is a finite $w \subseteq I$ such that:
if $J$ is a convex subset of $I$ disjoint to $w$ then $\langle \bar
a_t:t \in J\rangle$ is indiscernible over $A \cup B \cup \{\bar a_s:s
\in I \backslash J\}$
\sn 
\item "{$\circledast_2$}"  for any set $A \subseteq {\frak C}$ of
cardinality $\lambda$ and infinite linear orders $I_\alpha$ for
$\alpha < \lambda$ and $\bar a^\alpha_t (t \in I_\alpha,\alpha <
\lambda)$ such that $\langle \bar a^\alpha_t:t \in I_\alpha \rangle$ is an
indiscernible sequence over $A \cup \{\bar a^\beta_s:\beta \in \lambda
\backslash \{\alpha\},s \in I_\beta\}$ and finite $B \subseteq {\frak
C}$ \ub{there} is a finite $u \subseteq \lambda$ and $w_\alpha \in
[I_\alpha]^{< \aleph_0}$ for $\alpha \in u$ such that: if $\bar J =
\langle J_\alpha:\alpha < \lambda \rangle,J_\alpha$ is a convex subset
of $I_\alpha$, disjoint to $w_\alpha$ when $\alpha \in u$ then $\langle
\bar a^\alpha_t:t \in J_\alpha\rangle$ is indiscernible over $A \cup B
\cup \{\bar a^\beta_s:\beta \in \lambda \backslash \{\alpha\},s \in
J_\beta\}$ for every $\alpha < \lambda$.
\endroster
\endproclaim
\bigskip

\demo{Proof}  See this (and more) \cite[\S2]{Sh:863}.
\hfill$\square_{\scite{ss.4.5}}$ 
\enddemo
\bigskip

\definition{\stag{ss.5} Definition}  1) We say that $\vartheta(x_1,x_2;\bar
c)$ is a finite-to-finite function from $\varphi_1({\frak C},\bar
a_1)$ onto $\varphi_2({\frak C},\bar a_2)$ \ub{when}:
\mr
\item "{$(a)$}"  if $b_2 \in \varphi_2({\frak C},a_2)$ \ub{then} the set
$\{x:\vartheta(x,b_2,\bar c) \wedge \varphi_1(x,\bar a_1)\}$
satisfies:
{\roster
\itemitem{ $(i)$ }  it is finite \ub{but}
\sn
\itemitem{ $(ii)$ }  it is not empty except for finitely many such $b_2$'s
\endroster}
\item "{$(b)$}"  if $b_1 \in \varphi_1({\frak C},\bar a_1)$ \ub{then} the set
$\{x:\vartheta(b_1,x,\bar c) \wedge \varphi_2(x,\bar a_2)\}$
satisfies:
{\roster
\itemitem{ $(i)$ }  it is finite \ub{but}:
\sn
\itemitem{ $(ii)$ }  it is not empty except for finitely many such $b_1$'s.
\endroster}
\ermn
2) If we place ``onto $\varphi_2({\frak C},\bar a_2)$" by ``into
$\varphi_2({\frak C},\bar a_1)$" we mean that we require 
above only clauses $(a)(i),(b)(i),(ii)$.
\nl
3) We can replace $\varphi_1(x,\bar a_1),\varphi_2(x,\bar a_2)$ above
by types.
\enddefinition
\bigskip

\proclaim{\stag{ss.6} Claim}  If $T$ is strongly$^+$ dependent 
\ub{then} the following are impossible:
\mr
\item "{$(St)_1$}"  for some $\varphi(x,\bar a)$
{\roster
\itemitem{ $(a)$ }  $\varphi(x,\bar a)$ is not algebraic
\sn
\itemitem{ $(b)$ }  $E$ is a definable equivalence relation (in
${\frak C}$ by a first order formula possibly with parameters) with
domain $\subseteq \varphi({\frak C},\bar a)$ and infinitely many
equivalence classes
\sn 
\itemitem{ $(c)$ }  there is a formula $\vartheta(x,y,\bar z)$ such that for
every $b \in { \text{\rm Dom\/}}(E)$ for some $\bar c$, the formula
$\vartheta(x,y;\bar c)$ is a
finite to finite map from $\varphi({\frak C},\bar a)$ into $b/E$; 
\endroster} 
\item "{$(St)_2$}"  for some formulas
$\varphi(x),xEy,\vartheta(x,y,\bar z)$ possibly with parameters we have:
{\roster
\itemitem{ $(a)$ }  $\varphi(x)$ is non-algebraic
\sn
\itemitem{ $(b)$ }  $xEy \rightarrow \varphi(x) \wedge \varphi(y)$
\sn
\itemitem{ $(c)$ }  for uncountably many $c \in \varphi({\frak C})$
for some $\bar d$ the formula $\vartheta(x,y;\bar d)$ is a finite 
to finite function from $\varphi(x)$ into $xEc$
\sn
\itemitem{ $(d)$ } for some $k < \omega$, \ub{if} $b_1,\dotsc,b_k \in
\varphi({\frak C})$ are pairwise distinct then $\dsize
\bigwedge^k_{\ell =1} xEb_\ell$ is algebraic
\endroster}
\item "{$(St)_3$}"  similarly with $\varphi(x,\bar a)$ replaced by a
type, as well as $xEy$ replacing uncountable by cardinality of $\bar\kappa$
(and $x,y,z$ are replaced by $m$-tuples).
\endroster
\endproclaim
\bigskip

\demo{Proof}  The proof for $(St)_1$ is a special case of the proof of
$(St)_2$ and the proof for $(St)_3$ is similar but choosing $k^*$ and
$r(x)$ we should use for subtypes of $\varphi(x),x \in
y,\vartheta(x,y,\bar z)$.  So it is enough:
\sn
\ub{Proof of ``$(St)_2$ is impossible"}:

Without loss of generality in clause (c) of $(St)_2$ we have $\langle
c \rangle \triangleleft \bar d$ and let $\ell g(\bar d) = j$, i.e.,
$\vartheta = \vartheta(x,y,\bar z),\ell g(\bar z) = j$; also let $\bar
z^n = \langle z_{n,0},\dotsc,z_{n,j-1}\rangle$.

Clearly there is $k^*$ such that
\mr
\item "{$\boxdot_1$}"  for some uncountably $C \subseteq 
\varphi({\frak C})$ for every $c \in C$ for some $\bar d_c \in
{}^j{\frak C}$ \wilog \, $\langle c \rangle 
\triangleleft \bar d_c$ and $\vartheta(x,y,\bar d_c)$ is a finite
 to finite map from $\varphi({\frak C})$ into $x Ec$
and the size of the finite sets (see Definition \scite{ss.5}) is $<
k^*$ replacing the type by fixed finite sub-types;
\sn
\item "{$\boxdot_2$}"  moreover $C = r({\frak C})$ for some 
non-algebraic $r(x)$;
\sn
\item "{$\boxdot_3$}"  $k^*$ can serve as $k$ in clause (d) of $(St)_2$. 
\ermn
Let $\bar z_n = \bar z^0 \char 94 \ldots \char 94 \bar z^{n-1}$.  
We shall now define by
induction on $n < \omega$ formulas $\varphi_n(x,\bar z_n)$ and
$\vartheta_n(x_1,x_2,\bar z_n)$ also written $\varphi^n_{\bar z_n}(x),
\vartheta^n_{{\bar z}_n}(x_1,x_2)$.
\mn
\ub{Case 1}:  $n=0$.

So ($\bar z_n = <>$, and) $\varphi_n(x) = \varphi(x)$ and
$\vartheta_n(x_1,x_2) = (x_1 = x_2)$.
\mn
\ub{Case 2}:  $n = m+1$.

Let $\varphi^n_{{\bar z}_n}(x) := \varphi^m_{{\bar z}_m}(x) \and
(\exists x')[x'Ez_m \wedge \vartheta^m_{{\bar z}_m}(x',x)]$ and
$\vartheta^n_{\bar z_n}(x_1,x_2) := \varphi^m_{{\bar z}_m}(x_2) \and 
\varphi(x_1) \and (\exists x')[\vartheta(x_1,x',\bar z^m) \wedge
\vartheta^m_{\bar z_m}(x',x_2) \wedge x' Ez_m]$.

We now prove by induction on $n$ that:
\mr
\item "{$(*)_n$}"  if $\bar c = \langle c_\ell:\ell < n\rangle$ and $c_\ell
\in C \backslash ac \ell\{c_k:k < \ell\}$ for $\ell < n$ (so
$\vartheta(x,y,\bar d_{c_\ell})$ is a finite to finite 
function from $\varphi(x)$ into $x E c_\ell$ for $\ell < n$) and $\bar
d = \bar d_{c_0} \char 94 \bar d_{c_1} \char 94 \ldots \char 94 \bar
d_{c_{n-1}}$ then
{\roster
\itemitem{ $(\alpha)$ }  $\varphi^n_{\bar d}({\frak C})$ is an
infinite subset of $\varphi({\frak C})$ 
\sn
\itemitem{ $(\beta)$ } $\vartheta^n_{\bar d}(x_1,x_2)$ is a
finite to finite function from  a co-finite subset 
of $\varphi({\frak C})$ into a subset of $\varphi^n_{\bar d}({\frak C})$
\sn
\itemitem{ $(\gamma)$ }  if $n = m+1$ and 
$e \in \varphi^n_{\bar d}({\frak C})$ then 
$(\exists c' \in ac \ell(\bar c \restriction m \cup \{e\}))[c' E c_m]$
\sn
\itemitem{ $(\delta)$ }  $m<n \Rightarrow \varphi^m_{\bar d
\restriction jm}({\frak C}) \subseteq \varphi^n_{\bar d}({\frak C})$. 
\endroster}
\ermn
This is straight.  Let $I$ be a linear order such that any interval
has $< |T|$ members.

By $\boxdot_2,\boxdot_3$ there are $c_t \in C$ for $t \in I$
pairwise distinct, let $\bar d_t = \bar d_{c_t}$ so 
$\theta(x,y,\bar d_t)$ is a finite to finite
function from $\varphi(x)$ into $xEc_t$ such that
$\langle \bar d_t:t \in I \rangle$ is an indisernible sequence (e.g. use
$\boxdot$ above).

Now for every $<_I$-increasing sequence 
$\bar t = \langle t_n:n < \omega \rangle$ 
we consider $\bar c^n_{\bar t} = \bar c^n_{\bar t \restriction n} =
\bar d_{t_0} \char 94 \ldots \bar d_{t_{n-1}}$ and 
$p_{\bar t} = \{\varphi^n_{{\bar c}^n_{\bar t \restriction n}}(x):
n < \omega\}$.

Now
\mr
\item "{$\circledast_1$}"  for $\bar t$ as above $p_{\bar t}$ is
consistent. 
\ermn
[Why?  By $(*)_n(\alpha)$ there is an element $e \in 
\varphi^n_{\bar c^n_{\bar t}}({\frak C})$, by $(*)_n(\delta)$ the
element $e$ satisfies $\{\varphi_{\bar c^m_{\bar t \restriction m}}(x):
m \le n\}$.  As this holds for every $n$, the set $p_{\bar t} =
\{\varphi^n_{\bar d^n_{\bar t \restriction n}}(x):n < \omega\}$ is
finitely satisfiable as required.]
\mr
\item "{$\circledast_2$}"  if $e$ realizes $p_{\bar t}$ then for every
$n$ there is an element $e'$ algebraic over
$\{e,\bar d_{t_0},\dotsc,\bar d_{t_{n-1}}\}$ such that $e'
Eb^1_{t_n}$. 
\ermn
[Why?  By $(*)_n(\gamma)$.]
\mr
\item "{$\circledast_3$}"  if $e$ realizes $p_{\bar t}$ then for every
$n$ the set $\{s \in I$: there is $e'$ algebraic over
$\{e,\bar d^1_{t_0},\dotsc,\bar d^1_{t_{n-1}}\}$ such that $e' Ec^1_s\}$ has
$\le |T|$ members. 
\ermn
[Why?  There are $\le |T|$ such $e'$ and for each $e'$ by clause 
(d) of (St)$_2$ there are only finitely many such $s \in I$ 
(if we phrase it more carefully
we get that there are $< k_n(< \omega)$ many members).]
\nl
This is more than enough to show $T$ is not strongly$^+$ dependent.
\hfill$\square_{\scite{ss.6}}$
\enddemo
\bn
\ub{Discussion}:   We may phrase \scite{ss.6} for ideals of small formulas.
\bigskip

\proclaim{\stag{ss.8} Claim}   If $T$ is strongly$^1$ dependent and
$\ell = 1,2,3,4$, \ub{then} the statement $\circledast_\ell$ below 
is impossible where:
\mr
\item "{$\circledast_1(a)$}"  $\langle \bar a_\alpha:\alpha < \lambda
\rangle$ is an indiscernible sequence over $A$
\sn
\item "{${{}}(b)$}"  $u_n \subseteq \lambda$ is finite, [non-empty]
with $\langle u_n:n < \omega \rangle$ having pairwise 
disjoint convex hull 
\sn
\item "{${{}}(c)$}"  $\bar b \in {}^{\omega >} {\frak C}$
\sn
\item "{${{}}(d)$}"  for each $n$ for some $\alpha_n,k$ and $t^{\bold
t}_{n(0)} < \ldots < t^{\bold t}_{n(k-1)} \in u_n$ for $\bold t \in
\{$ {\bf false,truth}$\}$ and $\bar c_n \in {}^{\omega >} A$ and
$\varphi$ we have ${\frak C} \models \varphi(\bar c,
a_{t^{\bold t}_{n(0)}},\dotsc,a_{t^{\bold t}_{n(k-1)}},
\bar c_n)^{\bold t}$ for both values of $\bold t$
\sn
\item "{$\circledast_2\,\,\,\,\,\,\,\,$}"  like
$\circledast_1$ but allows $\bar
a_\alpha$ to be infinite 
\sn
\item "{$\circledast_3(a)$}"  $\langle \bar a^n_\alpha:\alpha <
\lambda \rangle$ is an indiscernible sequence over $A \cup
\{\bar a^m_\beta:m < \omega,m \ne n \text{ and } \beta < \lambda\}$
\sn
\item "{$(b)$}"  $\bar a^n_\alpha \ne \bar a^n_{\alpha +1}$
\sn
\item "{$(c)$}"  some $\bar a \in {\frak C}$ satisfies $n < \omega
\Rightarrow acl(A \cup \{a\}) \cap \{\bar a^n_\alpha:\alpha < \lambda\} \ne
\emptyset$
\sn
\item "{$\circledast_4\,\,\,\,\,\,\,\,$}"  like 
$\circledast_3$ but replace clause (c) by
\sn
\item "{$(c)'$}"   for some $\bar a
\in {\frak C}$ for every $n$ the sequence $\langle \bar a^n_\alpha:\alpha <
\lambda \rangle$ is not an indiscernible sequence over $A \cup \bar a$.
\endroster
\endproclaim
\bigskip

\demo{Proof}  Similar to the previous ones.
\enddemo
\bn
\margintag{ss.8.1}\ub{\stag{ss.8.1} Discussion}:  1) We have asked: 
show that the theory of the $p$-adic field is strongly dependent.

Udi Hrushovski has noted that the criterion (St)$_2$ from \scite{ss.6} from
this section apply so $T$ is not strongly$^2$ dependent.   Namely take the
following equivalence relation on $\Bbb Z_p$: val$(x-y) \ge 
\text{ val}(c)$, where $c$ is some fixed element with infinite valuation.
Given $x$, the map $y \mapsto (x + cy)$ is a bijection between $\Bbb
Z_p$ and the class.
\nl
2) By \cite{Sh:863} this theory is strongly$^1$ dependent.
\nl
3) Onshuus shows that also the theory of the field of the reals is not
   strongly$^2$ dependent (e.g. though Claim \scite{ss.6} does not
   apply but its proof works (using pairwise not too near $\bar b$'s,
   in general just an uncountable set of $\bar b$'s). 
\nl
4) In \cite{Sh:863} we prove reasonable existence of
indiscernibles for strongly dependent $T$ (and in \scite{ss.1} we can
use the case $\ell g(\bar x)=1$).
\bigskip

\proclaim{\stag{ss.8.2} Claim}  
Assume $x=1$ or $x=2$ or $x=1^*$ or $x=2^*$.
\nl
1) If $M \prec {\frak C},A \subseteq {\frak C}$ 
\ub{then} (the complete first order) $Th({\frak B}_{M,M,A}$ 
from \scite{t.3}(2) is strongly$^x$ dependent iff $T$ is strongly$^x$ 
dependent.  If $T$ is dependent this theory is equal to $T^*_{M,A}$,
see \scite{t.3}(4), \scite{t.4A}(4).
\nl
2) $\kappa_{\text{ict}}(T) = \kappa_{\text{ict}}(\text{\rm Th}
({\frak B}_{M,MA}))$ if $M \prec {\frak C},A \subseteq {\frak C}$
\nl
3) If $T_1 \subseteq T_2$ are complete first
order theories (so $\tau(T_1) \subseteq \tau(T_2))$ then
\mr
\item "{$(a)$}"  if $T_2$ is strongly$^x$ dependent then so is
$T_1$
\sn
\item "{$(b)$}"  $\kappa_{\text{ict}}(T_1) \le
\kappa_{\text{ict}}(T_2)$.
\ermn
4) If $T_1 \subseteq T_2$ are complete first order and $\tau(T_2)
\backslash \tau(T_1)$ consist of individual constants only \ub{then}
\mr
\item "{$(\alpha)$}"   $T_2$ is strongly$^x$ dependent iff $T_1$ is
strongly$^x$ dependent
\sn
\item "{$(\beta)$}"  $\kappa_{\text{ict}}(T_1) =
\kappa_{\text{ict}}(T_2)$.
\ermn
5) For $T$ is strongly$^x$ dependent iff $T^{\text{eq}}$
is strongly$^x$ dependent; similarly for strongly$^{\ell,*}$.
\nl
6) $\kappa_{\text{ict}}(T) = \kappa_{\text{ict}}(T^{\text{eq}})$.
\endproclaim
\bigskip

\demo{Proof}  Easy.
\enddemo
\newpage

\head {\S4 Definable groups} \endhead  \resetall \sectno=4
 \spuriousreset
\bigskip

\demo{\stag{ag.0} Context}
\mr
\item "{$(a)$}"  $T$ is a first order complete theory
\sn
\item "{$(b)$}"  ${\frak C}$ is a monster model of $T$.
\endroster
\enddemo
\bigskip

We try here to generalize the theorem on the existence of commutative
infinite subgroups for stable $T$ to dependent $T$.  Theorems
on definable groups in a monster ${\frak C}$, Th$({\frak C})$ stable,
are well known.
\definition{\stag{ag.1} Definition}  1) We say $G$ is a type-definable group (in
${\frak C}$) \ub{if} $G = (p,*,\text{inv}) =
(p^G,*^G,\text{inv}^G)$ where
\mr
\item "{$(a)$}"  $p = p(x)$ is a type
\sn
\item "{$(b)$}"  $*$ is a two-place function on ${\frak C}$, possibly
partial, definable (in ${\frak C}$),
we normally write $ab$ instead of $a*b$ or $*(a,b)$
\sn
\item "{$(c)$}"  $(p({\frak C}),*)$ is a group, we write $x \in G$ for
$x \in p({\frak C})$;
\sn
\item "{$(d)$}"  inv$^G$ is a (partial) unary function, definable (in
${\frak C}$), which on $p({\frak C})$ is 
the inverse, so if no confusion arises we shall write $(x)^{-1}$ for inv$(x)$.
\ermn
1A) We let $B^G_2$ be the set of parameters appearing in $p^G$; let
$B^G$ be the set of parameters appearing in $p^G$ or in the
definition of $*$ or of inv$^G$.
\nl
2) We say $G$ is a definable group \ub{if} $p(x)$ is a formula, i.e.,
a singleton. \nl
3) We say $G$ is an almost type definable group if $p(x)$ is replaced
by $\bar p = \langle p_i(x):i < \delta \rangle,p_i({\frak C})$
increasing with $i$ and $\bar p({\frak C})$ is defined as
$\cup\{p_i({\frak C}):i < \delta\}$.
\enddefinition
\bigskip

\remark{Remark}  Of course, we can use $p(\bar x)$ and/or work in
${\frak C}^{\text{eq}}$.
\endremark
\bigskip

\proclaim{\stag{sg.1} Claim}  Assume
\mr
\item "{$(a)$}"  $T$ is dependent
\sn
\item "{$(b)$}"  $G$ is a definable group in ${\frak C}$ or just
type-definable 
\sn
\item "{$(c)$}"  $A \subseteq G$ is a set of pairwise commuting
elements, $D$ a non-principal ultrafilter on $A$ or just
\sn
\item "{$(c)^-$}"   $A \subseteq G,D$ a non-principal 
ultrafilter on $A$ such that $(\forall^D a_1)(\forall^D
a_2)(a_1a_2 = a_2a_1)$, where $\forall^D x \varphi(x,\bar a)$ means
$\{b \in \text{\rm Dom}(D):{\frak C} \models \varphi[b,\bar a]\} \in D$.
\ermn
\ub{Then} there is a formula $\varphi(x,\bar a)$ such that:
\mr
\item "{$(\alpha)$}"  $\varphi(x,\bar a) \in {\text{\rm Av\/}}(\bar a,D)$
\sn
\item "{$(\beta)$}"  $G \cap \varphi({\frak C},\bar a)$ is an abelian
subgroup of $G$
\sn
\item "{$(\gamma)$}"  $\bar a \subseteq A \cup B^G \cup 
\{c:c \text{ realizes } {\text{\rm Av\/}}(A \cup B^G_2,D)\}$.
\endroster
\endproclaim
\bigskip

\remark{\stag{sg.1.3} Remark}  1) If $D$ is a 
principal ultrafilter say $\{a^*\} \in D$
then $\varphi(x,\bar a)$ is essentially Cm$_G(\text{Cm}_G(a^*))$ so
no new point, recalling 
(Cm$_G(A) = \{x \in A:x$ commutes with every $a \in A\}$. 
\nl
2) If $D$ is a non-principal ultrafilter, \ub{then} necessarily
$\varphi(x,\bar a)$ is not algebraic as it belongs to {\rm Av}$(\bar a,D)$.
\endremark
\bigskip

\demo{Proof}  We try to choose $a_n,b_n$ by induction on $n < \omega$
such that:
\mr
\item "{$(i)$}"  $a_n,b_n$ realizes $p_n(x) := 
{\text{\rm Av\/}}(A_n,D)$ where $A_n = A \cup B^G \cup\{a_k,b_k:
k < n\}$ so as $A \in D,A \subseteq G$ necessarily $p^G(x) \subseteq p_n(x)$
\sn
\item "{$(ii)$}"  $a_n,b_n$ does not commute (in $G$, they are in $G$ because
$p^G(x) \subseteq p_n$).
\endroster
\enddemo
\bn
\ub{Case 1}:  We succeed.

Assume $n < m < \omega,c' \in \{a_n,b_n\}$ and $c'' \in \{a_m,b_m\}$ clearly
$c',c''$ are in $G$.  Now we shall show that they commute because
$c''$ realizes Av$(A \cup B^G \cup\{c'\},D)$ and
$c'$ realizes Av$(A \cup B^G_2,D)$ recalling either assumption (c) 
about commuting in $A$ or assumption $(c)^-$.  
Hence if $k < \omega,n_0 < \ldots < n_{k-1} < \omega$ and $n < \omega$
then $c :=
b_{n_0}b_{n_1} \ldots b_{n_{k-1}}$ satisfies: $c,a_n$ commute iff $n
\notin \{n_0,\dotsc,n_{k-1}\}$, so $\varphi(x,y) = [xy = yx]$ has the
independence property contradicting assumption (a).  
\bn
\ub{Case 2}: We are stuck at $n < \omega$.

So $p_n(x) \cup p_n(y) \vdash (xy = yx)$, \ub{hence} there is a formula
$\psi(x,\bar a^*) \in {\text{\rm Av\/}}(A_n,D)$ such that
\mr
\item "{$(*)_1$}"  $\psi(x,\bar a^*) \wedge
\psi(y,\bar a^*) \vdash xy = yx$ (so both products are well defined).
\ermn
Let $p^G(x) = \{\vartheta(x,\bar a)\}$ or just $p(x) \vdash
\vartheta(x,\bar a)$ and $\theta(x,\bar a) \wedge \theta(y,\bar a)
\vdash$ ($xy$ well defined); \wilog \,
$\bar a \trianglelefteq \bar a^*$ and $\psi(x,\bar a^*) \vdash
\vartheta(x,\bar a)$ and
let $\vartheta^*(x,\bar a^*) = \vartheta(x,\bar a) \wedge 
(\forall y)[\psi(y,\bar a^*) \rightarrow yx = xy$ (so both well
defined)].  So $\psi(x,\bar a^*) \vdash \vartheta^*(x,\bar a^*)$. \nl
Let

$$
\varphi(x) = \varphi(x,\bar a^*) = \vartheta(x,\bar a) \wedge (\forall
y)[\vartheta^*(y,\bar a^*) \rightarrow xy=yx \text{ (both well defined)}].
$$
\mn
So $\psi({\frak C},\bar a^*) \subseteq \varphi({\frak C},\bar a^*)
\subseteq \vartheta^*({\frak C},\bar a^*) \subseteq \vartheta({\frak
C},\bar a)$.
So $\psi(x,\bar a^*) \vdash \varphi(x,\bar a)$ and $\varphi(x,\bar a) 
\vdash \vartheta^*(x,\bar a^*)$ hence the formula 
$\varphi(x,\bar a^*)$ belongs to the type 
$p_n(x)$ which is equal to $\text{Av}(A_n,D)$ hence
$\varphi(x,\bar a^*) \in \text{ Av}(\bar a^*,D)$ and $\bar a^* 
\subseteq A_n \subseteq A \cup B^G \cup \cup \{c:c$ 
realizes Av$(A \cup B^G,D)\}$.

We are done as $\varphi({\frak C},\bar a^*) \cap G$ is a 
subgroup and is abelian by the
definition of $\varphi(x)$.  \hfill$\square_{\scite{sg.1}}$
\bigskip

\proclaim{\stag{sg.2} Claim}  Assume
\mr
\item "{$(a)$}"  $G$ is a definable (infinite) group, (or just
type-definable) 
\sn
\item "{$(b)$}"  every element of $G \backslash \{e_G\}$ 
commutes with only finitely many
others
\sn
\item "{$(c)$}"  $G$ has infinitely many pairwise non-conjugate members.
\ermn
\ub{Then} $T$ is not strongly$^+$ dependent.
\endproclaim
\bigskip

\demo{Proof}  Assume first $p^G = \{\varphi(x)\}$.

Let $xEy := [x,y$ are conjugates], clearly it is an 
equivalence relation, and let

$$
\vartheta(x_1,x_2,y) := (x_1 = x_2 y x^{-1}_2).
$$
\mn
Note that: if $M \models 
\vartheta(x_1,z_1,y) \wedge \vartheta (x_1,z_2,y)$ then
$M \models z_1 y z^{-1}_1 = z_2 y z^{-1}_2$ hence $M \models (z^{-1}_2
z_1) y = y(z^{-1}_2 z)$ so $z^{-1}_2 z_1 \in \text{\rm Cm}_G(y)$ so
$\{z:\vartheta(x_1,z,y)\}$ is finite.  Trivially $\{x_1:
\vartheta(x_1,x_2,y)\}$ is finite.

We now get a contradiction by \scite{ss.6}: 
$\varphi(x),\vartheta(x_1,x_2,y)$ satisfies the demands in $(St)_1$
there,  which is impossible if $T$ is strongly$^+$ dependent; so we are done.

If $p^G$ is a 
type use $(St)_3$ of \scite{ss.6}.  \hfill$\square_{\scite{sg.2}}$
\enddemo
\bigskip

\definition{\stag{ag.3} Definition}  1) A place $\bold p$ is a tuple
$(p,B,D,*,\text{inv}) = (p^{\bold p},B^{\bold p},D^{\bold p},
*_{\bold p}$,inv$_{\bold p}) = (p[\bold p],B[\bold p],D[\bold
p],*[\bold p]$,inv$[\bold p])$ such that:
\mr
\item "{$(a)$}"  $B$ is a set $\subseteq {\frak C}$, $D$ is an 
ultrafilter on $B,p \subseteq \text{ Av}(B,D)$
\sn
\item "{$(b)$}"  $*$ is a partial two-place function defined with
parameters from $B$; we shall write $a *_{\bold p} b$ or, when clear from
the context, $a*b$ or $ab$
\sn
\item "{$(c)$}"  inv is a partial unary function definable from
parameters in $B$.
\ermn
1A) $\bold p$ is non-trival if for every $A$ the type   
Av$(A,D)$ is not algebraic. 
\nl
2) We say $\bold p$ is weakly a place in a definable group $G$ or type
definable group $G$ \ub{if} $\bold p$ is a place, $p^{\bold p}
\vdash p^G$, the set $B^{\bold p}$ includes Dom$(p^G)$
and the operations agree on $p_{\bold p}[{\frak C}]$ 
when the place operations are defined. \nl
2A) If those operations are the same, we say 
that $\bold p$ is strongly a place in $G$. \nl
3) We say $\bold p_1 \le \bold p_2$ if both are places, $B^{{\bold
p}_1} \subseteq B^{{\bold p}_2}$ and $p^{\bold p_2} \vdash 
p^{\bold p_1}$ and the operations are same. \nl
4) $\bold p \le_{\text{dir}} \bold q$ if $\bold p \le \bold q$ and
$B^{\bold q} \subseteq A \Rightarrow \text{ Av}(A,D^{\bold p}) =
\text{ Av}(A,D^{\bold q})$.
\enddefinition
\bigskip

\definition{\stag{ag.4} Definition}  1) A place $\bold p$ is
$\sigma$-closed \ub{when}:
\mr
\item "{$(a)$}"  $\sigma$ has the form $\sigma(\bar x_1;\ldots,
;\bar x_{n(*)})$, a term in the vocabulary of groups
\sn
\item "{$(b)$}"  if $\bar a_\ell \in {}^{(\ell g(\bar x_\ell)} {\frak C}$,
for $\ell = 1,\dotsc,n(*)$ and $B \subseteq A$, \ub{then} 
$\sigma(\bar a_1,\dotsc,\bar a_{n(*)})$ is well defined \footnote{so all the
stages in the computation of $\sigma(\bar a_0;\dotsc;\bar a_{n(*)})$
should be well defined} and realizes Av$(A,D)$ provided that
{\roster
\itemitem{ $(*)$ }  $n \le n(*) \and \ell < \ell g(\bar a_n)
\Rightarrow a_{n,\ell}$ realizes Av$(A \cup \bar a_1 \char 94 
\ldots \char 94 \bar a_{n-1},D)$.
\endroster}
\ermn
2) A place $\bold p$ is $(\sigma_1 = \sigma_2)$-good or satisfies
$(\sigma_1 = \sigma_2)$ when
\mr
\item "{$(a)$}"  $\sigma_\ell = \sigma_\ell(\bar x_1,\dotsc,\bar
x_{n(*)})$ a term in the vocabulary of groups for $\ell =1,2$ (so
e.g. $(x_1 x_2)x_3,x_1(x_2 x_3)$ are considered as different terms)
\sn
\item "{$(b)$}"  if $\bar a_\ell \in {}^{(\ell g({\bar x}_n)} {\frak C}$ for
$\ell \le n$ then $\sigma_1(\bar a_1;\dotsc;\bar a_{n(*)}) =
\sigma_2(\bar a_1;\ldots;\bar a_{n(*)})$ whenever $(*)$ of part (1)
holds for $A=B$; so both are well defined.
\ermn
3) We can replace $\sigma$ in part (1) by a set of terms.  Similarly
in part (2) for a set of pairs. \nl
4) We may write $x_\ell$ instead of $\langle x_\ell \rangle$.  So if we
write $\sigma(\bar x_1;\bar x_2) = \sigma(x_1;x_2) = x_1 x_2$ or
$\sigma = x_1 x_2$ we mean $x_1 = x_{1,0},x_2 = x_{2,0},\bar x_1 =
\langle x_{1,0} \rangle,\bar x_2 = \langle x_{2,0} \rangle$.  We may
use also $\sigma(\bar x;\bar y)$ instead of $\sigma(\bar x_1;\bar
x_2)$ and $\sigma(\bar x;\bar y;\bar z)$ similarly. 
\enddefinition
\bigskip

\definition{\stag{ag.5} Definition}  1) We say a place $\bold p$ is a poor
semi-group \ub{if} it is $\sigma$-closed for $\sigma = xy$ and
satisfies $(x_1x_2)x_3 = x_1(x_2x_3)$. 
\nl
2) We say a place $\bold p$ is a poor group \ub{if} it is a poor 
semi-group and is $\sigma$-closed for $\sigma = (x_1)^{-1}x_2$.
\nl
3) We say a place $\bold p$ is a quasi semi-group \ub{if} 
for any semi group term  $\sigma_*(\bar x),\bold p$ is $\sigma$-closed
for $\sigma(\bar x;y) = \sigma_*(\bar x)y$.
\nl
4) We say a place $\bold p$ is a quasi group \ub{if} for any semi-group terms
$\sigma_1(\bar x),\sigma_2(\bar x)$ we place $\bold p$ is $\sigma$-closed for 
$\sigma(\bar x;y) = \sigma_1(\bar x) y \sigma_2(\bar x)$.
\nl
5) We say $\bold p$ is abelian (or is commutative) 
if it is $(xy)$-closed and satisfies $xy = yx$. \nl
6) We say $\bold p$ is affine if $\bold p$ is $(xy^{-1}z)$-closed.
\nl
7) We say that a place $\bold p$ is a pseudo semi-group \ub{when}: if
the terms
$\sigma_1(x_1,\dotsc,x_n),\sigma_2(x_1,\dotsc,x_n)$ are equal in semi-groups
then $\bold p$ satisfies
$\sigma_1(x_1,\dotsc,x_n) = \sigma_2(x_1,\dotsc,x_n)$.
\nl
8) We say that a place $\bold p$ is a pseudo group if any term
$\sigma_1(x_1,\dotsc,x_n),\sigma_2(x_1,\dotsc,x_n)$ which are equal in groups,
$\bold p$ satisfies $\sigma_1(x_1,\dotsc,x_n) = \sigma_2(x_1,\dotsc,x_n)$.
\enddefinition
\bigskip

\definition{\stag{ag.6} Definition}  We say a place $\bold p$ is a
group if $G = G^{\bold p} = (\text{Av}(B^{\bold p},D),{}^* \bold
p$,inv$_{\bold p})$ is a group. 
\enddefinition
\bigskip

\proclaim{\stag{ag.7} Claim}  1) The obvious implications hold. \nl
2) If we use $\bar{\bold b}$ every $\bar{\bold b}'$ realizing the same
type has the same properties. 
\nl
3) For a place $\bold p$ the assertion ``$\bold p$ satisfies
$\sigma(\bar x_1,\dotsc,\bar x_{n(*)}) = \sigma(\bar
x_1,\dotsc,\bar x_{n(*)})$" means just that in Definition
\scite{ag.4} the term $\sigma(a_1,\dotsc,a_n)$ is well defined.
\endproclaim
\bn
We now note that there are places
\proclaim{\stag{ag.8} Claim}  1) Assume that $G$ is a definable group and
$a_n \in p^G[{\frak C}]$ for $n < \omega$. We define $a_{[u]} \in
p^G[{\frak C}]$ for any finite non-empty $u \subseteq \omega$ by
induction on $|u|$, if $u = \{n\}$ then $a_{[u]} = a_n$, if $|u| > 1$,
{\rm max}$(u) = n$ then $a_{[u]} = a_{[u \backslash \{n\}]} *^G a_n$ and we
are assuming they are all well defined and $a_{[u_1]} \ne a_{[u_2]}$
when $u_1 \triangleleft u_2$.  \ub{Then} we can find $D^*,\bold q$ such
that:
\mr
\item "{$(a)$}"  $\bold q$ is a place inside $G$
\sn
\item "{$(b)$}"  $\bold q$ is a poor semi-group and non-trivial
\sn
\item "{$(c)$}"  $B^{\bold q} = B^G \cup \bigcup\{a_{[u]}:u \subseteq
\omega$ is finite$\}$
\sn
\item "{$(d)$}"  $D^*$ is an ultrafilter on $[\omega]^{< \aleph_0}$ such
that $(\forall n)([\omega \backslash n]^{<\aleph_0} \in D^*)$ and for every $Y
\in D^*$ we can find $Y' \subseteq Y$ from $D^*$ closed under convex
union, i.e., if $u,v \in Y'$ and $\text{max}(u) < \text{\rm min}(v)$ then $u
\cup v \in Y'$
\sn
\item "{$(e)$}"  $D^{\bold q} = \{\{a_{[u]}:u \in Y\}:Y \in D^*\}$
\sn
\item "{$(f)$}"  if the $a_n$'s commute (i.e. $a_na_m = a_ma_n$ for $n
\ne m$) \ub{then} $\bold q$ is abelian.
\endroster
\endproclaim
\bigskip

\demo{Proof}  By a well known theorem of Glazer \footnote{his proof
uses the operations from clause (d) of \scite{ng.1} and \scite{ng.3}
below}, relative of
Hindman theorem saying $D^*$ as in clause (d) exists, see Comfort
\cite{Cmf77}.  \hfill$\square_{\scite{ag.8}}$
\enddemo
\bigskip

\remark{\stag{ag.8.2} Remark}  1) This can be combined naturally with \S1.
\endremark
\bigskip

\proclaim{\stag{ag.9} Claim}  1) Assume
\mr
\item "{$(a)$}"  $\bold p$ is a place in a type-definable group
(or much less)
\sn
\item "{$(b)$}"  the place $\bold p$ is a semi-group
\sn
\item "{$(c)$}"  $\bold p$ is commutative (in the sense of Definition
\scite{ag.4} + \scite{ag.5}, so satisfies $\sigma_1(x;y) = [x*y =y*x]$

\sn
\item "{$(d)$}"  if $A \supseteq B^{\bold p}$ then for some $b,c$
realizing {\rm Av}$(D^{\bold p},A),c *_G b,b *_G c$ are (necessarily well
defined, and) distinct.
\ermn
\ub{Then} $T$ has the independence property.
\nl
2) We can weaken clause (a) to
\mr
\item "{$(a)'$}"  $\bold p$ is a place such that for $n < \omega$
and $\langle a_1a'_1\rangle, \dotsc,\langle a_na'_n \rangle$ are as
in Definition \scite{ag.4} and $a_\ell \ne a'_\ell \Leftrightarrow \ell =
m$ then $a_1,\dotsc,a_{m-1},a_m a_{m+1} \ldots a_n
\ne a_1,\dotsc,a_{m-1} a'_m a_{m+1} \ldots a_m$.
\endroster
\endproclaim 
\bigskip

\remark{Remark}  This is related to the well known theorems on stable
theories (see Zilber and Hrushovski's works).
\endremark
\bigskip

\demo{Proof}  1) We choose $A_i,b_i,c_i$ by induction on $i < \omega$.

In stage $i$ first let $A_i = B^{\bold p} \cup \{b_j,c_j:j < i\}$ and
add $B^G$ if $B^G \nsubseteq B^{\bold p}$. \nl
Second, choose $b_i,c_i$ realizing Av$(A_i,D^{\bold p})$ such that
$b_i * c_i \ne c_i * b_i$.

Now if $i < j < \omega$ any $a' \in \{b_i,c_i\},a'' \in \{b_j,c_j\}$
then $a'$ realizes Av$(A_i,D^{\bold p})$ and $a''$ realizes
Av$(A_j,D^{\bold p})$ which include Av$(A_i \cup \{a'\},D^{\bold p})$.
So by assumption (c), the elements $a',a''$ commute in $G$.

So as is well known, for $n < \omega,i_0 < i_1 < \ldots < i_n$ the
element $b_{i_0} * b_{i_1} * \ldots * b_{i_{n-1}}$ commute in $G$ with
$a_j$ iff $j \notin \{i_0,\dotsc,i_{n-1}\}$ hence $T$ has the
independence property.  \nl
2) Similarly.    \hfill$\square_{\scite{ag.9}}$ 
\enddemo
\bn
Note that \scite{ag.10} is interesting for $G$ with a finite bound on
the order of elements as if $a \in G$ has infinite order then
Cm$_G(\text{\rm Gm}_G(a))$ is as desired.
\demo{\stag{ag.10} Conclusion}  [$T$ is dependent].

Assume $G$ is a definable group. \nl
1) If $\bold p$ is a commutative semi-group in $G$, non-trivial,
\ub{then} for some formula $\varphi(x,\bar a)$ such that $\varphi(\bar
x) \vdash ``x \in G"$ and $\varphi(x,\bar a) \in \text{ Av}(\bar a,D^{\bold
p})$ and $G \restriction \varphi({\frak C})$ is a commutative place. \nl
2) If $G$ has an infinite abelian subgroup, \ub{then} it has an
infinite definable commutative subgroup.
\enddemo
\bigskip

\demo{Proof}  1) By \scite{ag.9} for some $A \supseteq B^{\bold p}$
for every $b,c$ realizing $q := \text{ Av}(A,D^{\bold p})$ we have: 
the elements of $q({\frak C})$, which are all
in $G$, pairwise commute.  By compactness there is a formula
$\varphi_1(x) \in p[\bold p]$ such that the elements of
$\varphi_1({\frak C}) \cap G$ pairwise commute and \wilog \,
$\varphi_1(x) \vdash [x \in G]$; note however that 
this set is not necessarily a
subgroup.  Let $\varphi_2(x) := [x \in G] \wedge (\forall y)
(\varphi_1(y) \rightarrow x*y = y*x]$.  Clearly $\varphi_1({\frak
C}) \subseteq \varphi_2({\frak C}) \subseteq G$ and every member of
$\varphi_2({\frak C})$ commutes with every member of $\varphi_1({\frak
C})$.  So $\varphi(z) := [z \in G] \wedge (\forall y)[\varphi_2(y)
\rightarrow yz = zy]$ is first order and defines the center of $G
\restriction \varphi_2[{\frak C}]$ which includes $\varphi_1({\frak
C})$, so we are done. \nl
2) Let $G' \subseteq G$ be infinite abelian.  Choose by induction on
$n < \omega,a_n \in G'$ as required in \scite{ag.9} and then apply it.
\hfill$\square_{\scite{ag.10}}$
\enddemo
\bigskip

\remark{\stag{ag.11} Remark}  So \scite{ag.10} tells us that having
some commutativity implies having alot.  If in \scite{ag.9} every
$a_{[u]}$ is not in any ``small" definable set defined with parameters in
$B^{\bold p} \cup \{a_n:n < \max(u)\}$, \ub{then} also $\varphi(x,\bar
a)$ is not small where small means some reasonably definable ideal. 
\endremark
\bn
\centerline {$* \qquad * \qquad *$}
\bn
\definition{\stag{ng.1} Definition}  Assume
\mr
\item "{$(a)$}"  $G$ is a type definable semi-group
\sn
\item "{$(b)$}"  $M \supseteq B^G$ is $(|T| + |B^G|)^+$-saturated
\sn
\item "{$(c)$}"  ${\frak D}$ is the set of ultrafilters $D$ on $M$
such that $p^G \subseteq \text{ Av}(M,D)$
\sn
\item "{$(d)$}"  on ${\frak D} = {\frak D}_{G,M}$ we define an
operation
$$
\align
D_1 * D_2 = D_3 \text{ \ub{iff}} &\text{ for any } A 
\supseteq M \text{ and } a
\text{ realizing} \\
  &\text{ Av}(A,D_1) \text{ and } b \text{ realizing Av}(A+a,D_2)
\text{ the element} \\
  &a*b \text{ realizes Av}(A,D_3).
\endalign
$$
\sn
\item "{$(e)$}"  $ID_{G,M} = \{D \in {\frak D}_{G,M}:D*D = D\}$
\sn
\item "{$(f)$}"  $H^{\text{left}}_{G,M} = \{a \in G$: for every $D \in
{\frak D}$, and $A \supseteq M$ if $b$ realizes 
Av$(A+a,D)$ then $a*b$ realizes Av$(A,D)\}$
\sn
\item "{$(g)$}"  $H^{\text{right}}_{G,M}$ similarly using $b*a$
\sn
\item "{$(h)$}"   $H_{G,M} = H^{\text{left}}_{G,M} \cap
H^{\text{right}}_{G,M}$.
\ermn
The following as in \scite{ag.8}.  
\enddefinition
\bn
\margintag{ng.3}\ub{\stag{ng.3} Fact}:  ${\frak D}$ is a semi-group, i.e., associativity
holds and the operation is continuous in the second variable hence there is
an idempotent (even every non-empty 
subset closed under $*$ and topologically closed has an idempotent). 
\bn
Note \nl
\margintag{ng.4}\ub{\stag{ng.4} Fact}:  1) If $G$ is a group, \ub{then}
\mr
\item "{$(a)$}"  $H^{\text{left}}_{G,M}$ is a subgroup of $G$, with
bounded index, and is of the form $\cup\{q({\frak C}):q \in 
\bold S^{\text{left}}_{G,M}\}$ for some $\bold S^{\text{left}}_{G,M} 
\subseteq \bold S(M)$
\sn
\item "{$(b)$}"  Similarly
$H^{\text{right}}_{G,M},H_{G,M} = H^{\text{right}}_{G,M} \cap
H^{\text{left}}_{G,M}$ with $\bold S^{\text{right}}_{G,M},
\bold S_{G,M}$.
\ermn
2) If 
$D \in {\frak D}$ is non principal and Av$(M,D) \in {\Cal
S}^{\text{right}}_{G,M}$, \ub{then} for any $A \supseteq M$ and
element $a$ realizing Av$(A,D)$ and $b$ realizing Av$(A+a,D)$ we have
\mr
\item "{$(\alpha)$}"  $a *_G b$ realizes Av$(A,D)$
\sn
\item "{$(\beta)$}"  also $a^{-1} *b \in D$.
\ermn
3) $\bold S^{\text{left}}_{G,M} \subseteq ID_{G,M}$. \nl
4) Similarly for $\bold S^{\text{left}}_{G,M},b *_G a$. \nl
5) If $D \in {\frak D},p = \text{ Av}(M,D) \in \bold S_{G,M}$
\ub{then}
\mr
\item "{$(a)$}"  $\bold p = (M,D,*,\text{inv})$ is a quasi group
\sn
\item "{$(b)$}"  $\{a^{-1} b:a,b \in p(M)\}$ is a subgroup of $G$ with
bound index, in fact is $\{a \in {\frak C}:\text{tp}(a,M) \in {\Cal
S}_{G,M}\}$.
\endroster
\newpage

\head {\S5 Non-forking} \endhead  \resetall \sectno=5
 \spuriousreset
\bigskip

\demo{\stag{10.0} Hypothesis}  $T$ is dependent.
\enddemo
\bigskip

\definition{\stag{10.1} Definition}  \cite{Sh:93}  1) An $\alpha$-type
$p = p(\bar x)$ divides over $B$ \ub{if} some sequence $\bar{\bold b}$
and formula $\varphi(\bar x,\bar y)$ witness it which means
\mr
\item "{$(a)$}"  $\bar{\bold b} = \langle \bar b_n:n < \omega \rangle$ 
is an indiscernible sequence over $B$
\sn
\item "{$(b)$}"   $\varphi(\bar x,\bar y)$ is a formula with $\ell
g(\bar y) = \ell g(\bar b_n)$
\sn
\item "{$(c)$}"  $p \vdash \varphi(\bar x,\bar b_0)$
\sn
\item "{$(d)$}"  $\{\varphi(\bar x,\bar b_n):n < \omega\}$ is
contradictory.
\ermn
1A) Above we say $\varphi(\bar x,\bar b_0)$ explicitly divide over $B$.
\nl
1B) An $\alpha$-type $p=p(\bar x)$ splits strongly over $B$ \ub{when} for
some sequence $\bold{\bar b}$ and formula 
$\varphi(\bar x,\bar y)$ witness it which means:
\mr
\item "{$(a),(b)$}"  as above
\sn
\item "{$(c)$}"  $\varphi(\bar x,\bar b_0),\neg \varphi(\bar x,\bar b_1) 
\in p$.
\ermn
2) An $\alpha$-type $p$ forks over $B$ \ub{if} for some $\langle
\varphi_\ell(\bar x,\bar a_\ell):\ell < k \rangle$ we have $p \vdash
\dsize \bigvee_{\ell < k} \varphi_\ell(\bar x,\bar a_\ell)$ and
$\{\varphi_\ell(\bar x,\bar a_\ell)\}$ divides over $B$ for each $\ell < k$
(note: though $\bar x$ may be infinite, the formulas are finitary). 

We say $p(\bar x)$ exactly forks (or ex-forks) over $B$ \ub{when} some
$\varphi(\bar x,\bar b) \in p$ does exactly fork over $B$, \ub{which}
means that for some $\langle \varphi_\ell(\bar x,\bar b):\ell <
k\rangle$ we have: $\varphi(\bar x,\bar b) \vdash \dsize \bigvee_{\ell
< k} \varphi_\ell(\bar x,\bar b)$ and each 
$\varphi_\ell(\bar x,\bar b)$ explicitly divides over $B$; so we are
\underbar{not} allowed to add dummy parameters.
\nl
3) We say $C/A$ does not fork over $B$ \ub{if} letting $\bold{\bar c}$ list
$C$, tp$(\bar{\bold c},A)$ does not fork over $B$, or what is
equivalent $\bar c \in {}^{\omega >} C \Rightarrow \text{ tp}(\bar
c,A)$ does not fork over $B$ (so below we may write claims
for $\bar c$ and use them for $C$). 
\nl
4) The $m$-type $p$ is f.s. (finitely satisfiable) in $A$ \ub{if} every
finite $q \subseteq p$ is realized by some $\bar b \subseteq A$. \nl
5) The $\Delta$-multiplicity of $p$ over $B$ is Mult$_\Delta(p,B) = 
\text{ sup}\{|\{q \restriction \Delta:p \subseteq 
q \in \bold S^m(M),q$ does not fork over $B\}|:
M \supseteq B \cup \text{ Dom}(p)\}$.
\nl
Omitting $\Delta$ means $\Bbb L(\tau_T)$, omitting $B$ we mean Dom$(p)$.
\enddefinition
\bigskip

\definition{\stag{10.1B} Definition}  1) Let $p = p(\bar x)$ be an
$\alpha$-type and $\Delta$ be a set of $\Bbb L(\tau_T)$-formulas of the
form $\varphi(\bar x,\bar y)$ and $k \le \omega$.  For a type
$p(\bar x)$ we say that it $(\Delta,k)$-divides over $A$ \ub{when} 
some $\bar{\bold b},\varphi(\bar x,\bar y)$ witness it which means
\mr
\item "{$(a)$}"  $\bar{\bold b} = \langle \bar b_n:n < 2k+1 \rangle$
is $\Delta$-indiscernible
\sn
\item "{$(b)$}"  $\varphi(\bar x,\bar y) \in \Bbb L(\tau_T)$
\sn
\item "{$(c)$}"  $p \vdash \varphi(\bar x,\bar b_0)$
\sn
\item "{$(d)$}"  $\{\varphi(\bar x,\bar b_n):n < 2k+1\}$ is
$k$-contradictory. 
\ermn
2) For a type $p(\bar x)$ we say that it 
$(\Delta,k)$-forks over $B$ \ub{when} $p \vdash \dsize \bigvee_{\ell < n}
\varphi_\ell(x,\bar a_\ell)$ for some $n,\varphi_\ell(\bar x,\bar y)$
and $\bar a_\ell$, where each $\varphi_\ell(\bar x,\bar a_\ell)$ does 
$(\Delta,k)$-divides over $B$.
\enddefinition
\bn
\margintag{10.3}\ub{\stag{10.3} Observation}:  0) In Definition \scite{10.1}(1), if $p
= \{\varphi(x,\bar b)\}$ \ub{then} without loss of generality
$\bar b_0 = \bar b$.  If $p$ divides over $B$ \ub{then} $p$ forks over
$B$. \nl
0A)  Forking is preserved by permuting and repeating the 
variables.  If tp$(\bar b \char 94 \bar c,A)$ does not fork 
over $B$ then so does tp$(\bar b,A)$ and both do not divide over $B$.
Similarly for dividing and for exact forking (and later versions).
\nl
1) If $p \in \bold S^m(A)$ is finitely satisfiable
in $B$, \ub{then} $p$ does not fork over $B$; hence every type over
$M$ does not fork over $M$. \nl
2) If $p \in \bold S^m(A)$ does not fork or just does not divide 
over $B \subseteq A$,
\ub{then} $p$ does not split strongly over $B$. (Of course, if $p$
divides over $A$ then $p$ forks over $A$).

The type $p(\bar x)$ divides over $B$ \ub{iff} for some $k < \omega$
and $\varphi_\ell(\bar x,\bar c_\ell) \in p(\bar x)$ for $\ell < k$,
letting $\bar c = \bar c_0 \char 94 \ldots \char 94 \bar c_{k-1}$
the formula $\varphi(\bar x,\bar c) = \dsize \bigwedge_{\ell < k}
\varphi_\ell(\bar x,\bar c_\ell)$ explicitly divides over $B$.  Assume
the type $\{\varphi(\bar x,\bar b)\}$ or $p(\bar x) \in \bold S^m(A)$
for some $A$ or just $p(\bar x)$ is directed by $\vdash$; 
i.e. for every finite $q(\bar x) \subseteq p(\bar x)$ there is $\psi(\bar
x,\bar b) \in p(\bar x)$ such that $\psi(\bar x,\bar x) 
\vdash q(\bar x)$,  divides over $B$ iff $\varphi(\bar x,\bar
b)$ explicitly divides over $B$, and they imply that in Definition
5.1(1), we can choose $\bar b_0 = \bar b$ and that $\{\varphi(\bar x,\bar b)\}$
forks over $B$.  If $p(\bar x) \in \bold S^m(A)$ or just $p(\bar x)$
is closed under conjunctions (or just is directed by $\vdash$,  
\ub{then} $p(\bar x)$ forks over $B$ iff some $\varphi(\bar
x,\bar a) \in p(\bar x)$ forks over $B$.  The $m$-type $p(\bar x)$
forks over $B$ iff there is $\varphi(\bar x,\bar a)$ which exactly
forks over $B$ such that $p(\bar x) \vdash \varphi(\bar x,\bar a)$.
\nl
3) [extension property] 
If an $m$-type $p$ is over $A$ and does not fork over $B$,
\ub{then} some extension $q \in \bold S(A)$ of $p$ does not fork
over $B$. 
\nl
4) [few non-forking types]
For $B \subseteq A$ the set $\{p \in \bold S^m(A):p$ does not fork over
$B$ (or just does not split strongly) over $B\}$ 
has cardinality $\le 2^{2^{|B|+|T|}}$. 

If $p(\bar x)$ does not fork
over $M$, \ub{then} it does not split over $M$.
\nl
5) [monotonicity in the sets]
If $B_1 \subseteq B_2 \subseteq A_2 \subseteq A_1$ and $p \in 
\bold S(A_1)$ does not fork/divide over 
$B_1$, \ub{then} $p \restriction A_2$ does not fork/divide over
$B_2$. 
\nl
6) [indiscernibility preservation] 
If $\bar{\bold b}$ is an infinite indiscernible sequence over $A_1$
and $B \subseteq A_1 \subseteq A_2$ and $\bar{\bold b} \subseteq A_2$
and tp$(\bar c,A_2)$ does not fork over $B$ or just does not divide
over $B$ or just does not split strongly over $B$ 
\ub{then} $\bar{\bold b}$ is an
(infinite) indiscernible sequence over $A_1 \cup \bar c$. 
\nl
7) [finite character]
If $p$ forks over $B$ \ub{iff} some finite $q \subseteq p$ does;
if $p$ is closed under conjunction (up to equivalence suffices)
\ub{then} we can demand $q = \{\varphi\}$.  Similarly for divides ...
\nl
8) [monotonicity in the type] 
If $p(\bar x) \subseteq q(\bar x)$ or just $q(\bar x) \vdash p(\bar x)$
and $p(\bar x)$ forks over $B$
\ub{then} $q(\bar x)$ forks over $B$; similarly for divides; for
split strongly this works only for $p(\bar x) \subseteq q(\bar x)$.
\nl
9) An $m$-type $p$ is finitely satisfiable in $A$ \ub{iff} for some
ultrafilter $D$ on ${}^m A$ we have 
$p \subseteq \text{ Av}(\text{Dom}(p),D)$. 
\bigskip

\remark{Remark}  1) Only parts (2), (4), (6) of \scite{10.3} use ``$T$ is
dependent".  \nl
2) If $T$ is unstable then for every $\kappa$ there are 
some $A$ and $p \in \bold S(A)$ such that $p$ 
divides over every $B \subseteq A$ of cardinality
$< \kappa$ (use a Dedekind cut with both cofinalities $\ge \kappa$).
\endremark
\bigskip

\demo{Proof}  0), 0A), 1) Easy. 

The proof of part (1) is included in the proof of part (2).
\nl
2) Assume toward contradiction that $p$ splits strongly, then for some
infinite indiscernible sequence $\langle \bar b_n:n < \omega \rangle$
over $B$ and $n < m$ we \footnote{recallng $[\varphi_1 \equiv
\varphi_2]$ is the formula $(\varphi_1 \wedge \varphi_2) \vee (\neg
\varphi_1 \wedge \neg \varphi_2)$}
have $[\varphi(\bar x,\bar b_n) \equiv \neg
\varphi(x,\bar b_m)] \in p$ (really $p \vdash [\varphi(\bar x,\bar b_n)
\equiv \neg \varphi(\bar x,\bar b_m)]$ suffices).  By renaming, \wilog
\, $n = 0,m=1$.  Let $\bar c_n = \bar b_{2n} \char 94 \bar
b_{2n+1},\psi(\bar x,\bar c_n) = [\varphi(\bar x,\bar b_{2n}) \equiv
\neg \varphi(\bar x,\bar b_{2n+1})]$.  Clearly $\langle \bar c_n:n <
\omega \rangle$ is an indiscernible sequence over $B,p \vdash
\psi(\bar x,\bar c_0)$ and $\{\psi(\bar x,\bar c_n):n < \omega\}$ is
contradictory as $T$ is dependent. 

This proves the first sentence.  The second is by the definitions and
the third sentence.  
For the third, the ``if" part is obvious, hence let us prove the
``only if", so assume that $p(\bar x)$
divides over $B$, we can find $\varphi(\bar x,\bar b_0),\langle \bar
b_n:n < \omega\rangle$ as in Definition 5.2(1), i.e. satisfies
clauses (a)-(d) there.  As $p(\bar x) \vdash \varphi(\bar x,\bar b_0)$,
necessarily there is a finite subset $p'(\bar x)$ of $p(\bar x)$ such
that $p'(\bar x) \vdash \varphi(\bar x,\bar b_0)$.  Let $\langle
\varphi_\ell(\bar x,\bar c_\ell):\ell < k\rangle$ list $p'(x)$ and ...
Now for each $n < \omega$, the sequences $\bar b_n,\bar b_0$
realize the same type over $B$, hence there is a sequence $\bar c^n
\in {}^{\ell g(\bar c)}{\frak C}$ such that the sequences $\bar b_0
\char 94 \bar c,\bar b_n \char 94 \bar c^n$ realize the same type
over $B$ and \wilog \, $\bar c^0 = \bar c$.  By Ramsey theorem and
compactness we can find $\langle \bar d_n:n < \omega\rangle$ such that
$\bar b_n \char 94 \bar d_n$ realizes the same type as $\bar b_0 \char
94 \bar c$ over $B$ and $\langle \bar b_n \char 94 \bar d_n:n <
\omega\rangle$ is an indiscernible sequence over $B$.  So let $F$ be
an automorphism of ${\frak C}$ over $B$ which maps $\bar b_0 \char 94
\bar d_0$ to $\bar b_0 \char 94 \bar c$.   So $\langle F(\bar d_n):n <
\omega\rangle$ is an indiscernible sequence over $B$ and $F(\bar d_0)
= \bar c$ so $\psi(\bar x,\bar d_0) = \psi(\bar x,\bar c) = \dsize
\bigwedge_{\ell < k} \varphi'_\ell(\bar x,\bar c) \equiv \dsize
\bigwedge_{\ell < k} \varphi_\ell(\bar x,\bar c_\ell) \vdash
\varphi(\bar x,\bar b_0) = \varphi(\bar x,F(\bar b_0))$.

Necessarily also $n < \omega \Rightarrow \psi(\bar x,\bar d_n) \vdash
\varphi(x,\bar b_n)$ and as $\{\varphi(\bar x,\bar b_n):n < \omega\}$
is contradictory, so is $\{\psi(\bar x,\bar d_n):n < \omega\}$.  So
$\langle F(\bar d_n):n < \omega\rangle$ examplifies that $\psi(\bar
x,\bar d_0) = \psi(\bar x,\bar c)$ explicitly divides over $B$ as
promised.

The fourth and fifth sentences are obvious.
\nl
3) By the definitions (or see \cite{Sh:93}). \nl
4) Easy or see \cite{Sh:3}; e.g. by part (3) \wilog \, $B=M,A=|N|$ is
$\|M\|^+$-saturated.  Now if $\bar a_\ell \in {}^m N$ realizes the
same type over $M$ for $\ell=1,2$ then for some $\bar c_n \in {}^m N$ for
$n=1,2,\dotsc,\langle \bar a_\ell \rangle \char 94 \langle \bar
c_1,\bar c_2,\ldots\rangle$ is indiscernible over $M$.
\nl
5) Easy. 
\nl
6) By part (2) and transitivity of ``equality of types" and Fact
\scite{10.3X} below.  \nl
7), 8), 9)  Easy.   \hfill$\square_{\scite{10.3}}$
\enddemo
\bn
We implicitly use the trivial:
\demo{\stag{10.3X} Fact}  1) If $I$ is a linear order, $\bar s_0,\bar
s_1$ are increasing $n$-tuples from $I$ \ub{then}
\mr
\item "{$\circledast_\omega$}"  there is a linear order $J \supseteq
I$ such that for $\ell \in \{0,1\}$ there is an indiscernible sequence
$\langle \bar t^\ell_k:k < \omega \rangle$ of increasing $n$-tuples from
$J$ such that $\bar t^0_{k+1} =\bar t^1_{k+1}$ for $k < \omega$ and
$\ell=0,1 \Rightarrow \bar s_\ell = \bar t^\ell_0$;
indiscernible means for quantifier free formulas in the order
language, i.e., in the vocabulary $\{<\}$ is satisfaction in $J$.  If
$I$ has no last element or no first element then we can take $I=J$.
\ermn
2) Similar for $\langle \bar b_t:t \in I \rangle$ an infinite indiscernible
sequence over $A$ in ${\frak C}$. 
\enddemo
\bigskip

\demo{Proof}  1) Let $J \supseteq I$ be with no last element.
Choose for $k=1,2,\ldots$ an increasing sequence $\bar t_k$ of length
$n$ from $J$ such that $2 \le k < \omega \Rightarrow$ 
Rang$(\bar s_0 \char 94 \bar s_1) < \text{ Rang}(\bar t_k) < 
\text{ Rang}(\bar t_{k+1})$.  So $\langle \bar s_\ell
\rangle \char 94 \langle \bar t_1,\bar t_2,\ldots \rangle$ is an
indiscernible sequence in $J$ for $\ell=0,1$.  \nl
2) Easy.    \hfill$\square_{\scite{10.3X}}$
\enddemo
\bigskip

\definition{\stag{10.1A} Definition}  1) Let $p$ be an $m$-type satisfying $p
\restriction B_2 \in \bold S^m(B_2)$.  We say that $p$ strictly does not
divide over $(B_1,B_2)$, (normally $B_1 \subseteq B_2$; 
when $B_1 = B_2 = B$ we may write ``over $B$") \ub{when}:
\mr
\item "{$(a)$}"  $p$ does not divide over $B_1$
\sn
\item "{$(b)$}"  if $\langle \bar c_n:n < \omega \rangle$ is an
indiscernible sequence over $B_2$ such that $\bar c_0$ realizes $p$
and $A$ is any set satisfying Dom$(p) \cup B_1 \subseteq A$, \ub{then} there 
is an indiscernible sequence $\langle \bar c'_n:n <
\omega \rangle$ over $A$ such that $\bar c'_0$ realizes $p$ and
tp$(\langle \bar c_n:n < \omega \rangle,B_2) = \text{ tp}(\langle \bar
c'_n:n < \omega \rangle,B_2)$.
\ermn
1A) ``Strictly divide" is the negation. 
\nl
2) We say that $p$ strictly forks over $(B_1,B_2)$ \ub{iff} $p \vdash
\dsize \bigvee_{\ell <n} \varphi_\ell$ for some $\langle
\varphi_\ell:\ell < n \rangle$ such that $(p \restriction B_2) \cup
\{\varphi_\ell\}$ strictly divides over $(B_1,B_2)$ for each $\ell < n$.
\nl
3) An $m$-type $p(\bar x)$ strictly does not fork over $(B_1,B_2)$
\ub{when}: $p(\bar x)$ does not fork over $B_1$, so $p(\bar x) \rest
B_2 \in \bold S^m(B_2)$ and if 
$\langle \bar c_n:n < \omega\rangle$ is an indiscernible
sequence over $B_2$ of sequences realizing $p(\bar x)$ and $C
\supseteq B_1 \cup \text{ Dom}(p(\bar x))$ and $q(\bar x) \in \bold
 S^m(C)$ extend $p(\bar x)$ and does not fork over $B_1$ \ub{then}
 there is an indiscernible sequence $\langle \bar c'_n:n <
 \omega\rangle$ over $C$ realizing tp$(\langle \bar c_n:n <
 \omega\rangle,B_2)$ such that $\bar c'_0$ realizes $q(\bar x)$;
note that ``strictly does not fork" is not defined as ``does not
strictly forks"; to stress we may write ``strictly$^*$ does not fork".
\enddefinition
\bn
We shall need some statements concerning ``strictly$^*$ does not fork"
parallel to those on ``does not fork".
\demo{\stag{10.3A} Observation}  0) In clause (b) of Definition
\scite{10.1A}(1) we can weaken the assumption ``$\bar c_0$ realizes $p$" to
``$\bar c_0$ realizes $p \restriction B_2$".
\nl
1) ``Strictly does not divide/fork
over $(B_1,B_2)$" is perserved by permuting the variables, repeating variables
and by automorphisms of
${\frak C}$ and if it holds for tp$(\bar b \char 94 \bar c,A)$ so $B_2
\subseteq A$ then it
holds for tp$(\bar b,A)$.  Similarly for ``does not strictly fork".
\nl
1A) The $m$-type $p(\bar x)$ strictly does not divide over $(B_1,B_2)$
\ub{iff} $p(\bar x) \rest B_2 \in \bold S^m(B_2)$ and $(p(\bar x)
\rest B_2) \cup q(\bar x)$ strictly does not divide over $(B_1,B_2)$
for every finite $q(\bar x) \subseteq p(\bar x)$.
\nl
2) If $p$ strictly$^*$ does not fork over $(B_1,B_2)$ then $p$ does not
   strictly fork over $(B_1,B_2)$ which \ub{implies} that $p$
strictly does not divide over $(B_1,B_2)$.  
\nl
3) If $p$ strictly does not divide over $(B_1,B_2)$ \ub{then} $p$ does not
divide over $B_1$. 
\nl
4) If $p$ does not strictly$^*$ fork over $B$ \ub{then} $p$ does not fork
over $B$. 
\nl
5) If $p$ is an $m$-type which strictly$^*$ does not fork over $(B_1,B_2)$
and $B_1 \cup \text{ Dom}(p) \subseteq A$ 
\ub{then} there is $q \in \bold S^m(A)$
extending $p$ which strictly$^*$ does not fork over
$(B_1,B_2)$.  Similarly does not strictly fork.
If $p_1(\bar x) \subseteq 
p_2(\bar x)$ and $p_1(\bar x)$ strictly does not fork over
$(B_1,B_2)$ and $p_2(\bar x)$ does not fork over $B_1$
 \underbar{then} $p_2(\bar x)$ strictly does not fork over $(B_1,B_2)$.
\nl
6) If $B_1 \subseteq B'_1 \subseteq B'_2 = B_2$ and $p(\bar x) \vdash
p'(\bar x)$ and $p(\bar x)$ does not strictly divide/fork over
$(B_1,B_2)$ and $p' \restriction B'_2$ is complete, i.e. $\in \bold
S^m(B'_2)$ \ub{then} $p'(\bar x)$ does not strictly divide/fork over 
$(B'_1,B'_2)$. 
\nl
7) In Definition \scite{10.1A}(1), clause (b) the case $A = \text{
Dom}(p) \cup B_2$ suffice. 
\nl
8) If $p$ strictly forks over $(B_1,B_2)$ \ub{then} for some finite $q
\subseteq p$ the type $q \cup (p \restriction B_2)$ strictly forks
over $(B_1,B_2)$.  Moreover, for some finite $B'_2 \subseteq B_2$,
($p$ is an $m$-type), if $B_1 \cup B'_2 \subseteq B''_2,p'$ is 
an $m$-type extending $q$ and $p'
\restriction B_2 \in \bold S^m(B_2)$ then $p'$ strictly forks over
$(B_1,B''_2)$.  Similarly for strictly divide.
\nl
9) If $M \subseteq A,p = \text{ tp}(\bar b,A)$ and tp$(A,M + \bar b)$ 
is finitely satisfiable in $M$, \ub{then} $p$ strictly does
not fork over $M$. 
\enddemo
\bigskip

\demo{Proof}  Easy, e.g., \nl
0) The new version is stronger hence it implies the one from the
definition.  

So assume that $p$ is an $m$-type, $p \restriction B_2 \in 
\bold S^m(B_2)$ and $p$ strictly does not divide over $(B_1,B_2)$ and we
shall prove the new version of clause (b).  I.e., we have $\langle
\bar c_n:n < \omega \rangle$ is an indiscernible sequence over $B_2$
and $\bar c_0$ realizes $p \restriction B_2$.  Let $\bar c''_0 \in
{}^m{\frak C}$ realizes $p$ hence it realizes $p \restriction B_2$,
but $p \restriction B_2 \in \bold S^m(B_2)$ so tp$(\bar c_0,B_2) =
\text{\rm tp}(\bar c''_0,B_2)$.  We can deduce that there is an
automorphism $F$ of ${\frak C}$ over $B_2$ which maps $\bar c_0$ to
$\bar c''_0$, and define $\bar c''_n = F(\bar c_n)$.

Now $\langle F(\bar c_n):n <\omega\rangle$ satisfies the assumption
of clause (b) from Definition \scite{10.1A}(1) hence there is 
an indiscernible sequence $\langle \bar c'_n:n < \omega \rangle$ over $A$ such
that tp$(\langle \bar c'_n:n < \omega \rangle,B_2) = \text{\rm
tp}(\langle F(\bar c'_n):n < \omega \rangle,B_2)$, but the latter is
equal to tp$(\langle c_n:n < \omega \rangle,B_2)$ so we are done.
\nl
1A) Let $p'(\bar x) := p(\bar x) \rest B_2$ so $p'(\bar x) \in \bold
S^m(B_2)$ and let $A = \text{ Dom}(p) \cup B_1$.  Now by parts (0),(1)
and the definitions
\mr
\item "{$\boxtimes_1$}"  $p(\bar x)$ strictly does not divide above
$B_2$ \ub{iff} the following holds
{\roster
\itemitem{ $(a)$ }  $p(\bar x)$ does not divide over $B_1$
\sn
\itemitem{ $(b)$ }  if $\bar{\bold c} = \langle \bar c_n:n <
\omega\rangle$ is an indiscernible sequence over $B_2$ then (where
$\ell g(\bar x_n) = m$) the set $\cup\{p(\bar x_n):n < \omega\} \cup
\Gamma^1_{\bar{\bold c}} \cup \Gamma^2_{\bar{\bold c}}$ is finitely
satisfiable where:
\sn
\itemitem{ ${{}}$ }  $(\alpha) \quad \Gamma^1_{\bar{\bold c}} =
\{\varphi(\bar x_0,\dotsc,\bar x_{n-1},\bar b):n < \omega,\varphi(\bar
x_0,\dotsc,\bar x_{n-1},\bar y) \in \Bbb L(\tau_1)$,
\nl

\hskip35pt $\bar b \subseteq
B_2$ and $\models \varphi[\bar c_0,\dotsc,\bar c_{n-1},\bar b]\}$
\sn
\itemitem{ ${{}}$ }  $(\beta) \quad \Gamma^2_{\bold c} =
\{\varphi(\bar x_{k_0},\dotsc,\bar x_{k_{n-1}},\bar a) \equiv 
\varphi(\bar x_{\ell_0},\dotsc,\bar x_{\ell_{n-1}},\bar a):n < \omega$
and
\nl

\hskip35pt  $k_0 < \ldots < k_{n-1}$ and $\ell_0 < \ldots < \ell_{n-1}$ and
$\bar a \subseteq A$ and
\nl

\hskip35pt  $\varphi(\bar x_0,\dotsc,\bar x_{n-1},\bar y)\in \Bbb L(\tau_T)\}$.
\endroster}
\ermn
But the demand in $\boxtimes$ has finite character.
\nl
2) For the first implication assume $p(\bar x)$ strictly$^*$ does not
fork over $(B_1,B_2)$ but $p(\bar x)$ strictly forks over
$(B_1,B_2)$ so $p(\bar x) \vdash \dsize \bigvee_{\ell < n}
\varphi_\ell(\bar x,\bar a_i)$ and $(p(\bar x) \rest B_2)  \cup
\{\varphi_\ell(\bar x,\bar a)$ strictly does not divide over
$(B_1,B_2)$ hence by part (1A) the type $p_\ell(\bar x) := p(\bar x)
\cup \{\varphi_\ell(\bar x,\bar a_\ell)\}$ strictly divides over
$(B_1,B_2)$ for each $\ell < m$.

But by Definition \scite{10.1A}(3) the type 
$p(\bar x)$ does not fork over $B_1$ hence we can choose $\ell < n$
such that $p_\ell(\bar x)$ does not fork over $B_1$.  As $p_\ell$
strictly divides over $(B_1,B_2)$ in Definition 5.6(1),
clause (a) or clause (b) there fail, for $p_\ell(\bar x)$, but clause
(a) holds by the choice of $\ell$, hence clause (b) fails and let
$\langle \bar c_n:n < \omega\rangle$ and $A$ exemplifies it; so \wilog
\, $A \supseteq \bar a_\ell \cup B_1 \cup B_2 \cup \text{ Dom}(p)$.  
Now let $q(\bar x) \in \bold S^m(A)$ extends
$p_\ell(\bar x)$ and does not fork over $B_1$ and apply Definition
5.6(3) with $p(\bar x),A,q(\bar x),\langle \bar c_n:n <
\omega\rangle$ here standing for $p(\bar x),C,q(\bar x),\langle \bar
c_n:n < \omega\rangle$ there; so we can find $\langle \bar c'_n:n <
\omega\rangle$ indiscernible over $A$ realizing tp$(\langle \bar c_n:n
< \omega\rangle,B_2)$ such that $\bar c'_0$ realizing $q(\bar x)$.
This contradicts the choice of $A,\langle \bar c_n:n <
\omega\rangle$.  The second implication holds by the earlier parts.
\nl
5) The third sentence holds by the definition.  The first sentence
   follows by the third and \scite{10.3}(3).  The second sentence
   follows by the definition.
\nl
6) Without loss of generality Dom$(p) \cup B'_2 \subseteq A$ and it
suffices to prove the case of ``strictly does not divide".
Recall that by part (1) in Claim \scite{10.3}(A), 
clause (b) we can demand only
``$\bar c_0$ realizes $p \restriction B_2$" and for any such $\langle
\bar c_n:n < \omega \rangle$ choose is $\bar c''_0$ realizing $p$ hence
$\bar c_0$ and $\bar c''_0$ realizes the same type over 
$B_2$ hence there is automorphism
$F$ of ${\frak C}$ over $B_2$ mapping $\bar c_0$ to $\bar c''_0$ and use the
definition for $\langle F(\bar c_n):n < \omega \rangle$.
\nl
\nl
7) By Ramsey theorem and compactness.
\nl
9) Use an ultrafilter $D$.  \hfill$\square_{\scite{10.3A}}$
\enddemo
\bn
The next claim is a parallel of: every type over $A$ does not fork over
some ``small" $B \subseteq A$.  If we have ``$p$ is over $A$ implies $p$ does
not fork over $A$" we could have improvement.

More elaborately, note that if $M$ is a dense linear order and $p \in
\bold S(M)$, then $p$ actually corresponds to a Dedekind cut of $M$.
So though in general $p$ is not definable, $p \restriction \{c \in M:c
\notin (a,b)\}$ is definable whenever $(a,b)$ is an interval of $M$
which includes the cut.  So $p$ is definable in large pieces.  The
following (as well as \scite{10.10}) realizes the hope that something in this
direction holds for every dependent theory.
\bigskip

\proclaim{\stag{10.2A} Claim}  If $p \in 
\bold S^m(A)$ and $B \subseteq A$, \ub{then} we can find
$C \subseteq A$ of cardinality $\le |T|$ such that:
\mr
\item "{$\circledast$}"  if $\langle \bar a_n:n < \omega \rangle$ is
an indiscernible sequence over $B \cup C$ such that $\bar a_0
\subseteq A$ and ${\text{\rm tp\/}}(\bar a_0,B \cup C)$ strictly$^*$ does 
not fork over $B$ and
$\{\varphi(\bar x,\bar a_n):n < \omega\}$ is contradictory or just
$\varphi(\bar x,\bar a_0)$ exactly forks 
over $B \cup C$, \ub{then} $\neg \varphi(x,\bar a_0) \in p$.
\endroster
\endproclaim
\bigskip

\remark{Remark}  The first case in $\circledast$, (the one on
``$\{\varphi(\bar x,\bar a_n):n <\omega\}$ is contradictory, says,
in other words no formula in $p$ divides over $B \cup C$
when the type of the sequence of parameters over $B \cup C$ does not
fork over $B$.
\endremark
\bigskip

\demo{\stag{10.2B} Conclusion}  1) For every $p \in \bold S^m(A)$ 
and $B \subseteq A$,  we can find 
$C \subseteq A,|C| \le |T|$ such that:
\mr
\item "{$\circledast$}"   if $\langle \bar a_n:n
< \omega \rangle$ is an indiscernible sequence over $B \cup C$ satisfying
$\bar a_0 \cup \bar a_1 \subseteq A$ and ${\text{\rm tp\/}}(\bar a_0 
\char 94 \bar a_1,B \cup C)$ strictly$^*$
does not fork over $B$, \ub{then} for any $\varphi$
{\roster
\itemitem{ $(*)$ }  $\varphi(x,\bar a_0) \in 
p \text{ iff } \varphi(x,\bar a_1) \in p$.
\endroster}
\ermn
2) For every $\bar x = \langle x_\ell:\ell < m \rangle$ and formula
$\varphi = \varphi(\bar x;\bar y)$ 
for some finite $\Delta \subseteq \Bbb L(T)$ we have:
\block
if $p \in \bold S^m(A),B \subseteq A$, \ub{then} for some finite $C
\subseteq A$ (in fact $|C| < f(m,\varphi,T)$ for some function $f$), we
have:
\endblock
\sn
\block
if $\langle \bar a_\ell:\ell < k \rangle$ is
$\Delta$-indiscernible sequence over $B \cup C$ and tp$_\Delta(\bar
a_0 \char 94 \bar a_1,B \cup C)$ does not fork
over $A$, \ub{then} $\varphi(x,\bar a_0) \in p \Leftrightarrow
\varphi(x,\bar a_1) \in p$.
\endblock
3) The local version of \scite{10.2A} holds with a priori finite bound on $C$.
\enddemo
\bigskip

\demo{Proof of \scite{10.2A}}  By induction on 
$\alpha < |T|^+$ we try to choose $C_\alpha,\bar a_\alpha,k_\alpha$ and
$\langle \bar a^k_{\alpha,n}:n < \omega \rangle$ and
$\varphi_\alpha(\bar x,\bar y_\alpha),\varphi_{\alpha,k}
(\bar x,\bar y_{\alpha,k})$ such that:
\mr
\item "{$\boxtimes$}"   $(a) \quad C_\alpha = 
\cup \{\bar a_\beta:\beta < \alpha\} \cup B$
\sn
\item "{${{}}$}"   $(b) \quad \langle 
\bar a^k_{\alpha,n}:n < \omega \rangle$ is an
indiscernible sequence over $C_\alpha$ for $k < k_\alpha$
\sn
\item "{${{}}$}"   $(c) \quad \bar a_\alpha \subseteq A$ and $\bar
a_\alpha = \bar a^k_{\alpha,0}$ for $k < k_\alpha$
\sn
\item "{${{}}$}"   $(d) \quad \varphi_\alpha(\bar x,\bar a_\alpha) \in p$
\sn
\item "{${{}}$}"  $(e) \quad 
\{\varphi_{\alpha,k}(\bar x,\bar a^k_{\alpha,n}):n <
\omega\}$ is contradictory
\sn
\item "{${{}}$}"  $(f) \quad$ 
tp$(\bar a_\alpha,B \cup C_\alpha)$ strictly$^*$ does not fork
over $B$
\sn
\item "{${{}}$}"  $(g) \quad\varphi_\alpha(\bar x,\bar a_\alpha) \vdash \dsize
\bigvee_{k < k_\alpha} \varphi_{\alpha,k}(\bar x,\bar a^k_{\alpha,0})$.
\ermn
If for some $\alpha < |T|^+$ we are stuck, $C = C_\alpha \backslash B$ is as
required.  So assume that we have carried the induction and we shall
eventually get a contradiction.

By induction on $\alpha < |T|^+$ we choose $D_\alpha,F_\alpha,\bar b_\beta,
\langle
\bar b^k_{\beta,n}:n < \omega \rangle$ for $\beta < \alpha$ such that
(but $\bar b^k_{\alpha,n}$ are defined in the $(\alpha +1)$-th stage):
\mr
\item "{$\boxplus$}"   $(\alpha) \quad F_\alpha$ is an 
elementary mapping, increasing continuous with $\alpha$
\sn
\item "{${{}}$}"   $(\beta) \quad$ Dom$(F_\alpha) = C_\alpha$ and 
Rang$(F_\alpha) \subseteq D_\alpha$
\sn
\item "{${{}}$}"   $(\gamma) \quad D_\alpha = 
\text{ Rang}(F_\alpha) \cup \bigcup 
\{\bar b^k_{\beta,n}:\beta < \alpha,k < k_\alpha$ and 
$n < \omega\}$ so $D_\alpha \subseteq {\frak C}$ is increasing continuous
\sn
\item "{${{}}$}"   $(\delta) \quad
\langle \bar b^k_{\alpha,n}:n < \omega \rangle$
is an indiscernible sequence over $D_\alpha \supseteq
F_\alpha(C_\alpha)$ and $\bar b^k_{\alpha,0} = \bar b_\alpha$
\sn
\item "{${{}}$}"   $(\varepsilon) \quad F_{\alpha +1}(\bar a_\alpha) = \bar
b_\alpha$ and tp$(\bar b_\alpha,D_\alpha)$ does not fork over $F_\alpha(B)$
\sn
\item "{${{}}$}"   $(\zeta) \quad$ some automorphism 
$F^k_{\alpha +1} \supseteq F_{\alpha +1}$ of 
${\frak C}$ maps $\bar a^k_{\alpha,n}$ to $\bar b^k_{\alpha,n}$ for 
$n < \omega,k < k_\alpha$.
\ermn
For $\alpha=0,\alpha$ limit this is trivial.  For $\alpha = \beta +1$,
clearly $F_\alpha(\text{tp}(\bar a_\alpha,C_\alpha))$ 
is a type in $\bold S^{< \omega}(F_\alpha(C_\alpha))$
which strictly$^*$ 
does not fork over $F_\alpha(B) = F_0(B)$ hence by \scite{10.3}(3) has an
extension $q_\alpha \in \bold S^{< \omega}(D_\alpha)$ which does not
fork over $F_0(B)$ and let $\bar b_\alpha$ realize it.  
Let $F_{\alpha +1} \supseteq F_\alpha +1$ be
the elementary mapping extending $F_\alpha$ with domain $C_{\alpha
+1}$ mapping $\bar a_\alpha$ to $\bar b_\alpha$.  Let $F^k_{\alpha +1}
\supseteq F_\alpha$ be an automorphism of ${\frak C}$ 
as required by clauses $(\delta) + (\zeta)$; $F^k_{\alpha +1}$
exists as tp$(\bar a_\alpha,B \cup C_\alpha)$ strictly$^*$ does not fork
over $B$ and let $\bar b^k_{\alpha,n} = F^k_{\alpha +1}
(\bar a^k_{\alpha,n})$ for 
$n < \omega,k < k_\alpha$.  So $D_{\alpha +1}$ and
$F_{\alpha +1}$ are well defined.

Having carried the induction let $F \supseteq \bigcup\{F_\alpha:\alpha <
|T|^+\}$ be an automorphism of ${\frak C}$.  We claim that for each
$\alpha < |T|^+$ and $k < k_\alpha$, for every $\beta \in 
[\alpha,|T|^+]$ we have
\mr
\item "{$(*)_{\alpha,\beta}$}"  $\langle \bar b^k_{\alpha,n}:
n < \omega \rangle$ is an indiscernible sequence over $D_\alpha \cup 
\bigcup \{\bar b_\gamma:\gamma \in [\alpha +1,\beta)\}$.
\ermn
We prove this by induction on $\beta$.  For $\beta = \alpha$ this
holds by clause $(\delta)$, for $\beta \equiv \alpha +1$ this is the
same as for $\beta = \alpha$.  For 
$\beta$ limit use the definition of indiscernibility.  For 
$\beta = \zeta +1 > \alpha +1$ use tp$(\bar b_\zeta,D_\zeta)$ does not
fork over $F_0(B)$ hence over $D_\alpha \cup \{F_{\gamma
+1}(\bar b_\gamma):\alpha < \gamma <
\zeta\}$ by \scite{10.3}(5); so by the induction hypothesis and
\scite{10.3}(6) clearly $(*)_{\alpha,\beta}$ holds.  

\relax From $\alpha < \beta \le |T|^+ \Rightarrow (*)_{\alpha,\beta}$ we can conclude
\mr
\item "{$(**)$}"  for any $n < \omega$ and $\alpha_0 < \ldots <
\alpha_{n-1} < |T|^+$ and $\nu \in \dsize \prod_{\ell < n} k_{\alpha_\ell}$
and $\eta \in {}^n 2$ the sequences $\bar b^{\nu(0)}_{\alpha_0,0} 
\char 94 \bar b^{\nu(1)}_{\alpha_1,0} \char 94 \ldots 
\char 94 \bar b^{\nu(n-1)}_{\alpha_{n-1},0}$ and 
$\bar b^{\nu(0)}_{\alpha_0,\eta(0)} \char 94
\bar b^{\nu(1)}_{\alpha_1,\eta(1)} \char 94 \ldots \char 94 \bar
b^{\nu(n-1)}_{\alpha_{n-1},\eta(n-1)}$ realize the same type over $B$. \nl
[Why?  By induction on $\ell(*) = \text{ max}\{\ell:\eta(\ell) = 1$ 
or $\ell = -1\}$.
If $\ell(*) = -1$ then the two sequences are the same so $(**)$ holds
trivially.  Let $\rho \in {}^n 2$ be defined by $\rho(\ell)$ is 0 if
$\ell \ne \ell(*)$ and $\eta(\ell)$ otherwise, so the induction
hypothesis apply to $\rho$ so it suffices to prove that $\bar
b^\nu_\eta,\bar b^\nu_\rho$ realizes the same type over $B$.  Now  
assume $\ell(*) \in \{0,\dotsc,n-1\}$ and use
$(*)_{\alpha_{\ell(*)},|T|^+}$ for $k = \nu(\ell_*)$, it says that the
sequence $\langle \bar b^{\nu(\ell(*))}_{\alpha_{\ell(*),n}}:n <
\omega\rangle$ is an indiscernible sequence over $D_{\alpha_{\ell(*)}} \cup
\{\bar b_\gamma:\gamma \in [\alpha_{\ell(*)+1},|T|^+)\}$.

The second part in the union includes

$$
\{\bar b^{\nu(\ell(*)-1)}_{\alpha_{\ell(*)+1},0},\dotsc,\bar
b^{\nu(n-1)}_{\alpha_{n-1},0}\} = \{\bar b^{\nu(\ell(*))}
_{\alpha_{\ell(*)+1},\eta(\ell(*)+1)},\dotsc,
\bar b^{\nu(n-1)}_{\alpha_{n-1},\eta(n-1)}\}$$

$$
\hskip50pt = \{\bar
b^{\nu(\ell(*))}_{\alpha_{\ell(*)+1}},\ell(*)+1,\ldots,\bar
b^{\nu(n-1)}_{\alpha_{n-1},\rho(n-1)}\} 
$$
\mn
by the choice of $\ell(*)$,
and the first part of the union includes the rest so it suffices to
have that $\bar b_{\alpha_{\ell(*)}} = \bar
b^{\nu(\ell(*))}_{\alpha_{\ell(*)},0} =
b^{\nu(\ell(*))}_{\alpha_{\ell(*)},\rho(\ell(*))}$ and $\pm \bar
b^{\nu(\ell(*))}_{\alpha_{\ell(*)},\eta(\ell(*))} =
b^{\nu(\ell(*))}_{\alpha_{\ell(*)},1}$ realizes the same type
$D_{\alpha_{\ell(*)}} \cup \{\bar b_\gamma:\gamma =
\alpha_{\ell(*)+1},\dotsc,\alpha_{n-1}$, which holds.]
\ermn
Let $\bar c$ realize $F(p)$.  For each $\alpha <
|T|^+,\varphi_\alpha(x,\bar a_\alpha) \in p$ hence
$\varphi_\alpha(x,\bar b_\alpha) \in F(p)$.  Also
$\varphi_\alpha(x,\bar a_\alpha) \vdash \dsize \bigvee_{k<k_\alpha}
\varphi_{\alpha,k}(x,a^k_{\alpha,0})$ hence by clause $(\zeta)$ we have
$\varphi_\alpha(x,\bar b_\alpha) \vdash \dsize \bigvee_{k<k_\alpha}
\varphi_{\alpha,k}(x,a^k_{\alpha,0})$ hence we can choose $k(\alpha) <
k_\alpha$ such that ${\frak C} \models 
\varphi[\bar c,\bar a^{k(\alpha)}_{\alpha,0}]$.

Now as $\{\varphi_{\alpha,k(\alpha)}(x,\bar b^{k(\alpha)}_{\alpha,n}):
n < \omega\}$ is
contradictory there is $n = n[\alpha] < \omega$ such that ${\frak C} \models
\neg \varphi_{\alpha,k(\alpha)}(\bar c,\bar
b^{k(\alpha)}_{\alpha,n})$ whereas ${\frak C} \models
\varphi_{\alpha,k(\alpha)}[\bar c,\bar b^{k(\alpha)}_{\alpha,0}]$; 
by renaming without loss of generality ${\frak C} \models \neg
\varphi_{\alpha,k(\alpha)}[c,\bar b^{k(\alpha)}_{\alpha,n}]$ for 
$\alpha < |T|^+,n \in
[1,\omega)$.  Now if $n < \omega,\alpha_0 < \ldots < \alpha_{n-1} <
|T|^+$ and $\eta \in {}^n 2$ then ${\frak C} \models \dsize
\bigwedge_{\ell < m} \varphi_{\alpha_\ell,k(\alpha_\ell)}(\bar c,
\bar b^{k(\alpha_\ell)}_{\alpha_\ell,\eta(\ell)})^{\text{if}(\eta(\ell)=0)}$ 
hence ${\frak C} \models
(\exists \bar x)[\dsize \bigwedge_{\ell < n}
\varphi_{\alpha_\ell,k(\alpha_\ell)}
(\bar x,\bar b_{\alpha_\ell,\eta(\ell)})^{\text{if}(\eta(\ell)=0)}]$
hence by $(**)$ we have 
\nl
${\frak C} \models (\exists \bar x)[\dsize \bigwedge_{\ell < n}
\varphi_{\alpha_\ell,k(\alpha_\ell)}
(\bar x,\bar
b^{k(\alpha_\ell)}_{\alpha_\ell,0})^{\text{if}(\eta(\ell)=0)}]$.  
\nl
Hence the independence property holds, contradiction. 
\hfill$\square_{\scite{10.2A}}$
\enddemo
\bigskip

\demo{Proof of \scite{10.2B}}  1) Follows from \scite{10.2A} by
\scite{10.3}(2). 
\nl
2) By \scite{10.2A} and compactness or repeating the proof. 
\nl
3) Similarly.    \hfill$\square_{\scite{10.2B}}$
\enddemo
\bigskip

\proclaim{\stag{10.4} Claim}  1) Assume $p$ is a type, $B \subseteq M,
{\text{\rm Dom\/}}(p) \subseteq M$ and 
$M$ is $|B|^+$-saturated.  \ub{Then}
\mr
\item "{$(A)$}"   $p$ does not fork
over $B$ \ub{iff} $p$ has a complete extension over $M$ which does not
fork over $B$ \ub{iff} $p$ has a complete extension over $M$ which
does not divide over $B$ \ub{iff} $p$ has a complete extension over
$M$ which does not split strongly over $B$
\sn
\item "{$(B)$}"   if $p = { \text{\rm tp\/}}(\bar c,M)$ and
$\varphi(\bar x,\bar a) \in p$ forks over $B$, \ub{then} for some 
$\langle \bar a_n:n < \omega
\rangle$ indiscernible over $B,\{\bar a_n:n < \omega\} \subseteq M,\bar
a_0 = a$ and $\neg \varphi(\bar x,\bar a_1) \in p$, and of course,
$\varphi(\bar x,\bar a_0) \in p$
\sn
\item "{$(C)$}"   if $p \in \bold S^m(M)$ \ub{then} $p$ forks over $B$
iff $p$ exactly forks over $B$.
\ermn
2) Assume ${\text{\rm tp\/}}(C_1/A)$ does not fork over 
$B \subseteq A$ and ${\text{\rm tp\/}}(C_2,(A \cup C_1))$ does not 
fork over $B \cup C_1$.  \ub{Then} ${\text{\rm tp\/}}(C_1 \cup C_2),A)$
does not fork over $B$.
\endproclaim
\bigskip

\demo{Proof}  1) Read the definitions. 
\sn
\ub{Clause $(A)$}: \nl
First implies second by \scite{10.3}(3), second implies third by
Definition \scite{10.1} or \scite{10.3}(2), 
third implies fourth by \scite{10.3}(2).  If the first fails, then 
$p \vdash \dsize \bigvee_{\ell < k} \varphi_\ell(\bar x,\bar a_\ell)$
for some $k$ where each $\varphi_\ell(\bar x,\bar a_\ell)$ divides over $B$;
let $\langle \bar a_{\ell,n}:n < \omega \rangle$ witness this hence by
\scite{10.3}(2) \wilog \,
$\bar a_\ell = \bar a_{\ell,0}$.  As $M$ is $|B|^+$-saturated, \wilog
\,  $\bar a_{\ell,n} \subseteq M$.   So for every 
$q \in \bold S^m(M)$ extending $p$, for some $\ell <
k,\varphi_\ell(\bar x,\bar a_\ell) \in q$ but for every large enough
$n,\neg \varphi_\ell(\bar x,\bar a_{\ell,n}) \in q$, so $q$ splits
strongly, i.e., fourth fails.  So fourth implies first, ``closing the circle".
\sn
\ub{Clause $(B)$}: \nl
Similar. 
\nl
2) Let $M$ be $|B|^+$-saturated model such that $A \subseteq M$.  By
\scite{10.3}(3) there is an elementary mapping $f_1$ such that $f_1
\restriction A = \text{ id}_A$ and Dom$(f_1) = C_1 \cup A$ 
and $f_1(C_1)/M$ does not fork over $B$.  Similarly
we can find an elementary mapping $f \supseteq f_1$ such that Dom$(f)
= C_1 \cup C_2 \cup A$ and $f(C_2)/(M \cup f(C_1))$ does not fork over $A
\cup f(C_1)$.  By \scite{10.3}(2), $f_1(C_1)/M$ does not split
strongly over $B$.  Again by \scite{10.3}(2), $f(C_2)/(M \cup
f_1(C_1))$ does not split strongly over $B \cup f_1(C_1)$.  Together
they imply that
if $\bar{\bold b} \subseteq M$ is an infinite indiscernible sequence
over $B$ then it is an indiscernible sequence over  $f(C_1) \cup B$ and
even over $f(C_2) \cup (f(C_1) \cup B)$ (use the two previous
sentences and \scite{10.3}(6)).  
But this means that $f(C_1) \cup f(C_2)/M$ does not split
strongly over $B$, (here the exact version of strong splitting we
choose is immaterial as $M$ is $|B|^+$-saturated).  So by
\scite{10.4}(1) we get that $f(C_1) \cup f(C_2)/M$ does not fork over
$B$ hence $f(C_1 \cup C_2)/A$ does not fork over $B$ but $f \supseteq
\text{ id}_A$ so by \scite{10.3}(2) 
also $C_1 \cup C_2/A$ does not split strongly over $B$.
\hfill$\square_{\scite{10.4}}$
\enddemo
\bigskip

\demo{\stag{10.4.3} Conclusion}  1) If $M$ is 
$|B|^+$-saturated and $B \subseteq M$ and $p \in \bold S^n(M)$ 
\ub{then} $p$ does not fork over $B$ iff $p$ does not
strongly split over $B$.
\nl
2) If $A = |M|$, then in Conclusion \scite{10.2B}(1) we can replace 
strong splitting by dividing.
\enddemo
\bigskip

\demo{Proof}  1) By \scite{10.4}(1A).
\nl
2) By part (1).  \hfill$\square_{\scite{10.4.3}}$
\enddemo
\bigskip

\definition{\stag{10.5} Definition}  1) We say $\langle \bar a_t:t \in
J \rangle$ is a non-forking sequence over $(B,A)$ \ub{when} $B \subseteq A$
and for every $t \in J$ the type tp$(\bar a_t,A \cup \bigcup \{\bar a_s:s
<_J t\})$ does not fork over $B$.  \nl
2) We say that $\langle \bar a_t:t \in J \rangle$ is a strict
non-forking sequence over $(B_1,B_2,A)$ \ub{if} $B_1 \subseteq B_2
\subseteq A$ and for every $t \in J$ the type tp$(\bar a_t,A \cup
\bigcup\{\bar a_s:s <_J t\})$ strictly$^*$ does not fork over $(B_1,B_2)$,
see Definition \scite{10.1A}(3). 
If $B_1 = B_2$ we may write $(B_1,A)$ instead of $(B_1,B_1,A)$.
\nl
3) We say ${\Cal A} = (A,\langle (\bar a_\alpha,B_\alpha):\alpha <
\alpha^* \rangle)$ is an $\bold F^f_\kappa$-construction or $\langle
(\bar a_\alpha,B_\alpha):\alpha < \alpha^* \rangle$ an $\bold
F^f_\kappa$-construction over $A$ \ub{if} $B_\alpha \subseteq A_\alpha := A
\cup \bigcup \{\bar a_\beta:\beta < \alpha\}$ has cardinality $< \kappa$
and tp$(\bar a_\alpha,A_\alpha)$ does not fork over $B_\alpha$.
\nl
4) We can above replace $\bar a_t$ by $A_t$ meaning for some/every
$\bar a_t$ listing $A_t$ the demand holds.
\enddefinition
\bigskip

\proclaim{\stag{10.6} Claim}  
1) Assume
\mr
\item "{$(a)$}"  $\langle \bar a_t:t \in J \rangle$ is a strict
non-forking sequence over $(B,B,A)$
\sn
\item "{$(b)$}"  $\langle \bar b^\varepsilon_{t,n}:n < \omega \rangle$ is an
indiscernible sequence over $A$, each $\bar b_{t,n}$ realizing 
${\text{\rm tp\/}}(\bar a_t,A)$ for each $t \in J,\varepsilon < \zeta_t$.
\ermn
\ub{Then} we can find $\bar a_{t,n}$ for $t \in J,n <
\omega,\varepsilon < \zeta_t$ such that
\mr
\item "{$(\alpha)$}"  $\langle \bar a_{t,n}:n < \omega \rangle$ is an
indiscernible sequence over $A_1 \cup \{\bar a_{s,n}:n < \omega,s \in J
\backslash \{t\}\}$
\sn
\item "{$(\beta)$}"  $\text{\rm tp}(\langle \bar a_{t,n}:
n < \omega \rangle,A) = \text{\rm tp}(\langle \bar b_{t,n}:
n < \omega \rangle,A)$
\sn
\item "{$(\gamma)$}"  $\bar a_{t,0} = \bar a_t$.
\endroster
\endproclaim
\bigskip

\demo{Proof}  We prove by induction on $|J|$, both parts.
\sn 
\ub{Case 1}:   $J$ is finite.  

We prove this by induction on $n = |J|$, for $n=0,1$ this is trivial; assume
we have proved for $n$ and we shall prove for $n+1$.  Let $\lambda =
(|A| + |T|)^+$.

So let $J = \{t_\ell:\ell \le n\}$ with $t_\ell$ increasing with
$\ell$.  First we can find an
indiscernible sequence $\langle \bar c^\varepsilon_{t_0,\alpha}:
\alpha < \lambda
\rangle$ over $A$ such that $\bar c^\varepsilon_{t_0,0} = \bar a_{t_0}$ and for
some automorphism $F_\varepsilon$ of ${\frak C}$ over $B$ we 
have $k < \omega \Rightarrow F(\bar b^\varepsilon_{t_0,k}) 
= \bar c^\varepsilon_{t_0,k}$ and let $A' := A \cup 
\{\bar c^\varepsilon_{t_0,\alpha}:\alpha < \lambda\}$.  
[This is possible by Definition \scite{10.1A}.] 
\nl
Second, we can choose $\bar a'_{t_\ell}$ by induction on $\ell$ such
that $\bar a'_{t_0} = \bar a_{t_0}$ and if $\ell >0$ then tp$(\bar
a'_{t_\ell},A' \cup \bigcup \{\bar a'_{t_m}:m=1,\dotsc,\ell-1\})$ 
strictly$^*$ does not fork over $B$ and the two sequences
 $\bar a_{t_0} \char 94 \ldots \char 94 \bar a_{t_\ell},
\bar a'_{t_0} \char 94 \ldots \char 94 \bar a'_{t_\ell}$
realizes the same type over $A$.  We can do it by \scite{10.3A}(5) and
``strictly$^*$ does not fork" being preserved by 
elementary mapping.  By \scite{10.4}(2)
the type tp$(\bar a'_{t_1} \char 94 
\ldots \char 94 \bar a'_{t_n},A'\})$ does not
fork over $B$ hence by \scite{10.3}(6) the sequence $\langle \bar
a_{t_0,\alpha}:\alpha < \lambda \rangle$ is an indiscernible sequence
over $A \cup (\bar a'_{t_1} \char 94 \ldots \char 94 \bar a'_{t_n})$.

As tp$(\bar a_{t_\ell},A \cup \{\bar a_{t_m}:m < i\})$ strictly does
not fork over $(B,A)$ without loss of generality 
$\langle \bar b_{t_\ell,m}:m <
\omega\rangle$ is an indiscernible sequence over $A'$ such that each
member realizes tp$(\bar a_{t_\ell},A')$.

Now we use the induction hypothesis with $B,A',\langle \bar
a'_{t_\ell}:\ell =1,\dotsc,n \rangle,\langle \bar b_{t_\ell,m}:m < \omega
\rangle$ for $\ell = 1,\dotsc,n$ and let $\langle \bar a'_{t_\ell,n}:n <
\omega \rangle$ for $\ell =1,\dotsc,n$ be as in the claim. 

By \cite{Sh:715} for some $\alpha^* < \lambda$ the sequence $\langle \bar
c^\varepsilon_{t_0,\alpha}:\alpha \in [\alpha^*,
\alpha^* + \omega) \rangle$ is an
indiscernible sequence over $A \cup \bigcup \{\bar a'_{t_\ell,m}:m <
\omega,\ell=1,\dotsc,n\}$ and as $A' = A \cup \{\bar
a'_{t_0,\alpha}:\alpha < \lambda\}$ clearly for $\ell = 1,\dotsc,n$
the sequence $\langle \bar a'_{t_\ell,m}:m < \omega \rangle$ is
indiscernible over $A \cup \bigcup \{\bar a'_{t_k,m}:k \in \{1,\dotsc,n\}
\backslash \{\ell\}$ and $m < \omega\} \cup 
\bigcup \{\bar a'_{\alpha^* +m}:m < \omega\}$.
But we know that $\langle \bar
a'_{t_0,\alpha}:\alpha< \alpha^* + \omega \rangle$ is an indiscernible
sequence over $A \cup \{\bar a'_{t_\ell}:\ell =1,\dotsc,n\}$, hence
the sequence
$\bar a'_{t_\alpha,\alpha^*} \char 94 \bar a'_{t_1} \char 94 \ldots
\char 94 \bar a'_{t_n}$ realizes over $A$ the same type as $\bar
a'_{t_0,0} \char 94 \bar a'_{t_1} \char 94 \ldots \char 94 \bar
a'_{t_n}$ hence realizes over $A$ the same type 
as $\bar a_{t_0} \char 94 \bar a_{t_1} \char 94 \ldots
\char 94 \bar a_{t_n}$.  So for some automorphism $F$ of ${\frak C},F
\restriction A = \text{ id}_A,\ell =1,\dotsc,n \Rightarrow \bar
a_{t_\ell} = F(\bar a'_{t_\ell,0})$ and $\bar a_{t_0} = F(\bar
a'_{t_0,\alpha^*})$ and let $\bar a_{t_\ell,m} = F(\bar
a'_{t_\ell,m})$ for $\ell =1,\dotsc,n$ and $m < \omega$ and $\bar
a_{t_0,m} = F(\bar a'_{t_0,\alpha^* +m})$. \nl
So we are done.
\enddemo
\bn
\ub{Case 2}:  $J$ infinite.

By Case 1 + compactness.  \hfill$\square_{\scite{10.6}}$
\bigskip

\remark{Remark}  Can we use just no dividing?
\endremark
\bigskip

\proclaim{\stag{10.7} Claim}  1) Assume $\langle A_t:t \in J \rangle$
is a non-forking sequence over $(B,A)$ and $C_t \subseteq {\frak C}$
for $t \in J$.  \ub{Then} we can find $\langle f_t:t \in J \rangle$
such that
\mr
\item "{$(a)$}" $f_t$ is an elementary mapping with domain
$$ A \cup A_t \cup C_t
$$
\sn
\item "{$(b)$}"  $f_t \restriction (A \cup A_t)$ is the identity
\sn
\item "{$(c)$}"  ${\text{\rm tp\/}}(A_t,A \cup 
\bigcup \{A_s \cup f_s(C_s):s < t\}$ does
not fork over $B$.
\ermn
2) If in addition ${\text{\rm tp\/}}(C_t,A \cup A_t)$ does 
not fork over $A \cup A_t$ \ub{then} we can add
\mr
\item "{$(c)^+$}"  $\langle A_t \cup f_t(C_t):t \in J \rangle$ is a
non-forking sequence over $(B,A)$.
\endroster
\endproclaim
\bigskip

\remark{\stag{10.7Y} Remark}  1) We may consider $\bold
F^f$-construction, i.e., ${\Cal A} = (A,\langle a^{B_i}_\alpha:\alpha
  < \alpha^*\rangle)$ is an $\bold F^f$-construction, when
\mr
\item "{$(a)$}"  $B_i \subseteq A_i := A \cup \{a_j:j<i\}$
\sn
\item "{$(b)$}"  tp$(a_i,A_i)$ does not fork over $B_i$
\sn
\item "{$(c)$}"  $|B_i| < \kappa$.
\ermn
1A) We may replace above $\alpha$ by a linear order $I$, not
necessarily well founded.
\nl
2) In \scite{10.7}(2) we may weaken the assumption to: for every $A'
\supseteq A,A_t \cup C_t/A$ can be embedded to a complete
non-forking type over $A'$.
\endremark
\bigskip

\demo{Proof}  1) As in the proof of \scite{10.2A}. \nl
2) Similarly.
\enddemo
\bigskip

\proclaim{\stag{10.7A} Claim}  1) Assume
\mr
\item "{$(a)$}"  $\langle A_t:t \in J \rangle$ is a non-forking
sequence over $(B,A)$.  
\ermn
\ub{Then} for any initial segment $I$ of $J,{\text{\rm tp\/}}
(\bigcup\{A_t:t \in J \backslash I\},A \cup \bigcup \{A_t:t \in I\})$ does not
fork over $B$.
\nl
2) Assume (a) and
\mr
\item "{$(b)$}"  $\langle \bar a_{t,n}:n < \omega \rangle$ is an
indiscernible sequence over $A$,
\sn
\item "{$(c)$}"  $\bar a_{t,n} \in {}^{\omega >} (A_t)$
\sn
\item "{$(d)$}"  $\langle \bar a_{t,n}:n < \omega \rangle$ is an
indiscernible sequence over $A \cup \bigcup\{A_s:s <_J t\}$.
\ermn
\ub{Then} $\langle \langle \bar a_{t,n}:n < \omega \rangle:t \in J \rangle$
are mutually indiscernible over $A$.  Also for any non-zero $k <
\omega$ and $t_0 < \ldots < t_{k-1}$ in $J$ the sequences $\langle
\bar a_{t_\ell,n}:n < \omega \rangle$ for all $\ell < k$ are mutually
indiscernible over $A \cup \bigcup \{A_s:\neg(t_0 \le s \le t_{k-1})\}$.
\endproclaim
\bn
\margintag{10.7B}\ub{\stag{10.7B} Question}:  If $n_\ell < \omega$ for $\ell < n$ 
do the sequences
$\langle \bar a_{t_0,n_0} \char 94 \bar a_{t_1,n_1} \char 94 \ldots
\char 94 \bar a_{t_{k-1},n_{k-1}} \rangle$ and $\langle \bar a_{t_0,0} \char
94 \bar a_{t_1,0} \char 94 \ldots \char 94 \bar a_{t_{k-1},0}\rangle$
realize the same type over $A \cup \bigcup \{A_s:s <_J t_0$ or $s_J >
t_{k-1}\}$.  Need less?
\bigskip

\remark{Remark}  A statement similar to \scite{10.7A}(1) for $\bold
F^f_\kappa$-construction holds.
\endremark
\bigskip

\demo{Proof}  1) If $J \backslash I$ is finite, we prove this 
by induction on $|J \backslash I|$ using \scite{10.4}(2).  
The general case follows by \scite{10.3}(7).
\nl
2) It is enough to prove the second sentence.
For $k=1$ this follows by \scite{10.3}(6) + \scite{10.4}(2) 
using part (1) with $A
\cup \bigcup\{A_s:s < t\},\langle A_r:r \in J,r > t \rangle$ instead
$A,\langle A_r:r \in J \rangle$.

For $k+1 >1$, let us be given $t_0 <_J \ldots <_J t_k$.
Use the case $k=1$ for each $t_\ell$ and combine.  
\hfill$\square_{\scite{10.7A}}$
\enddemo

\remark{\stag{10.7A.2} Remark}  1) Recall that by \scite{10.4} if 
$p \in \bold S^m(M),M$ is 
quite saturated, \ub{then} dividing is the same as forking for the
type $p$.
\endremark
\bigskip

\proclaim{\stag{10.9} Claim}   Assume that for every 
set $B$, if $p(\bar x) \in \bold S^m(B)$ then
$p$ does not fork over $B$.

Assume that $\langle \bar a_t:t \in J \rangle$ is
a non-forking sequence over $(B,A)$ and $A = |M|$. 
\nl
1) For every (finite sequence) $\bar b$ the set 
$\{t:\bar b/(A \cup \bar a_t)$ forks over $\dbcu_{s<t} \bar a_s \cup A$ 
has cardinality $\le |T|$.  
\nl
2) For each $\varphi(\bar x,\bar y,\bar z)$ and $k < \omega$ for some $n =
n_{\varphi(\bar x,\bar y),k}$ the set $W^\varphi_{\bar b} :=
\{t:{\text{\rm tp\/}}_\varphi(\bar b,A \cup \bar a_t)$ has a subset
with $\le k$ members which forks over $\dbcu_{s<t} \bar a_s \cup A\}$ has
$\le n$ members.
\endproclaim
\bigskip

\demo{Proof}  1) By (2). 
\nl
2) Fix $k$.  Assume toward contradiction that this fails for $n$.  We
can find $t_0 <_I t_1 <_J \ldots <_J t_{n-1}$ from $W^\varphi_{\bar b}$.

Now for every $u \subseteq \{0,\dotsc,n-1\}$ there is
$\bar b_u$ realizing tp$(\bar b,A \cup\{\bar a_{t_\ell}:\ell \in u\})$
such that tp$(\bar b_u,A \cup \{\bar a_{t_\ell}:\ell < n\})$ does
not fork over $A \cup \{\bar b_{t_\ell}:\ell \in u\}$.  Now for each
$\ell < n$ we can find $q_\ell \subseteq \text{ tp}_\varphi(\bar b,A
\cup \bar a_{t_\ell})$ and with $\le k$ members which forks over
$\dbcu_{s<t} \bar a_s \cup A$; let $A_\ell = A \cup \bar a_{t_0} \cup
\ldots \cup \bar a_{t_\ell-1}$.
Clearly $\ell \in u \Rightarrow q_\ell \subseteq \text{ tp}(\bar b_u,A
\cup \{\bar a_{t_m}:m <n\})$.  Now if $\ell \in n \backslash u$, let
$i_{\ell,0} < \ldots < i_{\ell,m(\ell)-1} < n$ list $u \backslash
\ell$ so tp$(\bar a_{i_{\ell,m}},A_\ell \cup \bar a_{t_\ell} \cup \bar
a_{t_{i_{\ell,0}}} \ldots \cup \bar a_{t_{i_{\ell,m(\ell)-1}}})$ does not fork
over $A$ for $m < m(\ell)$ and tp$(\bar b_u,A \cup \bar a_{t_\ell}
\cup \bar a_{i_\ell} \cup \ldots \cup \bar a_{t_{i_{\ell,m(\ell)-1}}})$ does
not fork over $A \cup \{\bar b_{t_k}:k \in u\} \subseteq A_\ell \cup
\{\bar a_{t_{i_{\ell,0}}} \cup \ldots \cup \bar
a_{t_{i_{\ell,m(\ell)-1}}})$ 
hence by \scite{10.4}(2) + \scite{10.3}(0) the type tp$(\bar
b_u,A_\ell \cup \bar a_{t_\ell})$ does not fork over $A_\ell$.  Hence
now as tp$_\varphi(\bar b,A \cup \bar a_{t_\ell})$ has the subset
$q_\ell$ with $\le k$
members which forks over $\dbcu_{s<t} \bar a_s \cup A$, by
monotonicity $q_\ell$ it forks also over its subset $A_\ell$, hence $q_\ell$
forks over $A_\ell$ hence by the previous sentence $q_\ell \nsubseteq
\text{ tp}(b_u,A \cup \bar b_{t_\ell})$ so $\neg \wedge q_\ell \in
\text{ tp}(b_u,A \cup \bar b_{t_\ell})$.  
As for our fixed $k$ this holds for
every $n$, we get that $T$ has the indpendence property contradiction.
\hfill$\square_{\scite{10.9}}$
\enddemo
\bigskip

\proclaim{\stag{10.10} Claim}  Assume that $p(\bar x)$ is a type of cardinality
$< \kappa$ which does not fork over $A$.  \ub{Then} for some $B
\subseteq A$ of cardinality $< \kappa + |T|^+$, the type $p$ does not
fork over $B$.
\endproclaim
\bigskip

\demo{Proof}  Without loss of generality $p$ is closed under
conjunction.

For any finite sequence $\bar \varphi = \langle(\varphi_\ell(\bar
x,\bar y_\ell),n_\ell):\ell < n \rangle$ and formula 
$\psi(x,\bar c) \in p$ and set $B \subseteq A$ we define

$$
\align
\Gamma_{B,\bar \varphi,\psi(\bar x,\bar c)} = \{&(\forall
x)(\psi(x,\bar c) \rightarrow \dsize \bigvee_{\ell < n}
\varphi_\ell(\bar x,\bar y_{\ell,0}))\} \cup \\
  &\{\neg(\exists \bar x) \dsize \bigwedge_{n \in w} \varphi_\ell(\bar
x,y_{\ell,n}):\ell < n \text{ and } w \in [\omega]^{n_\ell}\} \cup \\
  &\{\vartheta(y_{\ell,m_1},\dotsc,y_{\ell,m_k},\bar b) =
\vartheta(y_{\ell,0},\dotsc,y_{\ell,k},\bar b): \\
  &\qquad \bar b \subseteq B,\vartheta \in \Bbb L(\tau_T),m_1 < \ldots
< m_k < \omega\}.
\endalign
$$
\mn
Now as $p$ does not fork over $A$, clearly for any $\bar \varphi$ as
above and
$\psi(\bar x,\bar c) \in p$ the set $\Gamma_{A,\bar \varphi,\psi(\bar
x,\bar c)}$ is inconsistent.  Hence for some finite set $B = B_{\bar
\varphi,\psi(x,\bar c)} \subseteq A$ the set $\Gamma_{B,\bar
\varphi,\psi(x,\bar c)}$ is inconsistent.  Now $B^* =
\cup\{B_{\bar \varphi,\psi(\bar x,\bar c)}:\psi(\bar x,\bar c) \in p$
and $\bar \varphi$ is as above$\}$ 
is as required.  \hfill$\square_{\scite{10.10}}$
\enddemo
\bn
The following is another substitute for ``every type $p$ does not fork
over a small subset of Dom$(p)$".
\proclaim{\stag{10.11A} Claim}  Assume that for 
every set $B$, if $p \in \bold S^{< \omega}(B)$ then
$p$ does not fork over $B$.

Assume $p \in \bold S^m(M)$ and $B
\subseteq M$.  \ub{Then} we can find $C$ such that
\mr
\item "{$(*)_1$}"  $C \subseteq M$ and $|C| \le |T|$ and
\sn
\item "{$(*)^p_{M,B,C}$}"  if $D \subseteq M$ and ${\text{\rm tp\/}}
(D/B \cup C)$
does not fork over $B$ \ub{then} $p \restriction (B \cup D)$ does not fork
over $B \cup C$.
\endroster
\endproclaim
\bigskip

\demo{Proof}  Follows by \scite{10.9}.
\enddemo
\bn
\centerline {$* \qquad * \qquad *$}
\bigskip

\definition{\stag{ns.1} Definition}  Assume that $C = |M|,M$ is
$\kappa$-saturated $A \subseteq M,|A| < \kappa$ and $p \in \bold
S^m(M)$ does not split over $A$.  For any set $B(\subseteq {\frak C})$
let $p^{[B,A]}$ be $q \restriction B$ where $q \in \bold S^m(M \cup
B)$ is the unique type in $\bold S^m(M \cup B)$ which does not split
over $A$.
\enddefinition
\bigskip

\demo{\stag{ns.2} Observation}  1) In Definition \scite{ns.1},
$p^{[B,A]}$ is well defined.
\nl
2) In \scite{ns.1} instead ``$C$ is $\kappa^+$-saturated; $|A| <
\kappa$" it suffices to assume that every $q \in 
\bold S^{<\omega}(B)$ is realized in $C$ that is $C$ is full over $B$.
\nl
3) $p^{[A,B_1]} \subseteq p^{[A,B_2]}$ if $B_1 \subseteq B_2$.
\enddemo
\bigskip

\proclaim{\stag{ns.3} Claim}  1) If the triple $(A,C,p)$ is as in
\scite{ns.2}(2), $A \subseteq A_0$ and $\bar a_n \in {}^m {\frak C}$
realizes $p^{[A_n,C]}$ for $n < \omega$ where $A_n= A_0 \cup
\bigcup\{\bar a_\ell:\ell < n\}$ \ub{then} $\langle \bar a_n:n <
\omega\rangle$ is an indiscernible sequence over $A_0$.  Also 
{\rm tp}$(\langle \bar a_n:n < \omega \rangle,A_0)$ is 
determined by $(A,C,p,A_0)$ and as $p$ determines $C$, 
we call it $p^{[A,A_0,\omega]}$.
\endproclaim
\bigskip

\demo{Proof}  See \cite[II,\S1]{Sh:c} or \cite{Sh:3}.
\enddemo
\bigskip

\proclaim{\stag{ns.4} Claim}  Assume that
\mr
\item "{$(a)$}"  $C \supseteq A$ is full over $A$
\sn
\item "{$(b)$}"  $p_0,p_1 \in \bold S^m(C)$ does not split over $A$
\sn
\item "{$(c)$}"  $p^{[A,A,\omega]}_0 = p^{[A,A,\omega]}_1$.
\ermn
\ub{Then} $p_0 = p_1$.
\endproclaim
\bigskip

\demo{Proof}  Let $\langle \bar a^\ell_n:n <\omega\rangle$ realize
$p^{[A,A,\omega]}_\ell$ so by clause (c) of the assumption
\mr
\item "{$(*)_1$}"  $\bar a^0_0,\dotsc,\bar a^0_{n-1}$ and $\bar a^1_0
\char 94 \ldots \char 94 \bar a^1_{n-1}$ realizes the same type over
$A$.
\ermn
If the conclusion fails, we can find $\bar c$ and 
$\varphi(\bar x,\bar y) \in \Bbb L(\tau_T)$ such that
\mr
\item "{$(*)_2$}"  $\neg \varphi(\bar x,\bar c) \in p_0$
and $\varphi(\bar x,\bar c) \in p_1$ so $\bar c \in
{}^{\ell g(\bar y)} C$.
\ermn
Now we choose by induction on $n < \omega$ a sequence $\bar a_n$ such that
\mr
\item "{$(*)_3$}"  if $\ell < 2$ and $n=\ell$ mod $2$ and we let
$A_n = \cup \bigcup\{\bar a_0,\dotsc,\bar a_{n-1}\}$ \ub{then}
tp$(\bar a_n,A_n \cup \bar c) = p^{[A,A_{k,n} \cup \bar c]}_\ell$.
\ermn
Now we can prove by induction on $n < \omega$ that
\mr
\item "{$(*)_4$}"  the sequences $\bar a^0_0 \char 94 \ldots \char 94
a^0_{n-1},\bar a^1_0 \char 94 \ldots \char 94 \bar a^1_{n-1}$ and
$\bar a_0 \char 94 \ldots \char 94 \bar a_{n-1}$ realizes the same
type over $A$.
\ermn
[Why?  The first two sequences realizes the same type by $(*)_1$.  For
the induction step, if $n=\ell$ mod $2$, by the definition \scite{ns.1},
we have $\bar a^\ell_0 \char 94 \ldots \char 94 \bar a^\ell_{n-1} \char 94
\bar a^\ell_n$ and $\bar a_0 \char 94 \ldots \char 94 \bar a_{n-1}
\char 94 \bar a_n$ realizes the same type over $A$.]

So $\langle \bar a_n:n < \omega\rangle$ is an indiscernible sequence
and ${\frak C} \models \varphi[\bar a_n,\bar c]$ iff $n$ is odd,
contradiction to ``$T$ is dependent".
\hfill$\square_{\scite{ns.4}}$
\enddemo
\bigskip

\demo{\stag{ns.5} Conclusion}  1) If $A \subseteq C$ and every $p \in
\bold S^{< \omega}(A)$ is realized in $C$ then 
$\{p \in \bold S^m(C):p$ does not split over $A\}$ has
cardinality $\le |\bold S^\omega(A)|$ which is 
$\le (\text{Ded}_r(|A|+|T|)^{|T|} \le 2^{|A|+|T|}$ recalling
Ded$_r(\mu) = \text{ Min}\{\lambda:\lambda$ is regular and every
linear order of density $\le \mu$ has cardinality $\le \lambda\}$.
\nl
2) Also for any finite $\Delta  \subseteq \Bbb L(\tau_T)$, the set
$\{p \restriction \Delta:p \in \bold S^m(C)$ does not split over $A\}$
has cardinality $\le \text{ Ded}_r(C)$.
\nl
3) If $p \in \bold S^m(C)$ is finitely satisfiable in $A \subseteq C$
then $p$ does not split over $A$.
\enddemo
\bigskip

\demo{Proof}  Should be clear.
\enddemo
\bn
\centerline{$* \qquad * \qquad *$}
\bigskip

\definition{\stag{10.31} Definition}  For $\ell \in \{1,2\}$, we 
say $\{\bar a_\alpha:\alpha <
\alpha^*\}$ is $\ell$-independent over $A$ if: we can find $\bar
a_{\alpha,n}$ (for $\alpha < \alpha^*,n < \omega \rangle$ such that:
\mr
\item "{$(a)$}"  $\bar a_\alpha = \bar a_{\alpha,0}$
\sn
\item "{$(b)$}"  $\langle \bar a_{\alpha,n}:n < \omega \rangle$ is an
indiscernible sequence over $A \cup \bigcup\{\bar a_{\beta,m}:\beta \in
\alpha^* \backslash \{\alpha\}$ and $m < \omega\}$
\sn
\item "{$(c)$}"  $(\alpha) \quad$ if 
$\ell =1$ then for some $\bar b_n \in A \,\,(n
< \omega)$ for every $\alpha < \alpha^*$ we 

\hskip25pt  have
$\langle \bar b_n:n < \omega \rangle \char 94
\langle \bar a_{\alpha,n}:n < \omega \rangle$ is an indiscernible
sequence
\sn
\item "{${{}}$}"  $(\beta) \quad$ if $\ell = 2$ then for some $\bar
b_{\alpha,n} \subseteq A$ (for $\alpha < \alpha^*,n < \omega$),
\nl

\hskip25pt  $\langle \bar b_{\alpha,n}:n < \omega \rangle \char 94
\langle \bar a_{\alpha,n}:n < \omega \rangle$ is an indiscernible
sequence.
\endroster
\enddefinition
\bn
We now show that even a very weak version of independence has
limitations.
\proclaim{\stag{10.32} Claim}  1) For every finite $\Delta \subseteq
\Bbb L(\tau_T)$ there is $n^* < \omega$ such that we cannot find $\bar
\varphi = \langle \varphi_n(\bar x,\bar a_n):n < n^* \rangle$ such
that
\mr
\item "{$(*)_{\bar \varphi}$}"  for each $n < n^*$ there are $m_n < \omega$
and $\langle \bar b^n_{m,\ell}:\ell < \omega,m < m_n \rangle$ and
$\langle \psi^n_m(\bar x,\bar y_n):m < m_n \rangle$ such that
{\roster
\itemitem{ $(\alpha)$ }  $\langle \bar b^n_{m,\ell}:\ell < \omega
\rangle$ is an indiscernible sequence over $\cup\{\bar a_k:k < n^*,k
\ne n\}$
\sn   
\itemitem{ $(\beta)$ }   $\bar b^n_{m,0} = \bar a_n$
\sn   
\itemitem{ $(\gamma)$ }  $\{\psi^n_m(x,\bar b^n_{m,\ell}):\ell <
\omega\}$ is contradictory for each $n$ and $m < m_n$
\sn   
\itemitem{ $(\delta)$ }   $\psi^n_m(\bar x,\bar y_n) \in \Delta$
\sn  
\itemitem{ $(\varepsilon)$ }    $\varphi_n(\bar x,\bar a_n) \vdash
\dsize \bigvee_{m < m_n} \psi^n_m(\bar x,\bar a_n)$
\sn  
\itemitem{ $(\zeta)$ }   $\models(\exists \bar x) \dsize
\bigwedge_{n<n^*} \varphi_n(\bar x,\bar a_n)$.
\endroster}
\ermn
2) We weaken $(\alpha)$ above to ${\text{\rm tp\/}}(\bar b^n_{m,\ell},
\bigcup\{\bar a_k:k < n^*,k \ne n\}) = { \text{\rm tp\/}}
(\bar a_n,\bigcup\{\bar a_k:k < n^*,k \ne n\})$.
\nl
3) Above for some finite $\Delta^+ \subseteq \Bbb L(\tau_T)$, we can
in $(\alpha)$ demand only $\Delta^+$-indiscernible; also \wilog \,\,
$\varphi_n(\bar x,\bar y_n) = \dsize \bigvee_{m < m_n} \psi^n_m(\bar
x,\bar y_n)$. 
\endproclaim
\bigskip

\demo{Proof}  1) [Close to \scite{10.2A}.]   Note
\mr
\item "{$\circledast$}"  if $\bar c \in {}^{\ell g(\bar x)}({\frak
C})$ and $n < n^*$ and $\models \varphi_n(\bar c,\bar a_n)$ then for some
$\bar c' \in {}^{\ell g(\bar x)}({\frak C})$ we have
{\roster
\itemitem{ $(i)$ }  tp$(\bar c',\bigcup\{\bar a_k:k < n^*,k \ne n\}) =
\text{ tp}(\bar c,\bigcup\{\bar a_k:k < n^*,k \ne n\})$
\sn
\itemitem{ $(ii)$ }  tp$_\Delta(\bar c,\bar a_n) \ne \text{
tp}_\Delta(\bar c',\bar a_n)$.
\endroster}
\ermn
[Why $\circledast$ holds?  Clearly it is enough to find $\bar b'_n$
such that
\mr
\item "{$(i)$}"  $\bar b_n,\bar b'_n$ realize the same type over
$\cup\{\bar a_k:k < n^*,k \ne n\}$
\sn
\item "{$(ii)$}"  for some $m < m_n$ we have $\psi^n_m(\bar b_n,
\bar a_n) \wedge \neg \psi^n_m(\bar b'_n,\bar a_n)$.
\ermn
Why does $\bar b'_n$ exist?  As $\models \varphi_n[\bar c,\bar a_n]$
by $(\varepsilon)$ for some $m < m_n,\models \psi^n_n[\bar c,\bar
a_n]$ and by $(\alpha) + (\gamma)$, for some $\ell < \omega,b'_n =
\bar b^n_{m,\ell}$ is as required.]

By repeated use of $\circledast$ we get $m^*_\ell < m_\ell$ such that
$\langle \psi^n_{m^*_\ell}(\bar x,\bar a_n):n < n^* \rangle$ is
independent but $\psi^n_{m^*_\ell}(\bar x,\bar y_n) \in \Delta$ is
finite, so $n^*$ as required exists.
\nl
2),3) Similarly.  \hfill$\square_{\scite{10.32}}$
\enddemo
\bigskip

\proclaim{\stag{10.33} Claim}  Assume
\mr
\item "{$(a)$}"  $\langle \bar b_n:n < \omega \rangle$ is
indiscernible over $M$
\sn
\item "{$(b)$}"  $\{\varphi(\bar x,\bar b_n):n < \omega\}$ is
contradictory
\sn
\item "{$(c)$}"  $M \prec N,p \in \bold S(N),\varphi(\bar x,\bar b_0)
\in p$ and $\neg \varphi(x,\bar b_n) \in p$ for $n > 0$
\sn
\item "{$(d)$}"  $N$ is $\|M\|^+$-saturated.
\ermn
\ub{Then} for some $\langle \bar b'_n:n < \omega \rangle$ we have
\mr
\item "{$(\alpha)$}"   $\langle \bar b'_n:n < \omega \rangle$ is
indiscernible over $M$ based on $M$, $\bar b'_n \subseteq N$
\sn
\item "{$(\beta)$}"  $\bar b'_0 \in \{\bar b_0,\bar b_1\}$
\sn
\item "{$(\gamma)$}"   $\varphi(\bar x,\bar b'_0) \equiv \neg
\varphi(\bar x,\bar b'_1)$ belongs to $p$.
\endroster
\endproclaim
\bigskip

\demo{Proof}  Easy.
\enddemo
\bigskip

\definition{\stag{10.34} Definition}  1) For $p \in \bold S^m(M)$ let
${\Cal E}(p)$ be the set of pairs $(\varphi(\bar x,\bar y),\bold e)$
such that
\mr
\item "{$(a)$}"  $\bold e$ is a definable equivalence relation on
${}^{\ell g(\bar y)}M$ in $M$
\sn
\item "{$(b)$}"  if $\bar b_1 \bold e \bar b_2$ then $\varphi(\bar
x,\bar b_1) \in p \Leftrightarrow \varphi(\bar x,\bar b_2) \in p$. 
\ermn
2) ${\Cal E}'_{\text{tp}}(p)$ is defined similarly by $\bold e$ is
definable by types. 
\enddefinition
\bigskip

\proclaim{\stag{10.35} Claim}  Assume $\varphi = \varphi(x,\bar y),M \prec
N,N$ is $\|M\|^+$-saturated and $p \in \bold S(N)$.  Then we cannot
find $\{D_i:i < n_\varphi\}$, a set of ultrafilters over 
${}^{\ell g(\bar y)}(N)$ pairwise orthogonal (see below)
with $p_i = { \text{\rm Av\/}}(N,D_i)$ such that $p(x) \cup p_i(\bar
y_0) \cup p_i(\bar y_1) \cup \{\varphi(\bar x,\bar y_1),\neg
\varphi(x,\bar y_0)\}$ is consistent for $i < n_\varphi$.
\endproclaim
\bn
\ub{Now we deal with Orthogonality}.
\definition{\stag{or.1} Definition}  
1) Two complete types $p(\bar x),q(\bar y)$ over $A$ are weakly
orthogonal if $p(\bar x) \cup q(\bar y)$ is a complete type over $A$. 
\nl
2) Assume $\bar{\bold b}_1,\bar{\bold
b}_2$ are endless indiscernible sequences.  We say $\bar{\bold
b}_1,\bar{\bold b}_2$ are orthogonal and write $\bar{\bold b}_1 \perp
\bar{\bold b}_2$ if:

\block
for every set $A$ which includes $\bar{\bold b}_1 \cup \bar{\bold
b}_2$, Av$(A,\bar{\bold b}_1)$, Av$(A,\bar{\bold b}_2)$ are weakly
orthogonal
\endblock
3) $\bar{\bold b}_1$ is strongly orthogonal to $\bar{\bold b}_2,
\bold b_1 \underset{\text{st}} {}\to \perp \bold b_2$ \ub{if} it
is orthogonal to every endless indiscernible sequence 
$\bar{\bold b}'_2$ of finite distance from 
$\bar{\bold b}_2$ (see \cite[1.11=np1.4B]{Sh:715}(2). \nl
4) An endless indiscernible sequence $\bar{\bold b}_1$ is orthogonal
to $\varphi(x,\bar a)$ \ub{if} it is orthogonal to every endless
indiscernible sequence $\bar{\bold b}_2 = \langle b_{2,\alpha}:\alpha
< \delta\rangle$ such that $b_{2,\alpha} \in \varphi({\frak C},\bar
a)$ for every $\alpha < \delta$.   Similarly for $\varphi_1(x,\bar
a_1),\varphi_2(x,\bar a_2)$.  If tp$(b_\alpha,A)$, tp$(\cup\{\bar
b_\beta:\beta \in \alpha^* \backslash \{\alpha\},)$ for $\alpha <
\alpha^*$ \underbar{then} every $p \in \bold S^m(A)$ is weakly
orthogonal to tp$(\bar b_\alpha,A)$ for all but $\le |T|$ of the
ordinal $\alpha < \alpha^+$.
\nl
5) $\bar{\bold b}$ is based on $A$ if $\bar{\bold b}$ is an
indiscernible sequence and $C_A(\bar{\bold b})$ (see 
\cite{Sh:715} or \cite{Sh:93}) has boundedly many conjugations over $A$. \nl
6) If $\bar{\bold b}_1 \underset{\text{st}} {}\to \perp 
\bar{\bold b}_2$ and $\bar{\bold b}'_\ell$
is a neighbor (see \cite[1.11=np1.4B]{Sh:715}) to 
$\bar{\bold b}_\ell$ \ub{then} $\bar{\bold b}'_1$ is strongly
orthogonal to $\bar{\bold b}'_2$. 
\enddefinition
\bigskip

\proclaim{\stag{or.2} Claim}  1) Orthogonality is symmetric
relation. \nl
2) If $\bold b_1,\bold b_2$ are orthogonal \ub{then} they are
perpendicular (see Definition \scite{np8.1.8}). 
\endproclaim
\bn
\margintag{or.2D}\ub{\stag{or.2D} Example}:  In Th$(\Bbb R,<)$, different initial segments are
orthogonal, even two disjoint intervals.  For $(\Bbb R,0,1,+,\times)$ the
situation is different: any two non trivial intervals are ``the same".
\bigskip

\proclaim{\stag{or.3} Claim}  1) 
Assume $\lambda = \lambda^{< \lambda},I$ is a dense
linear order with neither first nor last element and $\bar{\bold b} = \langle
\bar b_t:t \in I \rangle$ an indiscernible sequence.  If $|I| =
\lambda$, \ub{then} there is $M \supseteq \bar{\bold b}$ which is
$\lambda$-saturated and $\lambda$-atomic over $\bar{\bold b}$. \nl
2) If $p \in \bold S^m(\bar{\bold b})$ is $\lambda$-isolated
\ub{then} it is $|T|^+$-isolated. 
\endproclaim
\bn
\margintag{or.3D}\ub{\stag{or.3D} Question}:  $M$ is $\mu$-minimal 
over $\bar{\bold b}$, i.e. over $B :=
 \cup\{\bar b_t:t \in I\}$ (i.e. there is
no $N,\bar{\bold b} \subseteq N \prec M,N \ne M$ such that 
$N$ is $\mu^+$-saturated) \ub{iff} there is no indiscernible sequence
 over $\bar{\bold b}$ which is $|T|^+$-isolated over $\bar{\bold b}$.

\bn
\margintag{or.4}\ub{\stag{or.4} Question}:  If Av$(M,\bar{\bold b}_1)$, Av$(M,\bar{\bold b}_2)$
(or with $D$'s) are weakly orthogonal and are perpendicular, \ub{then}
are they orthogonal?
\bn
\margintag{or.5}\ub{\stag{or.5} Question}:  On the set of elements or sequences
realizing the type Av$(\bar{\bold b}_1,\bar{\bold b}),\bar{\bold b}$
an endless indiscernible sequence, can we define a dependence relation
similar enough to the stable case (so if $\bar{\bold b} \char 94
\bar{\bold b}_2$ is indiscernible then $\bar{\bold b}_1$ is
independent).
\bn
\margintag{or.6}\ub{\stag{or.6} Question}:  For each of the following conditions can we
characterize the dependent theories which satisfy it?
\mr
\item "{$(a)$}"  for any two non-trivial indiscernible sequences
$\bar{\bold b}_1,\bar{\bold b}_2$, we can find $\bar{\bold b}'_\ell$
of finite distance from $\bar{\bold b}_\ell$ (see \cite{Sh:715}, for
$\ell=1,2$) such that $\bar{\bold b}'_1,\bar{\bold b}'_2$ are not
orthogonal
\sn 
\item "{$(b)$}"  any two non-trivial indiscernible sequences of
singletons have finite distance?
\sn
\item "{$(c)$}"  $T$ is Th$(\Bbb F),\Bbb F$ a field (so this class
includes the $p$-adics various reasonable fields with valuations and
closed under finite extensions).
\endroster

\newpage
    
REFERENCES.  
\bibliographystyle{lit-plain}
\bibliography{lista,listb,listx,listf,liste}

\enddocument